\documentclass[reqno,a4paper]{siamltex}

\usepackage{graphicx}
\usepackage{color}

\long\def\drop#1{}
\usepackage{amsmath,amsfonts}
\usepackage{geompsfi}
\usepackage{mathrsfs}
\usepackage{a4wide}
\ifx\theorem\undefined 
\newtheorem{theorem}{Theorem}

\DeclareMathOperator\supp{supp}
\fi
\DeclareMathOperator\tr{tr}

\newcommand{\figref}[1]{Figure~\ref{#1}}

\newcounter{remark}[section]
\def\theremark{\thesection.\arabic{remark}.}
\newenvironment{remark}%
    {\par\medbreak\refstepcounter{remark}%
         {\noindent\bf Remark~\theremark\ }}%
    {\par\medbreak}

\def\ds{\displaystyle}

\def\MP{\mathrm{MP}}
\def\wMP{w_\MP^{}}

\long\def\meta#1{\texttt{#1}}
\def\pref#1{(\text{\ref{#1}})}
\def\mod#1{\left|#1\right|}
\def\Mod#1{\left\|#1\right\|}

\def\Z{\mathbb Z}
\def\R{\mathbb R}
\let\e\varepsilon
\let\d\delta
\def\blambda{\bar\lambda}
\def\limp{P_{\mathrm{imp}}}
\def\lambdacr{P_{\mathrm{cr}}}
\def\l{\lambda}
\def\bw{\overline w}
\def\bphi{\overline \phi}
\def\vKD{Von K\'arm\'an-Donnell}%
\def\Xint#1{\mathchoice
   {\XXint\displaystyle\textstyle{#1}}%
   {\XXint\textstyle\scriptstyle{#1}}%
   {\XXint\scriptstyle\scriptscriptstyle{#1}}%
   {\XXint\scriptscriptstyle\scriptscriptstyle{#1}}%
   \!\int}
\def\XXint#1#2#3{{\setbox0=\hbox{$#1{#2#3}{\int}$}
     \vcenter{\hbox{$#2#3$}}\kern-.5\wd0}}

\def\dashint{\Xint-}
\def\frparbcenter#1{\framebox(43,60)[t]{\parbox{43mm}{\begin{center}#1\end{center}}}}

\author{Ji\v r\'I Hor\'ak\thanks{Universit\"at K\"oln, Germany}
\and Gabriel J. Lord\thanks{Heriot-Watt University, Edinburgh, United Kingdom}
\and Mark A. Peletier\thanks{Technische Universiteit Eindhoven, The Netherlands}}
\title{Cylinder Buckling: The Mountain Pass as an Organizing Center}
\begin{document}
\maketitle

\begin{abstract}
We revisit the classical problem of the buckling of a long thin axially
compressed cylindrical shell. By examining the energy landscape of the
perfect cylinder we deduce an estimate of the sensitivity of the shell
to imperfections. Key to obtaining this is the existence of a mountain
pass point for the system.
We prove the existence on bounded domains of such solutions for all
most all loads and 
then numerically compute example 
mountain pass solutions. Numerically the mountain pass solution with
lowest energy has the form of a single dimple.
We interpret these results and validate
the lower bound against some experimental results available in the literature.
\end{abstract}


\pagestyle{myheadings}
\thispagestyle{plain}
\markboth{JI\v R\'I HOR\'AK, GABRIEL J. LORD, AND MARK A. PELETIER}{CYLINDER BUCKLING}

\section{Introduction}
\subsection{Buckling of cylinders under axial loading}

A classical problem in structural engineering is the prediction of the
load-carrying capacity of an axially-loaded cylinder.
As well as being a commonly used structural element, the axially-loaded
cylinder is also the archetype of unstable, imperfection-sensitive
buckling, and this has led to a large body of
theoretical and experimental research.

\begin{figure}[htb]
\centering
\centerline{$\vcenter{\psfig{width=8cm,figure=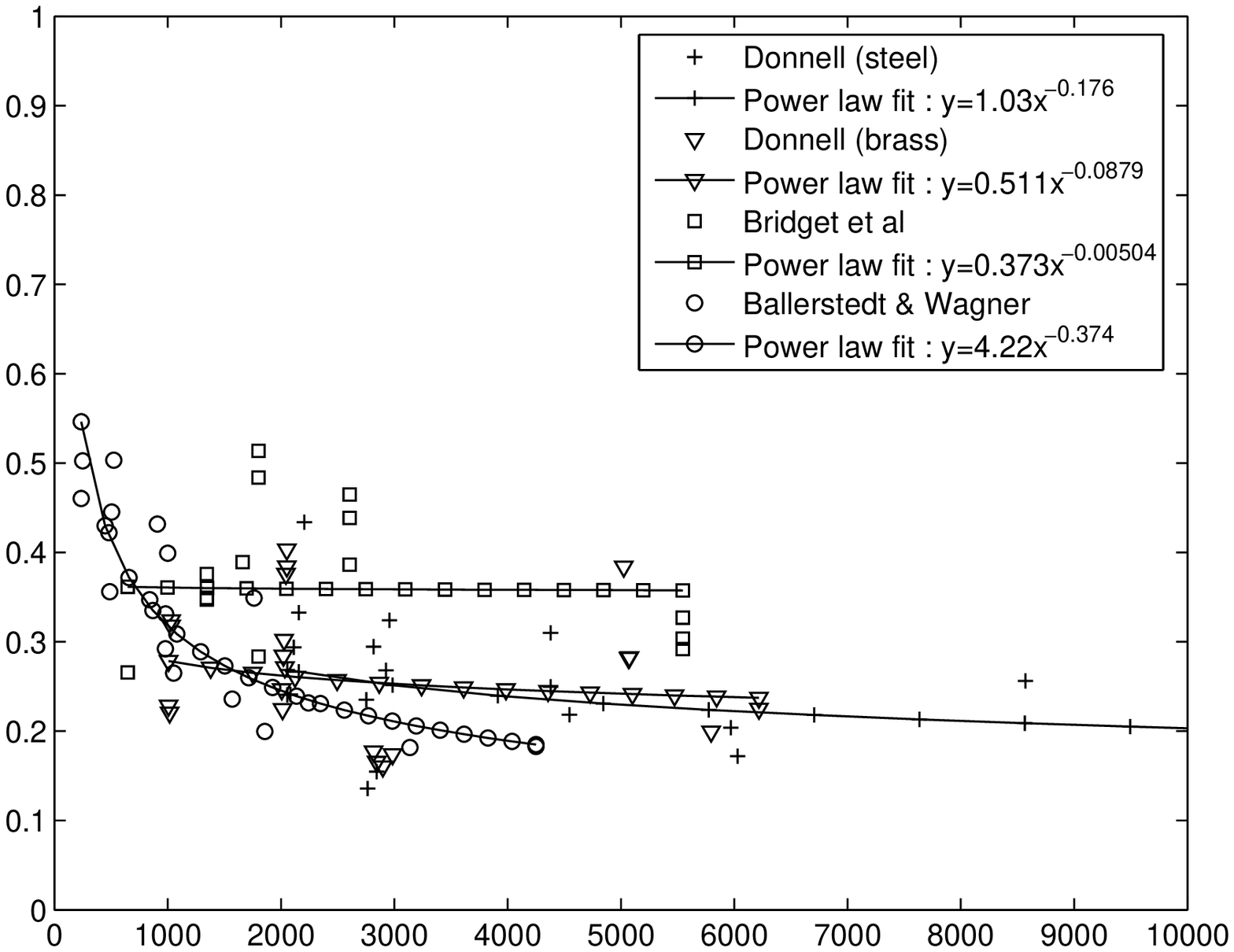}}\hskip1.5cm
\vcenter{\psfig{height=4cm,figure=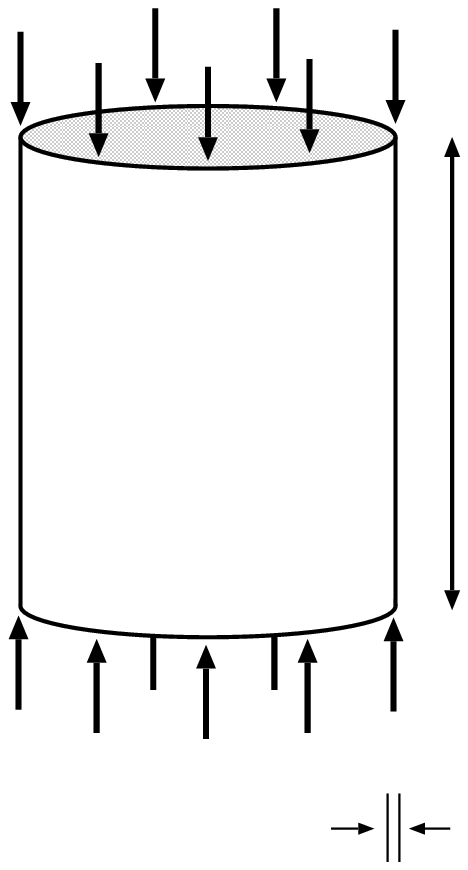}}$}
\smallskip
\caption{Experimental data from various research groups, all representing
failure loads of axially-loaded cylinders. The horizontal
axis is the ratio of cylinder length and wall thickness; the vertical axis 
is the ratio of the failure load and the theoretical critical load as predicted
for perfect cylinders. Note that all tested cylinders fail at loads significantly lower than
that predicted by theory; in some cases failure occurred at 
less than one-fifth of the theoretical load-carrying capacity. 
Power-law fitting lines are added to emphasize the dependence of failure load
on geometry. The data are from~\cite{Dnnll34,Bridgetetal,BallerstedtWagner}.}
\label{fig:experiments}
\end{figure}

In the decades before and after the second world war a central problem
was to understand the large discrepancy between theoretical
predictions and experimental observations, as shown in \figref{fig:experiments}. A variety of different explanations has been
put forward, but with the experimental work of
Tennyson~\cite{Tennyson64} and the theoretical work of
Almroth~\cite{Almroth63} it became clear that this discrepancy is
mostly due to imperfections in loading conditions and in the shape of the specimens.
Further experimental and theoretical work by many others has confirmed this
conclusion~\cite{GormanEvanIwanowski70,WeingartenMorganSeide65,Ymki84}.

For near-perfect cylinders the linear and weakly non-linear theories
(see Section~\ref{subsec:bif}) adequately describe 
the experimental buckling load\footnote{In this
  paper the terms ``experimental buckling
  load'' and ``failure load'' are used interchangeably} and the
deformation just before failure (see \emph{e.g.}~\cite{ArBa}). Cylinders used
in practical applications, however, are far too imperfect for this
weakly non-linear theory to apply, and from a practical point of view
the problem of predicting the failure load is still open.

In fact there is good reason to believe that 
it will never be possible to accurately predict failure loads 
for cylinders that are used in practice.
For simple materials, such as metals, it is believed
that current numerical methods can describe the local material
behaviour with enough accuracy that correct prediction of the
complete behaviour of the cylinder---including its failure---is
feasible.  This could be achieved provided the geometrical and
material imperfections as well as the loading conditions are
determined in sufficient detail.
%
  The difficulty
lies in the qualifier ``in sufficient detail'', since an extremely
accurate measurement of geometric imperfections would be
necessary~\cite{ArBa}, and in the design phase both
the loading conditions and the geometrical and material imperfections
in the finalized product
are only known in vague terms.  Therefore, in recent decades the
attention of theoretical research has turned to characterizing the
failure load in weaker ways, preferably in the form of a lower
(safe) bound.

\subsection{Characterizing sensitivity to imperfections}
\label{subsec:bif}
\begin{figure}[hbt]
\centerline{\psfig{figure=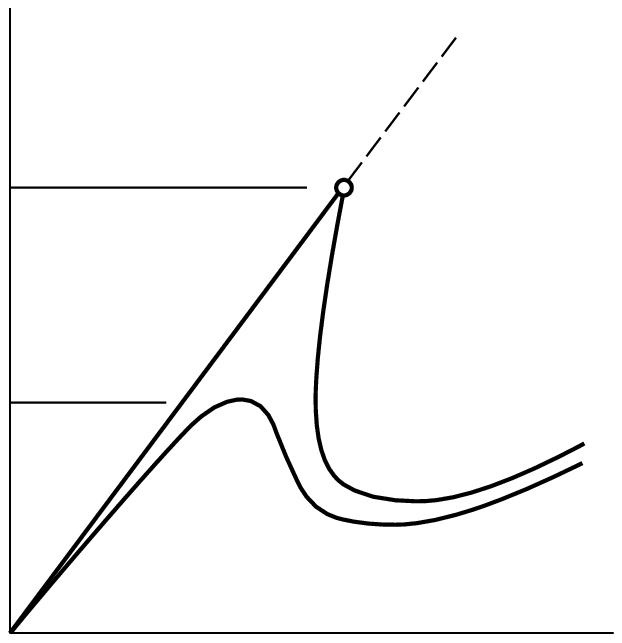,height=4cm,width=4cm}}
\caption{Illustration of perfect and imperfect bifurcation curves.}
\label{fig:bif}
\end{figure}

Viewed as a bifurcation problem, the buckling of the cylinder is a
subcritical symmetry-breaking pitchfork bifurcation (\figref{fig:bif}). Generically, imperfections in the structure eliminate the
bifurcation and round off the branch of solutions\footnote{Koiter
  actually used this elimination of a bifurcation point as a
  \emph{definition} of ``perfect system'' and ``perturbed
  system''~\cite{Ko:45}}, resulting in a turning-point at a load
$\limp$ strictly below the critical (bifurcation) load $\lambdacr$
of the perfect structure. In an experiment in which the load is slowly
increased, the system will fail (\emph{i.e.} make a large jump in state
space) at load $\limp$. 

Again, if the imperfections in geometry and loading are fully known,
then calculation of $\limp$ is a practical problem rather than a
theoretical one, and we do not address this problem here.  For the
more difficult question of characterizing $\limp$ under incomplete
information various strategies have been
proposed.  
A classical line of thought originates with Koiter~\cite{Ko:45}, in which
the imperfections are chosen \emph{a priori} within certain
finite-dimensional sets.  Common choices are the sets spanned by an
eigenvector at bifurcation, or by the eigenvectors associated with the
first $n$ bifurcations. This approach might be termed weakly
non-linear, as it is based on an expansion of the energy close to the
bifurcation point, in the directions suggested by the bifurcation
point itself. It gives predictions that are correct if the
imperfections are very small---much smaller than those 
encountered in
practice.

Since the \emph{a priori} choice of imperfections is a weak point of this
method, a natural step is to optimize over all possible perturbations.
Deml and Wunderlich pioneered this approach, in which a numerical
algorithm is used to find a ``worst geometric imperfection''~\cite{DemlWunderlich97}.
This ``worst imperfection'' is defined as the one that produces a
turning-point of minimal load. Some constraint on the magnitude of allowable
perturbations is necessary, of course, to prevent running a steam
roller over the cylinder being interpreted as an admissible
imperfection. The authors of~\cite{DemlWunderlich97}
and~\cite{WunderlichAlbertin00} first suggest constraining the
$L^\infty$-norm of the perturbation displacement, but they immediately
replace the $L^\infty$-norm by an $L^p$-norm for computational
convenience.  

This method has an interesting aspect that is
often glossed over in the engineering literature. 
By definition the failure load obtained by
this method is a lower bound for the failure load 
of all systems that have perturbations of less or equal magnitude.
The measure of magnitude is defined by the choice of constraint.
Therefore the choice of constraint on the imperfections
is a critical one, since it 
implicitly defines a class of imperfections that produce the same or
higher failure load.


\subsection{Main results}

In this paper we follow a related, but distinct, line of reasoning. 
Instead of studying actual behaviour of imperfect cylinders, we deduce
an estimate of the sensitivity to imperfections from the
energy landscape of the perfect cylinder. The final result is a lower bound 
on the failure load similar to above, and the approach gives additional insight into the problem.

The key result is the existence of
a \emph{mountain-pass point}, an equilibrium state that sits
astraddle in the energy landscape between two valleys;
one valley surrounds the unbuckled state, and the other contains many buckled,
large-deformation states. 

\begin{figure}[h]
\begin{center}
\setlength{\unitlength}{1mm}
\begin{picture}(137,113)
\color[rgb]{.5,.5,.5}
\put(0,53){\frparbcenter{\scalebox{.1}{\includegraphics[width=16in]{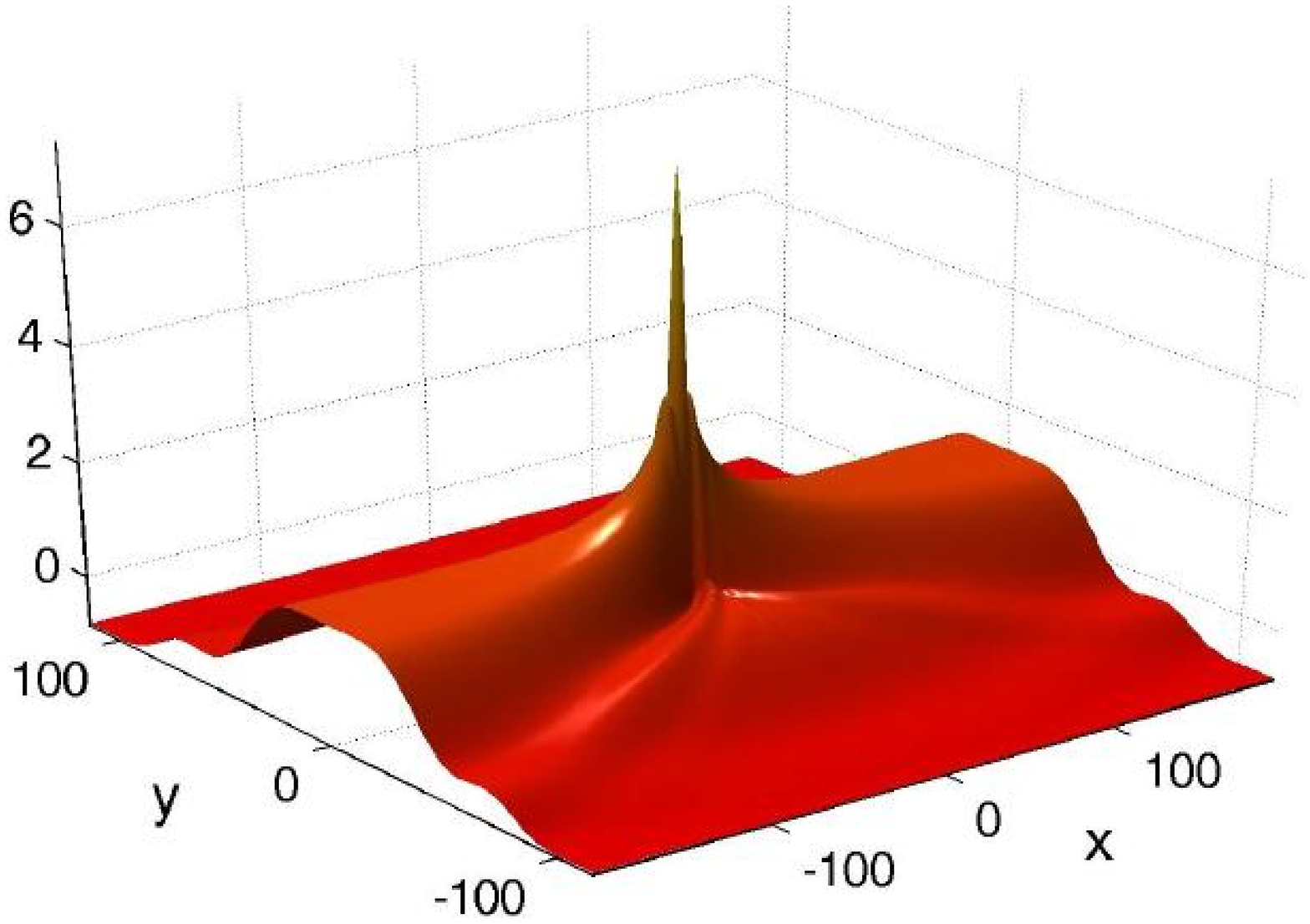}}\\\scalebox{.09}{\includegraphics[width=17in]{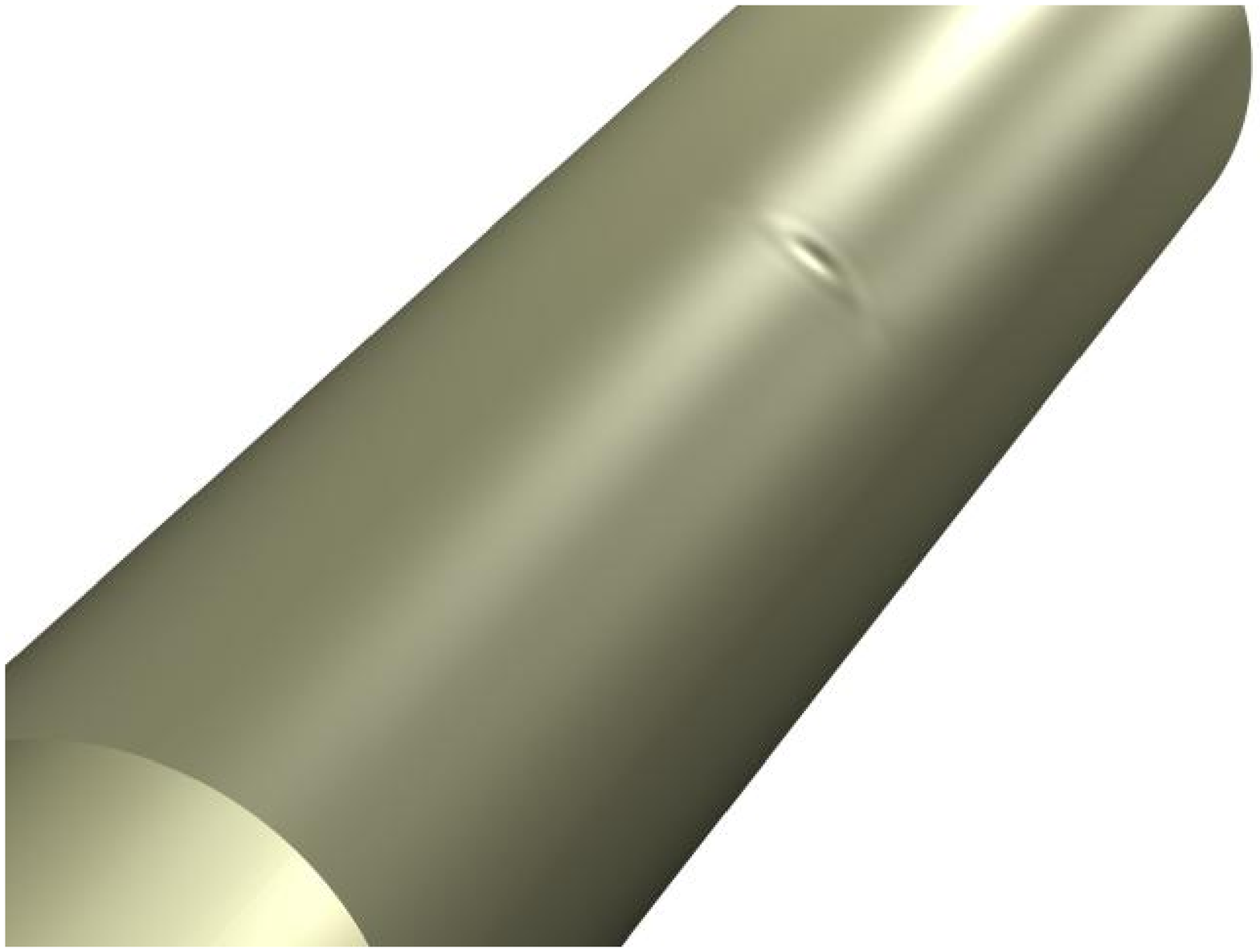}}}}
\put(47,53){\frparbcenter{\scalebox{.1}{\includegraphics[width=16in]{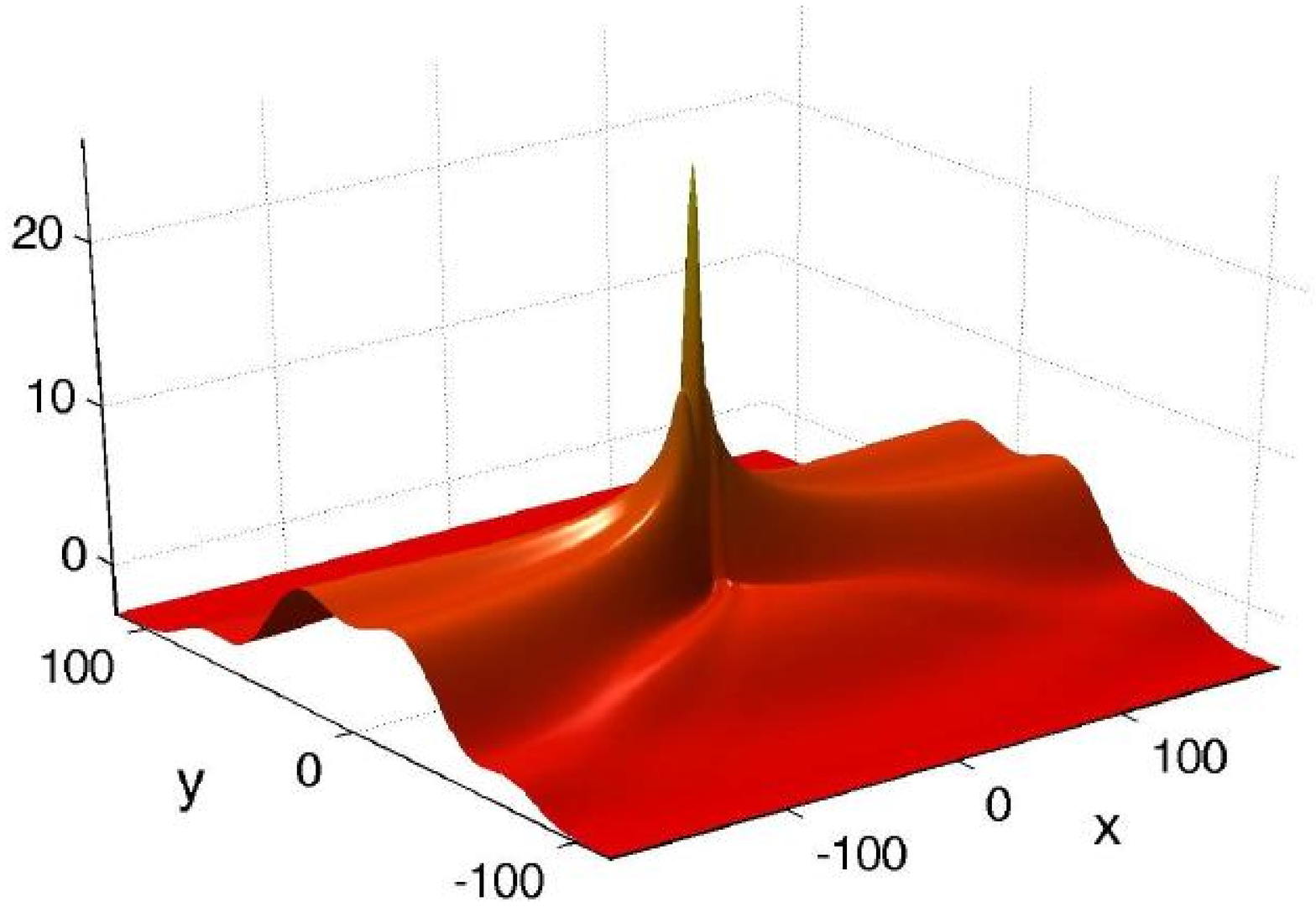}}\\\scalebox{.09}{\includegraphics[width=17in]{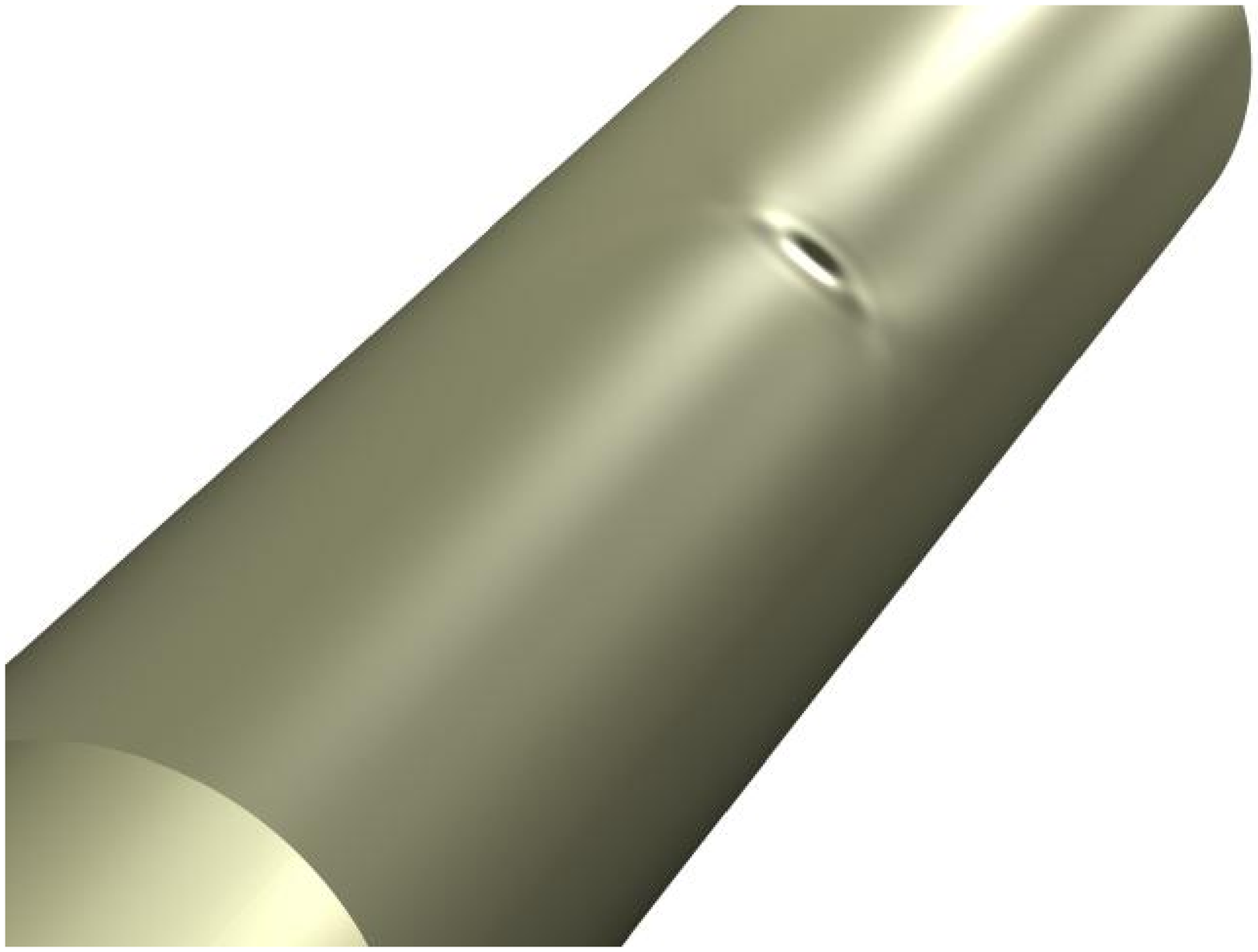}}}}
\put(94,53){\frparbcenter{\scalebox{.1}{\includegraphics[width=16in]{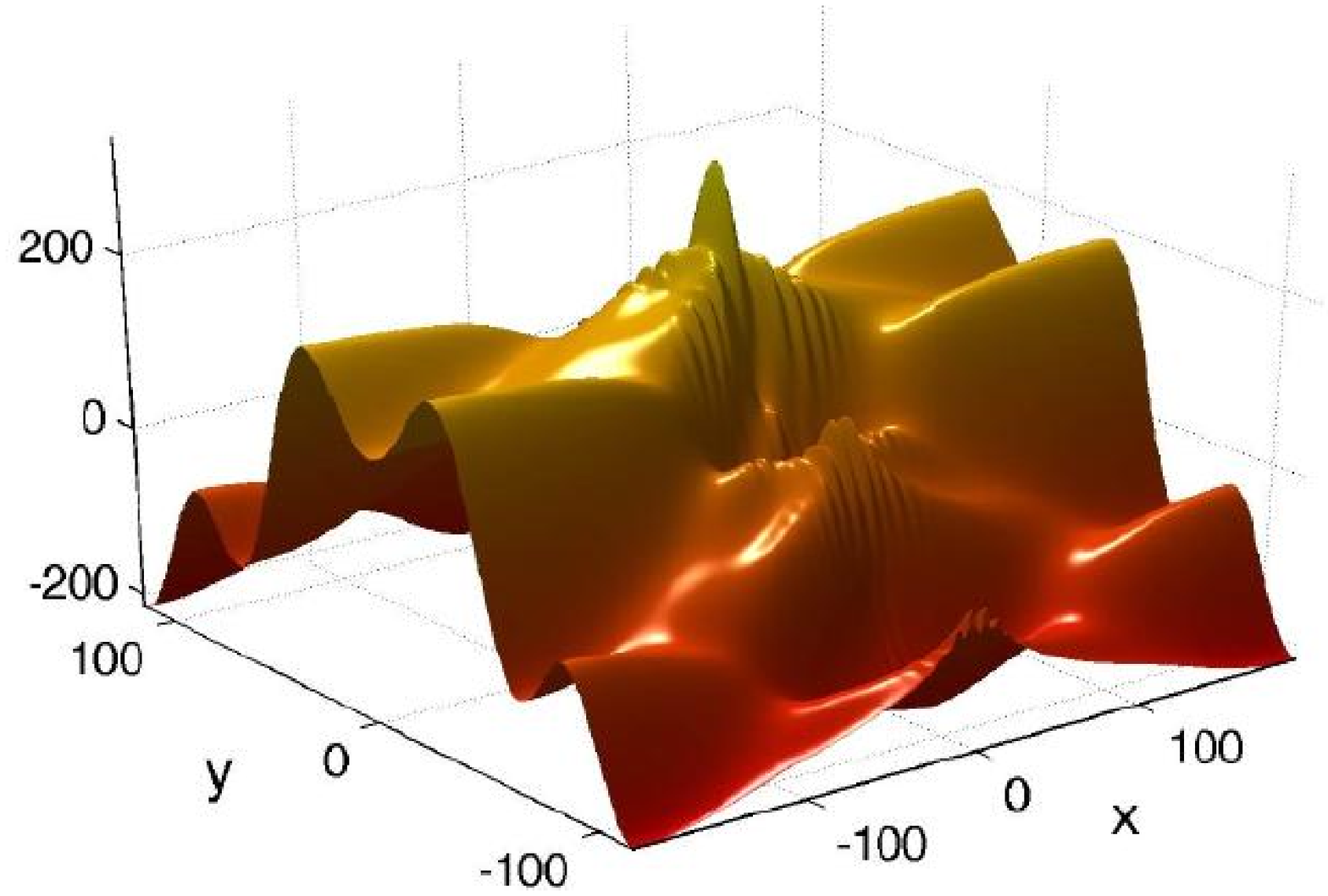}}\\\scalebox{.09}{\includegraphics[width=17in]{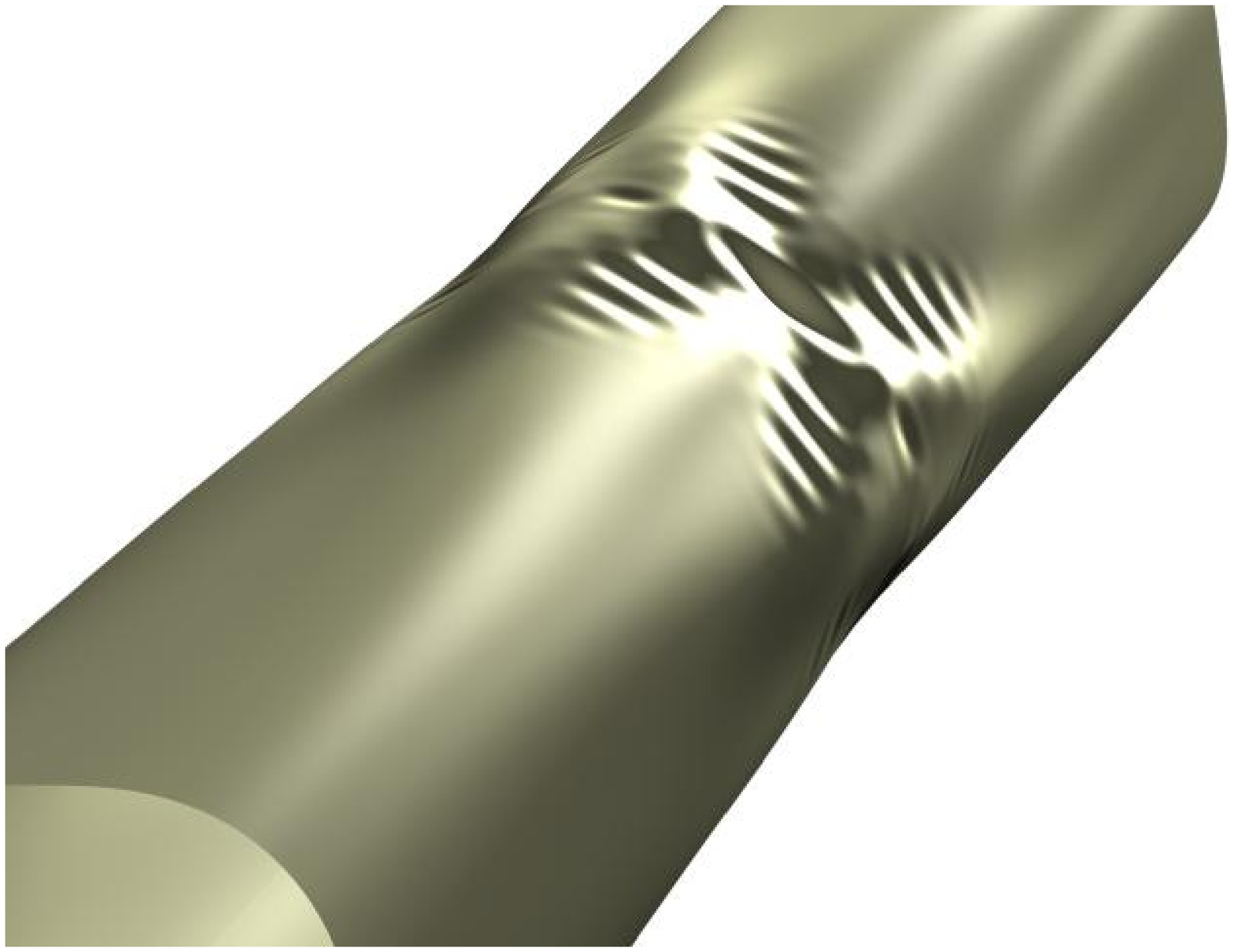}}}}
\color{black}
\put(37,55){(a)}
\put(4,79){$\wMP$}
\put(84,55){(b)}
\put(131,55){(c)}
\put(65,0){\scalebox{.644}{\includegraphics{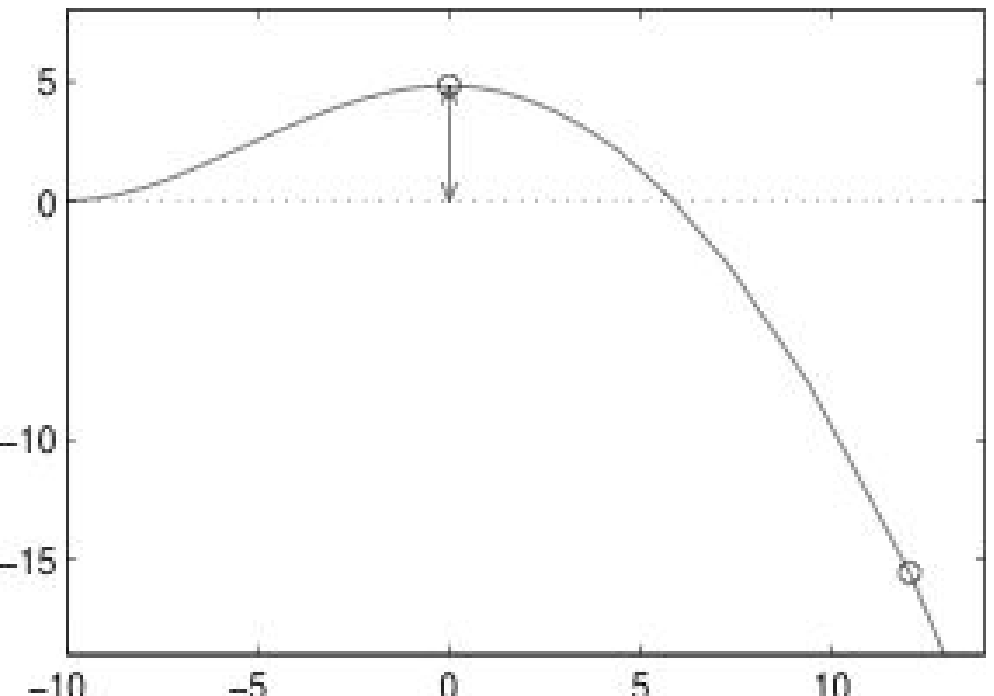}}}
\put(85.5,34){$V(\l)$}
\put(95.5,41.5){(a)}
\put(124,12){(b)}
\put(116,15){\vector(2,-3){5}}
\put(112,16){(c)}
\put(65,20){\rotatebox{90}{\hbox{\small$F_\lambda(w)$}}}
\put(90,-3){\small distance in $X$}
\put(9,25){\parbox{50mm}{(a) $\ldots V(\l)=F_\l(\wMP)\approx 4.84$\\
(b) $\ldots F_\l(w)\approx -15.57$\\
(c) $\ldots F_\l(w)\approx -5.5\cdot10^4$}}
\end{picture}
\end{center}
\caption{Part (a) shows a numerical computation of a solution $\wMP$ that is a mountain-pass point of the energy $F_\l$ for a load $\l=1.5$. We show both the
  graph of $\wMP(x,y)$ as well as its rendering on a cylinder. At
  $w_\MP$ there exist two directions in the state space $X$ in which
  the energy $F_\l$ decreases. By perturbing $\wMP$ in these
  directions and following a gradient flow of $F_\l$ we move away from
  $\wMP$. In one direction the dimple shrinks to zero (not shown) and
  in the opposite it grows in amplitude and extent (figures (b) and
  (c)).}
\label{fig:gf_mp}
\end{figure}

This mountain-pass point has a number of interesting properties:
\begin{enumerate}
\item It has the appearance of a \emph{single-dimple solution}, a
  small buckle in the form of a single dent (see
  \figref{fig:gf_mp} (a)).  Single-dimple deformations  have appeared
  in the engineering literature in a number of different ways (see
  Section~\ref{sec:conclusions}), but theoretical understanding of this
  phenomenon is still lacking. 
  Localization (concentration) of deformation is known
  commonly to appear in extended structures~\cite{HPCWWBL}, and in the
  cylinder localization in the axial direction has been studied
  theoretically and numerically~\cite{HLN1,LordChampneysHunt97,LCH2}.  
  Whether localization is
  possible in the tangential direction has been an open problem for
  some time; it is interesting that our simulations for the perfect
  structure show solutions that are localized in both axial and
  tangential directions.
\item Like all mountain-pass points this single-dimple solution is
  unstable, in the sense that there are directions in state space in
  which the energy decreases. In one direction the dimple roughly
  shrinks and disappears, and in the opposite direction it grows and
  multiplies (\figref{fig:gf_mp} (b-c)).  It is remarkable,
  however, that our numerical results indicate that the single-dimple
  solution has an alternative characterization as a \emph{constrained
    global minimizer} (a global minimizer of the strain energy under
  prescribed end shortening).
\item The equations can be rescaled so that the only remaining
  parameters are the load level and the domain. The geometry of the
  mountain-pass solution we calculate even appears to be independent
  of the domain size.
\end{enumerate}

This mountain-pass point is central in an estimate of the
sensitivity to perturbations.
For the system to escape from the neighbourhood of the unbuckled state, 
it must possess at least the energy associated with this mountain-pass point.
The mountain-pass energy level is therefore an indication of the degree of stability of the
unbuckled state.  Implicitly it defines a class of perturbations for which
the unbuckled state is stable.
This approach is related to the ``perturbation-energy'' approach first 
suggested by Kr\"oplin and co-workers~\cite{KroeplinDinklerHillmann85,DuddeckKroeplinDinklerHillmannWagenhuber89,WagenhuberDuddeck91}, but differs in 
some essential points; see the discussion in Section~\ref{sec:discussion}.

\subsection{Methods}
We use both analytical and numerical methods.
%
In Section~\ref{sec:vKD} we introduce the \vKD\ equations, which form the basis of this
paper, and rescale them in an appropriate manner. In Section~\ref{sec:MP} we 
present the functional setting that we use, show that the energy functional has the
geometry associated with a mountain pass, and prove the existence of
mountain-pass points (Lemma~\ref{lemma:PS2}).
There are certain interesting technical issues. 
By their localized nature, single-dimple solutions are most naturally defined on
an unbounded domain; however, we are only able to prove existence of mountain-pass points 
on bounded domains, and consequently we work on finite domains that become large
in the limit of thin shells. Similarly, the non-coercive nature of the energy functional 
implies that we 
prove existence of mountain-pass points for almost all load levels
(see Lemma~\ref{lemma:PS}). 

In Section~\ref{sec:numerics} we turn to numerical investigation. We use a variety of different
algorithms to find solutions of the discretized \vKD\ equations. With a discrete
mountain pass algorithm we find solutions that are, by construction, mountain-pass points.
The solution of \figref{fig:gf_mp} (a) was found in this manner. 
With a constrained
gradient flow we also find local minima of the strain energy under prescribed end shortening. 
Some of these solutions appear to coincide with those found by the discrete mountain pass
algorithm, and the mountain pass solutions are stable under this gradient flow. These
observations lead us to conjecture that the global mountain pass solution is
also a global constrained minimizer of the strain energy.
By a constrained version of the discrete mountain pass algorithm we also find critical points of higher Morse index.

Section~\ref{sec:discussion} is devoted to an interpretation of these results
in the context of imperfection sensitivity, as mentioned above, and
in Section~\ref{sec:conclusions} we wrap up with the main conclusions.

\section{The Von K\'arm\'an-Donnell equations}
\label{sec:vKD}
We consider a cylindrical shell of radius $R$, thickness $t$,
Young's modulus $E$, and Poisson's ratio $\nu$, that is subject to a compressive
axial force~$P$.  In dimensionless form the \vKD\ equations 
are given by
\begin{align}
\label{eq:main_w}
&\e^2 \Delta^2 \bw + \blambda \bw_{xx} - \bphi_{xx} - 2[\bw,\bphi] = 0, \\
&\Delta^2 \bphi + [\bw,\bw] + \bw_{xx} = 0,
\label{eq:main_phi}
\end{align}
where subscripts $x$ and $y$
denote differentiation with respect to the spatial variables, and 
the bracket $[\cdot,\cdot]$ is defined as
\[
[a,b] = \frac12 a_{xx}b_{yy} + \frac 12 a_{yy}b_{xx} - a_{xy}b_{xy}.
\]
The function $\bw$ is the inward radial displacement measured from a
unbuckled (fundamental) state, $\bphi$ is the Airy stress function,
$\e^2 = t^2(192\pi^4R^2(1-\nu^2))^{-1}$ and the nondimensional load parameter
is given by $\blambda = P(8\pi^3ERt)^{-1}$. 
The unknowns $\bw$ and $\bphi$ are defined on the two-dimensional spatial
domain $(-\ell,\ell)\times(-1/2,1/2)$, where $x\in(-\ell,\ell)$ is
the axial and $y\in(-1/2,1/2)$ is the tangential coordinate.  Since the
$y$-domain $(-1/2,1/2)$ represents the circumference of the cylinder,
the functions $\bw$ and $\bphi$ are periodic in $y$; at the axial ends
$x\in\{-\ell,\ell\}$ they satisfy the boundary
conditions
\[
  \bw_x = (\Delta \bw)_x = \bphi_x = (\Delta \bphi)_x = 0.
\]
In Appendix~\ref{app:derivation} we discuss the derivation of 
this formulation.

\medskip 

Equations~(\ref{eq:main_w}-\ref{eq:main_phi}) can be rescaled
to be parameter-independent in all but the domain of definition and
load parameter, by substituting
\begin{equation}
\label{scaling:epsilon}
\bw \mapsto \e w, \qquad \bphi \mapsto \e^2\phi, \qquad
x \mapsto \e^{1/2} x ,\qquad y \mapsto \e^{1/2} y, \qquad
\blambda \mapsto \e\lambda.
\end{equation}
so that the equations become
\begin{align}
\label{eq:main_w2}
&  \Delta^2  w +  \lambda  w_{ x x} 
  - \phi_{ x x} - 2\,[ w,\phi] = 0, \\
&\Delta^2 \phi + [ w, w] +  w_{ x x} = 0.
\label{eq:main_phi2}
\end{align}
The domain of definition of $w$ and $\phi$ is now 
\[
\Omega := (-\ell\e^{-1/2},\ell\e^{-1/2})\times (-\tfrac12\e^{-1/2},\tfrac12\e^{-1/2}),
\]
which expands to $\R^2$ as $\e\to0$. When not indicated otherwise we choose the aspect
ratio $2\ell=1$; in Section~\ref{subsec:domain_influence} we comment on 
the influence of domain size and aspect ratio.

The boundary conditions for $w$ and $\phi$ now are
\begin{subequations}
\label{def:BC}
\begin{align}
&w\text{ is periodic in $y$,\quad  and} \quad  w_x = (\Delta  w)_x = 0
\text{ at }x=\pm \tfrac12 \e^{-1/2},
\label{def:BCw} \\
&\phi\text{ is periodic in $y$,\quad  and} \quad  \phi_x = (\Delta  \phi)_x = 0
\text{ at }x=\pm \tfrac12 \e^{-1/2}.
\label{def:BCphi}
\end{align}
\end{subequations}

Equations~(\ref{eq:main_w2}--\ref{eq:main_phi2}) are related to 
the stored energy $E$ and the average axial shortening $S$ given by
\begin{equation}
\label{def:energy}
E(w) :=  \frac12\int_\Omega \left( \Delta  w^2 +  \Delta \phi^2\right), \qquad \text{and}
\qquad S(w):= \frac12 \int_\Omega  w_{ x}^2.
\end{equation}
Note that the function $\phi$ in~\pref{def:energy} is determined from
$w$ by solving~\pref{eq:main_phi2} with boundary
conditions~\pref{def:BCphi}. Furthermore, solutions
of~(\ref{eq:main_w2}--\ref{eq:main_phi2}) are stationary points of the 
total potential 
\begin{equation}
F_\lambda(w) := E(w) - \lambda S(w)
\label{def:Fl}
\end{equation}
and are also stationary points of $E$ under the constraint of constant
$S$. We use both properties below.

\section{The mountain pass: overview}
\label{sec:MP}

We briefly recall the general context of
the Mountain-Pass Theorem of Ambrosetti and Rabinowitz~\cite{AmbRab}.
Let $I$ be a functional defined on a Banach
space $X$, and let $w_1,w_2$ be two distinct points in $X$. 
Consider the
family $\Gamma$ of all paths in $X$ connecting $w_1$ and $w_2$ and
define
\begin{equation}
\label{def:c}
c = \inf_{\gamma\in\Gamma}\max_{w\in\gamma} I(w) \ ,
\end{equation}
that is the infimum of the maxima of the functional $I$ along paths in
$\Gamma$. If $c>\max\{I(w_1),I(w_2)\}$, then the paths have to cross a
``mountain range'' and one may conjecture that there exists a critical point $w_\MP$ of $I$ at the
level $c$, called a mountain pass point.

We will apply this idea to the \vKD-equations in the following way.
We take for $I$ the total potential $F_\l$ (see~\pref{def:Fl}) at some fixed value of $\l$,
and for the end point $w_1$ the origin. 
We will obtain a mountain-pass solution by the following steps:
\begin{itemize}
\item[MP1.] We first show that $w_1=0$ is a local minimizer: there 
are $\varrho, \alpha>0$ such that $F_\l(w)\geq \alpha$ for all $w$ 
with $\|w\|_X=\varrho$ (Lemma~\ref{lemma:origin});
\item[MP2.] If $\e$ is small enough, then 
there exists $w_2$ with $F_\l(w_2) \leq 0$ (Corollary~\ref{cor:periodics}).
\item[MP3.] Given a sequence of paths $\gamma_n$ that  approximates the infimum in~\pref{def:c},
we extract a (Palais-Smale) sequence of points $w_n\in \gamma_n$, each one close to the 
maximum along $\gamma_n$, and show that this sequence converges in an appropriate 
manner (Lemmas~\ref{lemma:PS} and~\ref{lemma:PS2}). 
\end{itemize}

In this way it follows that there
exists a mountain-pass critical point~$w$ with $F_\lambda(w)= c$,
provided that $\e$ is sufficiently small (or that the domain is sufficiently large).
For technical reasons (lack of coerciveness of the functional $F_\l$)
this procedure can be performed only for almost all $0<\l<2$ (see Lemma~\ref{lemma:PS}).


In the rest of this section we detail the steps outlined above.

\subsection{Choice of spatial domain} We are interested in mountain-pass solutions of the system of
equations (\ref{eq:main_w2}--\ref{eq:main_phi2}) 
that are quasi-independent of the domain size $\e^{-1/2}$, 
in the sense that they converge to 
a non-trivial solution on $\R^2$ as $\e\to0$. This point of view suggests
considering the problem on the whole of $\R^2$ rather than on a sequence
of domains of increasing size; however, there are two reasons not to do this.
To start with, the numerical calculations described below are necessarily done
on a bounded domain; more importantly, 
for the proof of existence of mountain-pass  points 
boundedness of the domain is necessary. For these 
reasons we concentrate on bounded domains, while
keeping the context of the unbounded domain in mind. 

\subsection{Functional setting and linearization}
\label{sec:linearization}

We introduce a functional setting for the functions $w$ that
is suggested by the linearization of the stored energy functional $E$.
Writing $\phi=\phi_1+\phi_2$, where
\begin{equation}
\Delta^2 \phi_1 = -w_{xx} \qquad\text{and}\qquad 
\Delta^2 \phi_2 = -[w,w],
\label{def:phi_1}
\end{equation}
we can expand the energy functional $E$ as 
\begin{equation}
\label{eq:E_develop}
E(w) = \frac12 \int_\Omega \Delta w^2 + \frac12 \int_\Omega \Delta \phi_1^2
 + \int_\Omega \Delta\phi_1\Delta\phi_2 + \frac12 \int_\Omega \Delta\phi_2^2.
\end{equation}
Since $\phi_2$ is quadratic in $w$, the second derivative of $E$ is 
given by
\[
d^2E(0) \cdot u \cdot v = \int_\Omega \Delta u\Delta v + \int_\Omega \Delta \phi^u_1\Delta\phi_1^v,
\]
where $\phi_1^{u,v}$ are obtained from $u$ and $v$ by~\pref{def:phi_1}.
Inspired by this linearization of $E$ we define
\[
X=\left\{\psi\in H^2(\Omega) :  \psi_x(\pm\tfrac12 \e^{-1/2},\cdot)=0,\ \psi \text{ is periodic in } y,
  \text{ and }\,\int_\Omega\psi = 0\right\}
\]
with norm
\[
  \|w\|^2_X = \int_\Omega\bigl(\,\Delta w^2 +\Delta\phi_1^2\,\bigr),
\]
where $\phi_1\in H^2(\Omega)$ is the unique solution of
\[
  \Delta^2\phi_1 = -w_{xx}, \qquad \phi_1 \text{ satisfies~\pref{def:BCphi}, }
\qquad \text{and}\qquad \int_\Omega \phi_1 = 0.
\]
This norm is equivalent to the $H^2$-norm on the set $X$, and with the
appropriate inner product the space $X$ is a Hilbert space. 

We now address the requirements of the mountain pass theorem mentioned above
in MP1--MP3.

\subsection{The origin is a local minimizer}

The norm in $X$ is related in a natural manner with the shortening~$S$,
as is demonstrated by the (sharp) estimate
\begin{equation}
2S(w) = \int_\Omega w_x^2 = -\int_\Omega ww_{xx} = \int_\Omega w\Delta^2 \phi_1
= \int_\Omega \Delta w\Delta \phi_1 \leq \frac12 \int_\Omega \Delta w  ^2 + \frac12 \int_\Omega \Delta \phi_1^2 
= \frac12 \Mod{w}_X^2.
\label{eq:tmp1}
\end{equation}
This inequality strongly suggests 
that for $\lambda<2$ the origin is
a strict local minimum for the  functional $F_\lambda(w) = E(w)-\lambda S(w)$. 

\begin{lemma}
\label{lemma:origin}
\begin{enumerate}
\item There exists a constant $C$, dependent on $\Omega$, such that 
\begin{equation}
\label{est:nonlinear}
\mod{\int_\Omega \Delta \phi_1\Delta \phi_2} \leq C \Mod{w}_X^3.
\end{equation}
\item \label{lemma:origin:2}
For any $0<\lambda<2$ there exists $\varrho>0$
such that
\[
\inf \bigl\{\;  F_\lambda(w): \Mod w_X = \varrho\; \bigr\} > 0.
\]
\end{enumerate}
\end{lemma}

\begin{proof}
Split
$\phi = \phi_1 + \phi_2$ as in~\pref{def:phi_1}, and 
note that the function $\nabla \phi_1$ is bounded in $L^\infty$ by the
Sobolev imbedding $\{\psi \in H^3: \int \psi = 0\} \hookrightarrow L^\infty$:
\[
\Mod{\nabla \phi_1}^2_{L^\infty} 
\leq C \Mod{\Delta^2\phi_1}^2_{L^2} = C  \Mod{w_{xx}}_{L^2}^2.
\]
The third term on the right-hand side of~\pref{eq:E_develop} can now be rewritten as
\[
\int_\Omega \Delta \phi_1\Delta \phi_2 = 
\int_\Omega \phi_1[w,w] = \int_\Omega \phi_1 (w_{xx}w_{yy} - w_{xy}^2)
  = \int_\Omega (\phi_{1y}w_xw_{xy} - \phi_{1x}w_xw_{yy}),
\]
which we estimate by
\[
\mod{\int_\Omega \Delta \phi_1\Delta \phi_2} \leq 2 \Mod{\nabla \phi_1}_{L^\infty}
    \Mod{w_x}_{L^2} \Mod{\Delta w}_{L^2}
  \leq C \Mod{w_x}_{L^2} \Mod{\Delta w}_{L^2}^2 \leq  C \sqrt{S(w)} \Mod w_X^2,
\]
thus proving the first part of the Lemma.

Since $\lambda< 2$, choose $0< \varrho < (2-\l)/2C$ and define
 $\eta = \tfrac12 (1-C\varrho-\l/2)>0$.
Then 
%
on the set $\mathcal C = \{w: \Mod w_X = \varrho \}$,
using \pref{eq:tmp1}, we find that 
\begin{align*}
F_\lambda(w) &= E(w) - \lambda S(w)\\
 &\geq \frac{1}{2}\Mod w_X^2  - C\sqrt {S(w)} \Mod w_X^2 - \lambda S(w) \\
 &\geq \frac{\varrho^2}{2}  - C\varrho^3 \frac12 - \lambda \frac {\varrho^2}4\\
 &= \frac{\varrho^2}{2}\left( 1 - C \varrho- \frac\lambda2\right) \\
 &= \eta \varrho^2.
\end{align*}
\end{proof}

Part~\ref{lemma:origin:2} of this lemma implies that by choosing $w_1$ to be
the origin we have shown condition MP1.

\begin{remark}
Although the inequality~\pref{eq:tmp1} suggests that the origin should 
be a local minimizer for any domain $\Omega$, bounded or not, the proof above 
only applies to bounded domains. F. Otto has constructed a proof
of this result that is valid on any domain (personal communication). Interestingly,
this proof uses not only the cubic energy term 
$\int_\Omega \Delta \phi_1\Delta \phi_2$ but also the quartic term $\int_\Omega \Delta \phi_2^2$,
and appears to break down without this latter term.
\end{remark}

\subsection{Periodic solutions exist with negative \boldmath $F_\l$}
\label{sec:Yoshimura}

To satisfy MP2 we show in this section that for any $\lambda>0$ functions
$w\in X$ exist for which $F_\lambda(w) = E(w)-\lambda S(w) < 0$. To do this we
construct a sequence of functions $w_\d$ with specific scaling properties: 
\begin{lemma}
\label{th:scaling_we}
There exists a sequence of functions $w_\d$, $1$-periodic on $\R^2$, 
such that 
\begin{multline*}
\int\limits_{[-1/2,1/2]^2} w_{\d x}^2 \sim 1,\qquad 
\int\limits_{[-1/2,1/2]^2} \Delta w_\d^2 = O(\d^{-1}),\\
\text{and}\qquad
\int\limits _{[-1/2,1/2]^2}\Delta \phi_\d^2 = O(\d^{2-\alpha}) \qquad
\text{as $\delta\to0$},
\end{multline*}
for any $\alpha>0$.
Here the functions $w_\d$ and $\phi_\d$ solve equation~\pref{eq:main_phi2} with periodic
boundary conditions. In addition, $w_\d$ and $\phi_\d$ satisfy the
boundary conditions~\pref{def:BC} on the boundary of $[-1/2,1/2]^2$.
\end{lemma}

\begin{figure}[ht]
\centering
\scalebox{.16}{\includegraphics[width=17in]{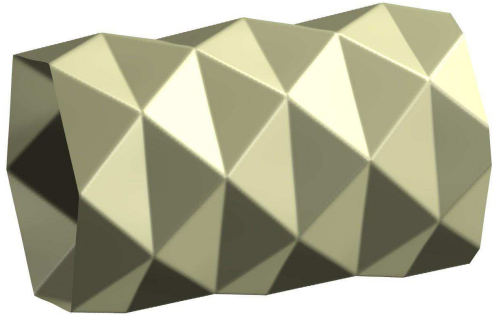}}
\caption{Yoshimura folding pattern.}
\label{fig:w1}
\end{figure}

The proof, given in the Appendix, is inspired by the so-called 
\emph{Yoshimura pattern}~\cite{Yoshimura}, a folding pattern by which
a flat sheet of 
paper, or a cylindrical sheet of thin material, can be folded into 
a macroscopically cylindrical structure with zero Gaussian curvature
but locally infinite total curvature (\figref{fig:w1}). The functions $w_\d$ are smoothed
versions of the Yoshimura pattern, adapted to the geometrically linear setting
of the von K\'arm\'an-Donnell equations, and $\d$ measures the
width of the fold.

\begin{corollary}
\label{cor:periodics}
\begin{enumerate}
\item Fix $\l>0$. If $\e$ is sufficiently small, then there exists $w\in X$ such that
$F_\l(w)<0$. 
\item Fix $\e$ sufficiently small. Then there exists $\lambda_0(\e)\in[0,2)$ such that for
all $\lambda> \lambda_0$ there exists $w\in X$ with $F_\lambda(w)<0$.
\end{enumerate}
\end{corollary}

\begin{proof}
By scaling the functions $w_\d$ of Lemma~\ref{th:scaling_we} the claims can be fulfilled:
let $\delta=\e^{2/3}$, and set
\[
\tilde w_\e(x,y) = \e^{-1}w_{\e^{2/3}}(x\e^{1/2},y\e^{1/2}),\qquad  \tilde \phi_\e(x,y) = \e^{-2} \phi_{\e^{2/3}}(x\e^{1/2},y\e^{1/2});
\]
then $\tilde w_\e\in X$,  and~\pref{eq:main_phi2} is invariant under this scaling;
in addition, choosing $\alpha=1/6$,
\[
Q_\e := \frac{\ds\int_{\Omega_\e} 
  \bigl[\Delta \tilde w_\e^2 + \Delta \tilde \phi_\e^2\bigr]}
  {\ds\int_{\Omega_\e} \tilde w_{\e x}^2}
= O\bigl(\e^{1/6}\bigr) \qquad \text{as $\e\to0$}.
\]
Therefore $\lim_{\e\to0} Q_\e = 0$, proving the first claim. 
For the second claim we fix $\e$ such that $Q_\e<2$;
then for all $\lambda>Q_\e$, $F_\lambda(\tilde w_\e)<0$.
\end{proof}

\medskip


\subsection{Convergence of selected sequences}

For given $\l\in(0,2)$ and for sufficiently small $\e>0$, 
the two previous sections provide two points: the origin $w_1=0$ that
satisfies MP1 and a point $w_2$ with $F_\lambda(w_2)<0$, such that
\begin{equation}
\label{def:c2}
c(\l) := \inf_{\gamma\in\Gamma}\max_{w\in\gamma} F_\lambda(w) >0,
\end{equation}
where $\Gamma$ is the set of curves connecting $0$ and $w_2$,
\[
\Gamma = \{\gamma\in C([0,1];X): \gamma(0) = 0,\ \gamma(1) = w_2\},
\]
(actually, $\Gamma$ depends on $\l$ through the dependence on $w_2$, but
$w_2$ can be taken independent of $\lambda$  in a neighbourhood of a given $\l\in(0,2)$).

We were unable to prove the classical Palais-Smale condition, which reads
\begin{quote}
For any sequence $w_n\in X$ such that $F_\lambda(w_n)\to c$
and $F_\lambda'(w_n) \to 0$ in $X'$,
there exists a subsequence that converges in $X$. 
\end{quote}
The difficulty lies in the lack of coerciveness of the functional $F_\l$: the quotient
$F_\l(w)/\Mod{w}^2_{X}$ is not bounded away from zero, implying that Palais-Smale
sequences may be unbounded in $X$. 

The ``Struwe monotonicity trick''~\cite{Struwe90} 
provides a way of proving the boundedness of Palais-Smale
sequences for at least \emph{almost all} $\l\in(0,2)$. The pertinent
observation is that for fixed~$w$, 
$F_\l(w)$ is decreasing in $\l$; by consequence $c(\l)$ is a decreasing function of $\l$
and therefore differentiable in almost all $\l\in(0,2)$. If $\gamma(t)$ is the
highest point of a near-optimal curve $\gamma$ at some $\l_0$, then $c'(\l_0)$
should be close to $-S(\gamma(t))$. Finiteness of $c'$ at $\l_0$ thus
implies that near-mountain-pass points have bounded $S$, and
this additional information suffices for the construction
of bounded sequences:

\begin{lemma}
\label{lemma:PS}
Let $\l\in(0,2)$ be such that $c'(\l)$ exists. Then there exists a bounded 
Palais-Smale sequence $w_n$, \emph{i.e.} a sequence that satisfies
\begin{enumerate}
\item $w_n$ is bounded in $X$;
\item $F_\l'(w_n)\to 0$ in $X'$ and $F_\l(w_n)\to c(\l)$;
\item there exists a sequence of curves $(\gamma_n)\subset \Gamma$ such 
that $w_n\in \gamma_n([0,1])$ and $\max_{t\in[0,1]} F_\l(\gamma_n(t))\to c(\l)$.
\label{mpcond3}
\end{enumerate}
\end{lemma}

In~\cite{SmetsVandenBerg02} this same argument was used to 
study mountain-pass points for the related one-dimensional functional
\[
J_\l(u) = \int_\R \Bigl\{\frac12 {u''}^2 - \frac\l2 {u'}^2 + F(u)\Bigr\},
\]
where $F$ is a non-negative double- or single-well potential.
The proof of Lemma~\ref{lemma:PS} is a word-for-word repeat of the 
proof of \cite[Proposition 5]{SmetsVandenBerg02},
and we omit it here.

\medskip

\begin{lemma}
\label{lemma:PS2}
The sequence $w_n$ given by Lemma~\ref{lemma:PS} is compact in $X$,
and a subsequence converges to a stationary point $w\in X$ of $F_\l$. 
\end{lemma}

Strictly speaking the stationary point given by this lemma may not be 
a mountain-pass point itself, in the sense that there may not be a curve $\gamma\in \Gamma$
of which $w$ is the highest point. Property~\ref{mpcond3} of Lemma~\ref{lemma:PS}, however, 
states that $w$ has an approximate mountain-pass character.

\begin{proof}
We extract a subsequence that converges weakly in $X$ and strongly in $H^1$ and $L^\infty$ to a limit $w$.
Defining $\phi_{1n}$ and $\phi_{2n}$ by~\pref{def:phi_1}, 
we find that the right-hand sides in~\pref{def:phi_1} are bounded in $L^2$ and $L^1$,
and therefore that $\phi_{1n}$ and $\phi_{2n}$ converge
strongly (up to extracting a subsequence, which we do without changing 
notation) in $H^2$ to
functions $\phi_{1,2}$. Both functions $\phi_{1,2}$ are again related
to $w$ by~\pref{def:phi_1}; for $\phi_2$ this follows from remarking that
for given $\zeta\in C_c^\infty(\Omega)$,
\[
\lim_{n\to\infty} \int_\Omega \zeta[w_n,w_n] = \lim_{n\to\infty} \int_\Omega w_n[\zeta,w_n] = \int_\Omega w[\zeta,w] = \int_\Omega \zeta[w,w],
\]
so that the right-hand side converges in the sense of distributions. Similarly
it follows from the strong $H^2$-convergence of 
$\phi_n = \phi_{1n} + \phi_{2n}$ that
\[
\lim_{n\to\infty} \int_\Omega \phi_n[w_n,w-w_n] 
= \lim_{n\to\infty}\int_\Omega w_n[\phi_n,w-w_n] = 0.
\]

To show that $w_n$ converges strongly in $X$, 
note that the derivative $F'_\l(w_n)$ can be characterized as
\[
 F_\l'(w_n)\cdot v = \int_\Omega \Delta w_n \Delta v 
  - \int_\Omega \phi_n \bigl(v_{xx} + 2[w_n,v]\bigr)
  - \l \int_\Omega w_{nx}v_x.
\]
We now calculate
\begin{align*}
\lim_{n\to\infty} \biggl\{\int_\Omega \Delta w^2 &- \int_\Omega \Delta w_n^2 \biggr\} 
= \lim_{n\to\infty} \int_\Omega \Delta w_n\Delta(w-w_n) \\
&= \lim_{n\to\infty} \bigg\{  F_\l'(w_n)\cdot(w-w_n) 
  + \int_\Omega \phi_n \bigl( (w-w_n)_{xx} + 2[w_n,w-w_n]\bigr) \\
&\hskip0.5\hsize  + \l\int_\Omega w_{nx}(w-w_n)_x \biggr\} \\
&= 0.
\end{align*}
The strong convergence of $w_n$ in $X$ now follows from the uniform convexity of $X$.
\end{proof}

\section{Numerical Results}
\label{sec:numerics}

\subsection{Description of the algorithm}
Our goal in this section is to find, numerically, critical points of
$F_\lambda$. Although we will focus on mountain-pass points described
above and sketch the method used to find them, numerical
approximations of other critical points of $F_\lambda$ will be shown
as well. More details on all the numerical methods used are given in a
companion paper~\cite{HoLoPe2}.

In order to employ the mountain-pass algorithm we
discretize~(\ref{eq:main_w2}--\ref{def:Fl}) using finite differences.
The algorithm was first proposed in \cite{ChMcK1} for a second order
elliptic problem in 1D. It was later used in \cite{HoMcK} for a fourth-order 
problem in 2D. 

\begin{figure}[htbp]
  \begin{center}
    \setlength{\unitlength}{1mm}
    \begin{picture}(100,67)
      \put(0,0){\scalebox{.2376}{\includegraphics[width=16.56in]{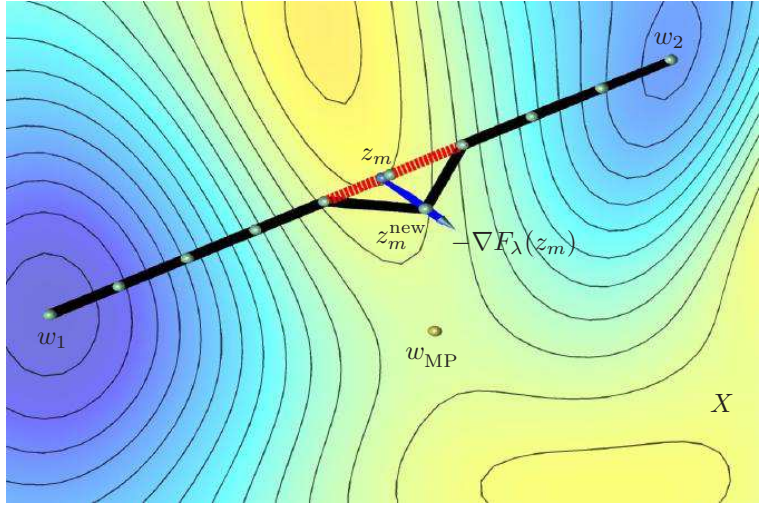}}}
      \put(4,21){$w_1$}
      \put(86,61){$w_2$}
      \put(47,45.5){$z_m$}
      \put(49,34.5){$z_m^\mathrm{new}$}
      \put(59,33.5){$-\nabla F_\lambda(z_m)$}
      \put(53,19){$w_\MP^{}$}
      \put(93.5,12){$X$}
    \end{picture}
  \end{center}
  \caption{Deforming the path in the main loop of the mountain pass
    algorithm: point $z_m$ is moved a small distance in the direction
    $-\nabla F_\lambda(z_m)$ and becomes $z_m^\mathrm{new}$. This step
    is repeated until the mountain pass point $\wMP$ is reached.}
  \label{fig:mpa}
\end{figure}

The main idea of the algorithm is illustrated in \figref{fig:mpa}.
We take a discretized path connecting $w_1=0$ with a point $w_2$ such
that $F_\lambda(w_2)<0$. 
After finding the point $z_m$ at which $F_\l$ is maximal along the path,
this point is moved
%
a small
distance in the direction of the steepest descent of $F_\lambda$.
Thus the path has been deformed and the maximum of $F_\lambda$
lowered. This deforming of the path is repeated until the maximum along
the path cannot be lowered any more: a critical point $\wMP$ has been
reached.

\begin{figure}[h!]
\begin{center}
\setlength{\unitlength}{1mm}
\begin{picture}(137,186)
\color[rgb]{.5,.5,.5}
\put(0,126){\frparbcenter{\scalebox{.1}{\includegraphics[width=16in]{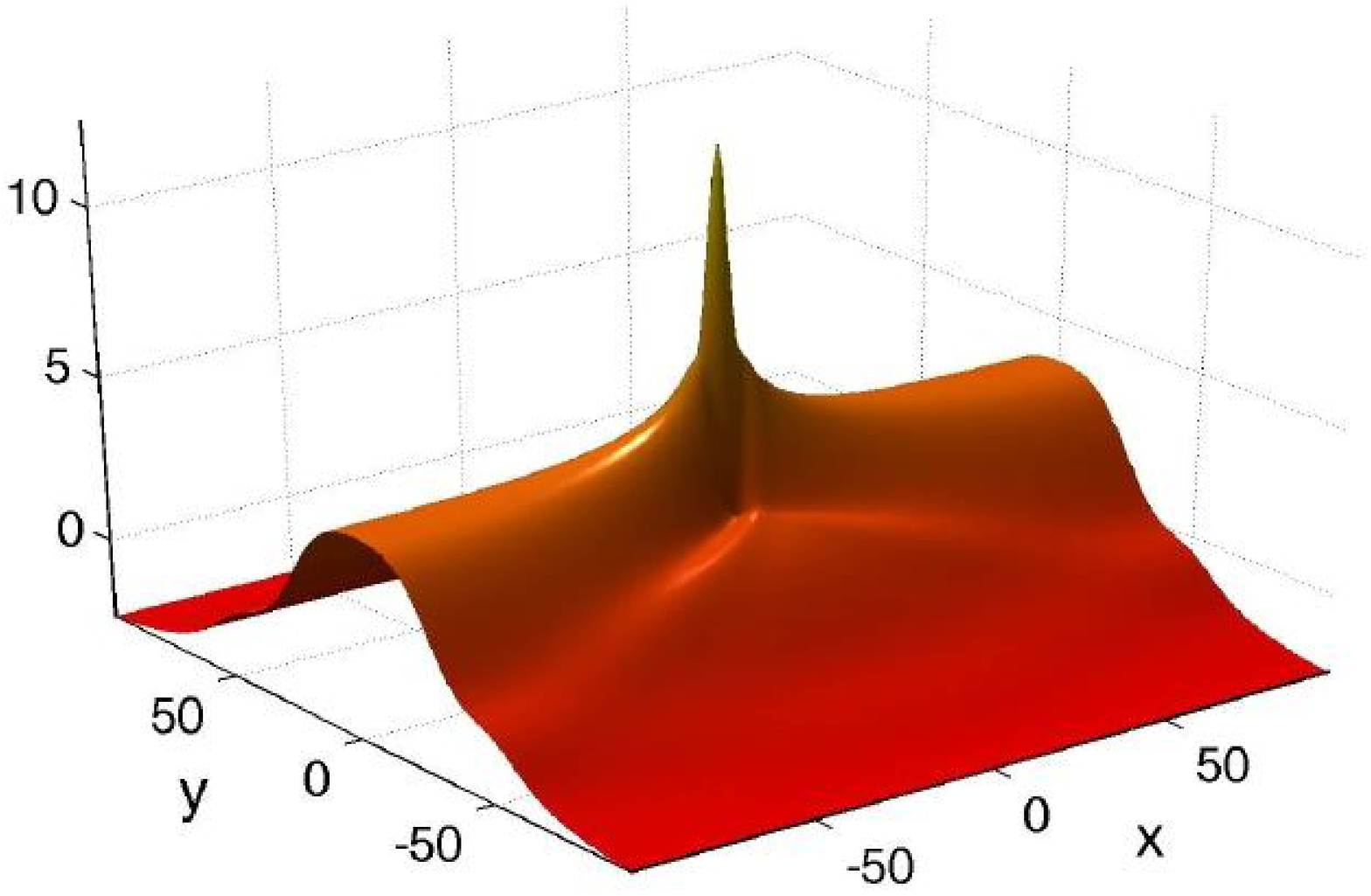}}\\\scalebox{.09}{\includegraphics[width=17in]{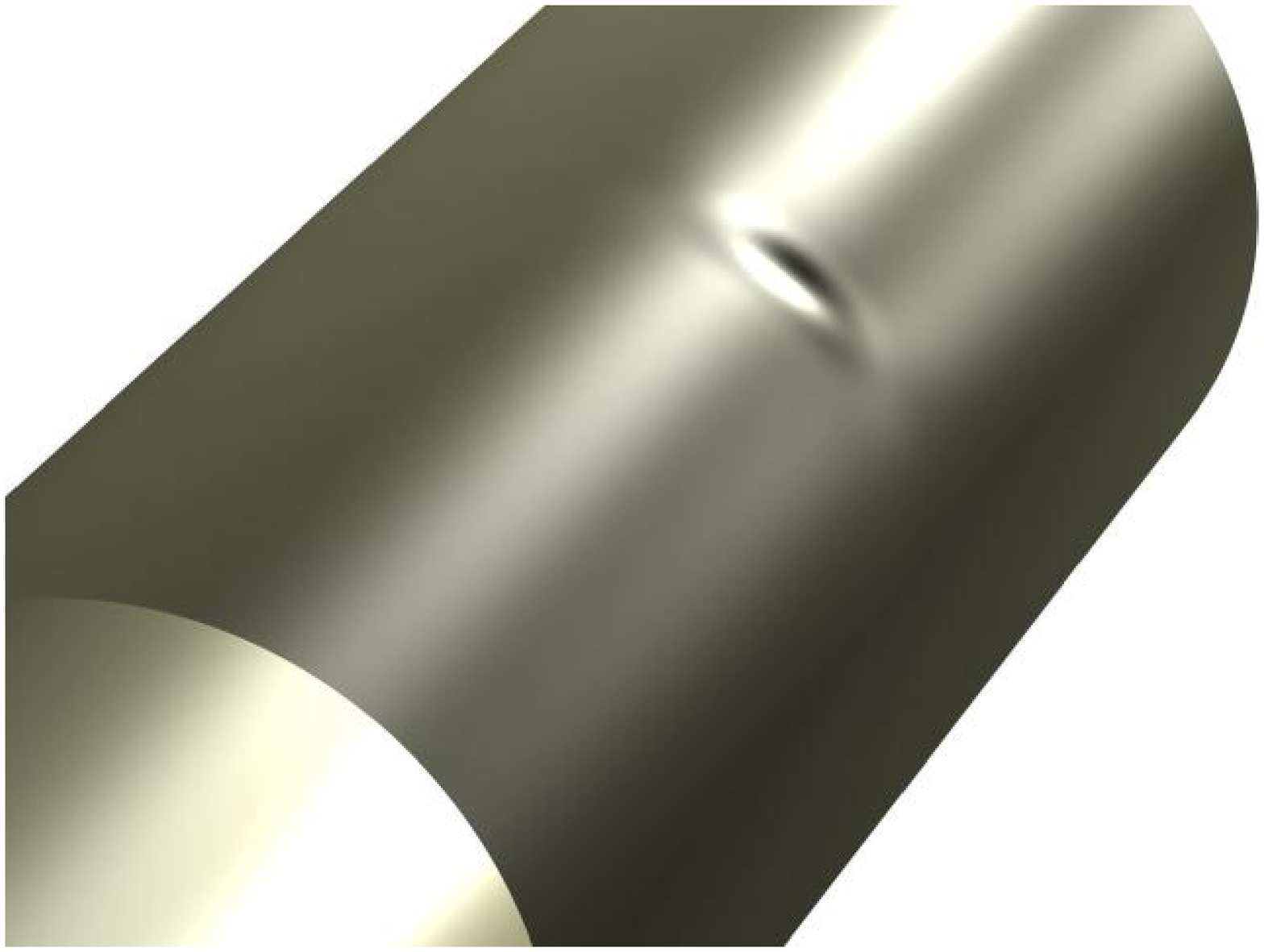}}}}
\put(47,126){\frparbcenter{\scalebox{.1}{\includegraphics[width=16in]{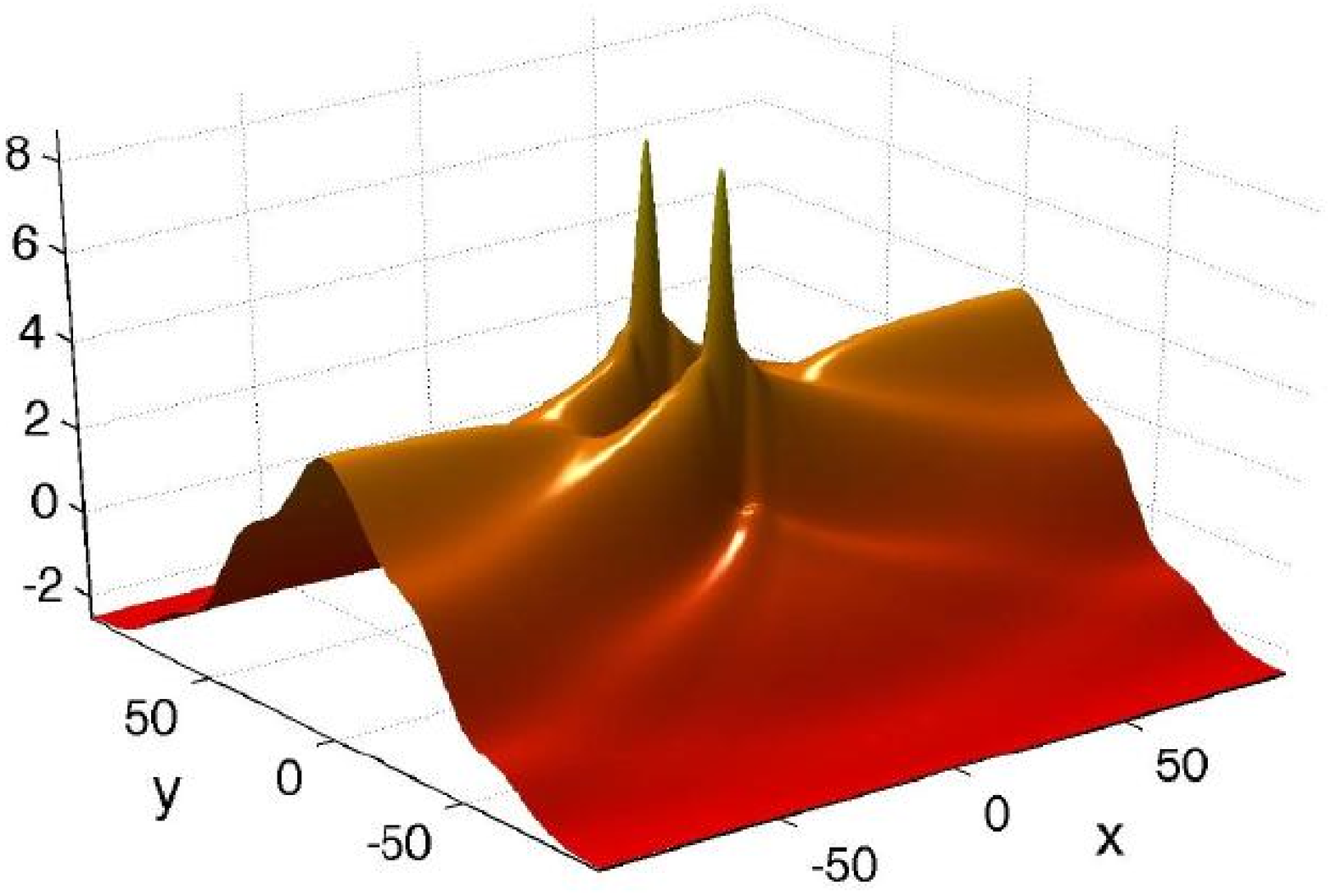}}\\\scalebox{.09}{\includegraphics[width=17in]{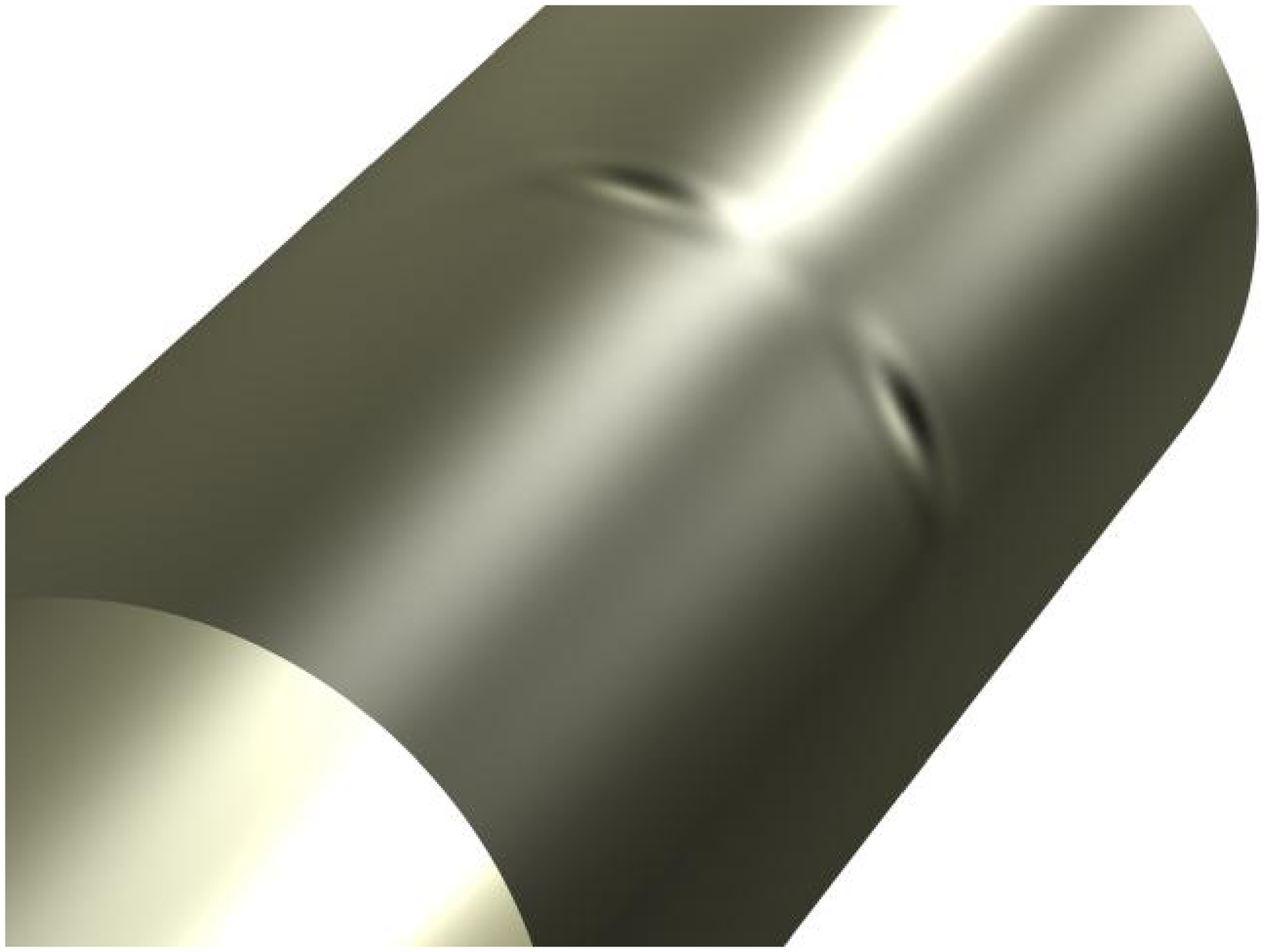}}}}
\put(94,126){\frparbcenter{\scalebox{.1}{\includegraphics[width=16in]{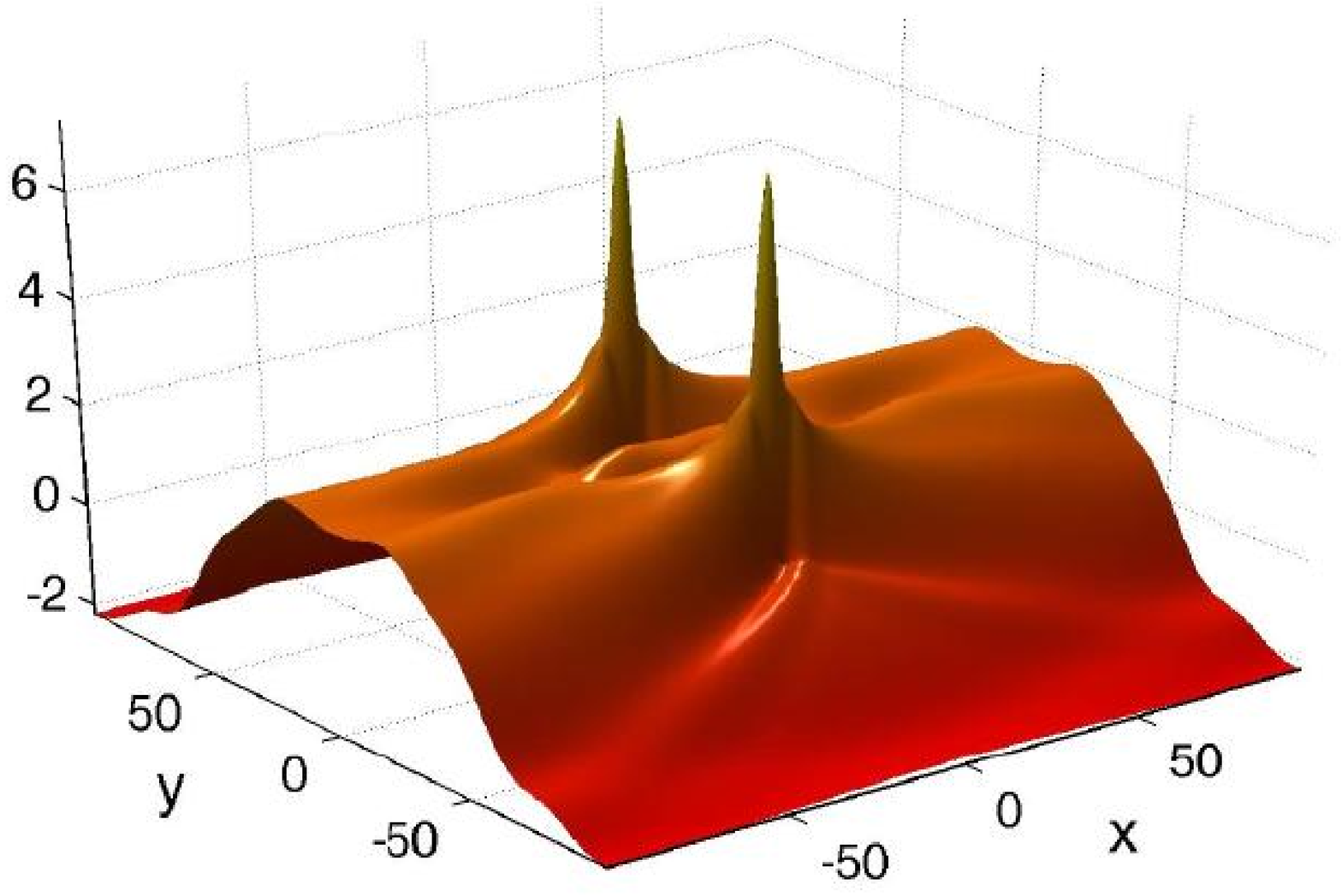}}\\\scalebox{.09}{\includegraphics[width=17in]{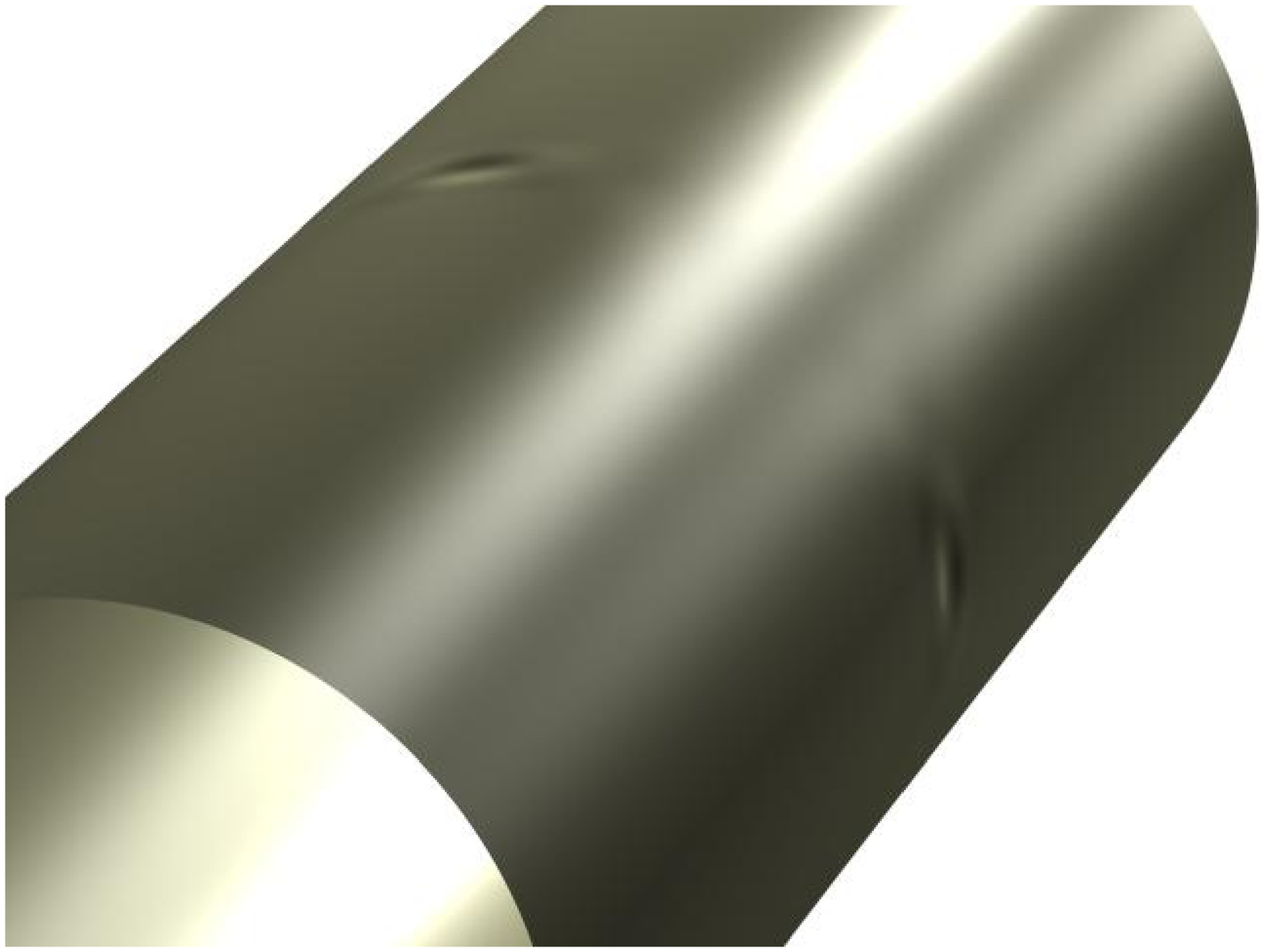}}}}
\put(0,63){\frparbcenter{\scalebox{.1}{\includegraphics[width=16in]{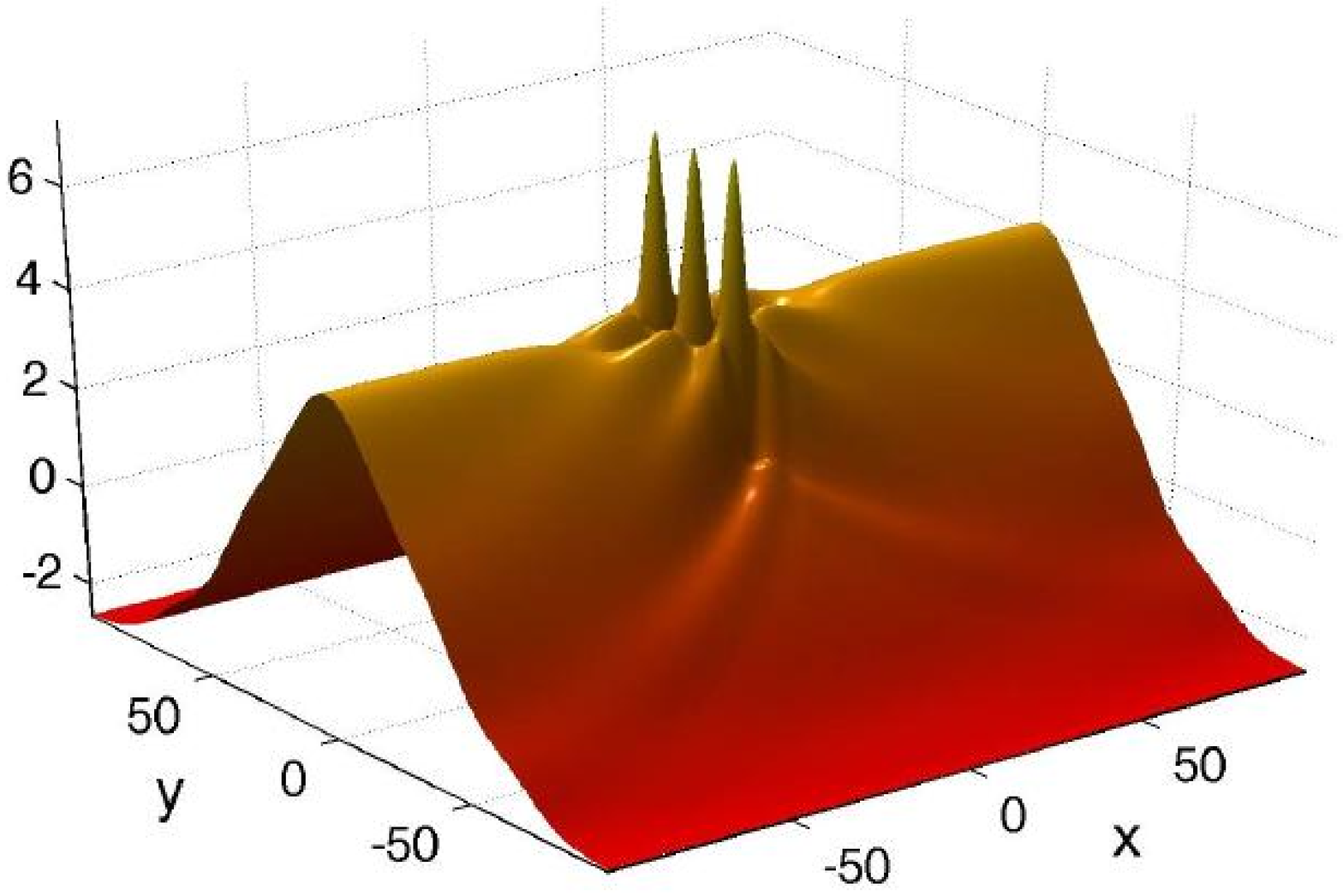}}\\\scalebox{.09}{\includegraphics[width=17in]{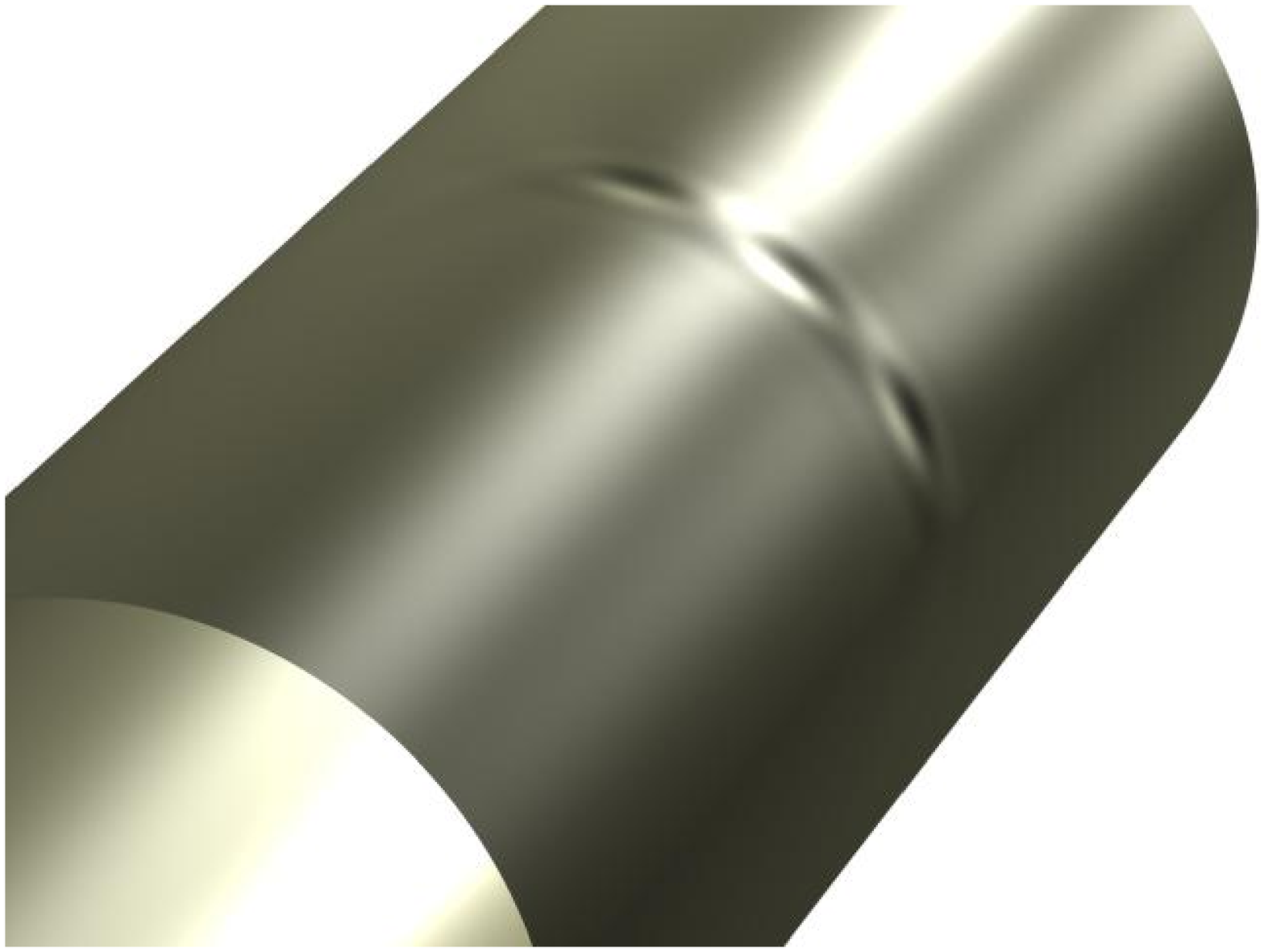}}}}
\put(47,63){\frparbcenter{\scalebox{.1}{\includegraphics[width=16in]{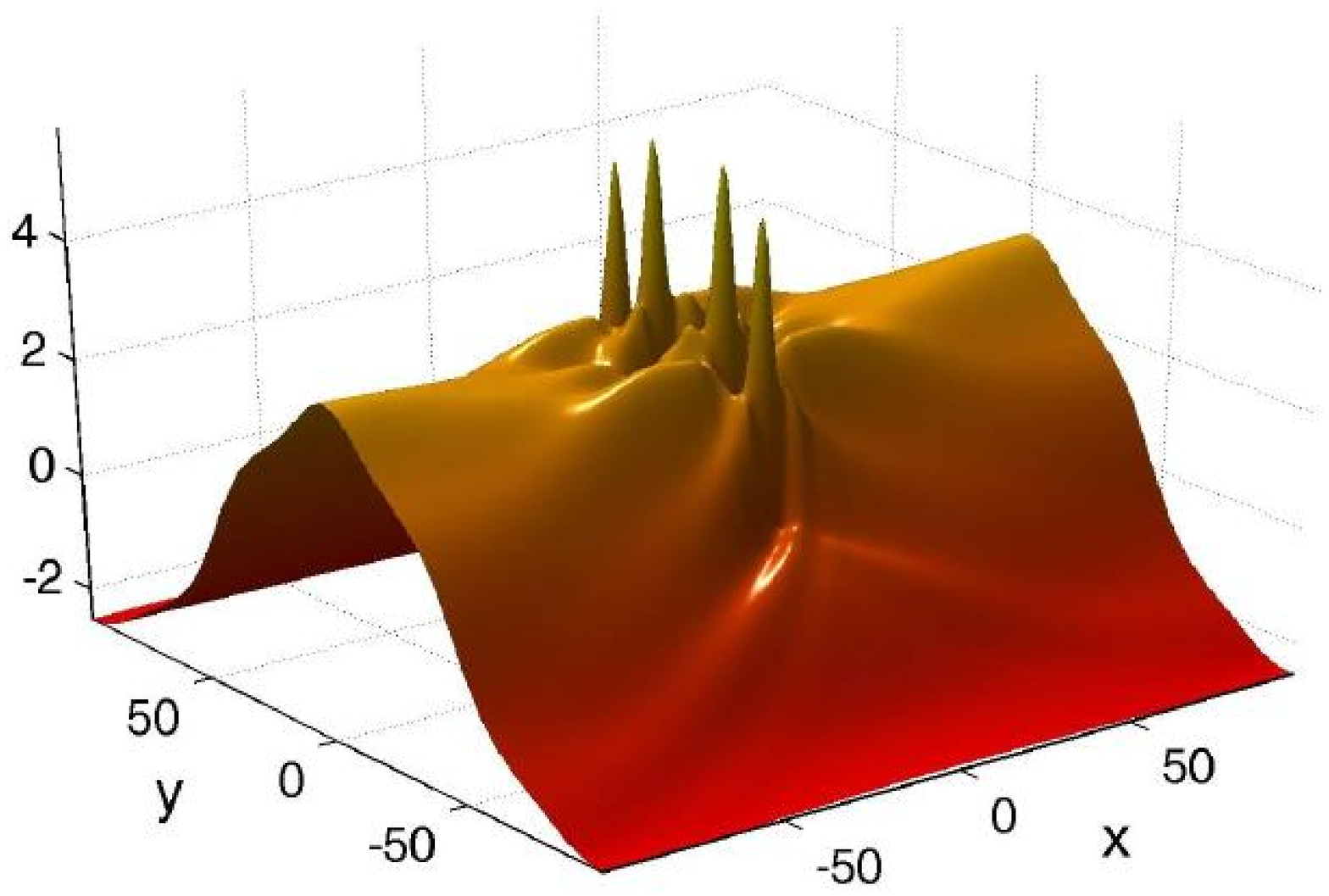}}\\\scalebox{.09}{\includegraphics[width=17in]{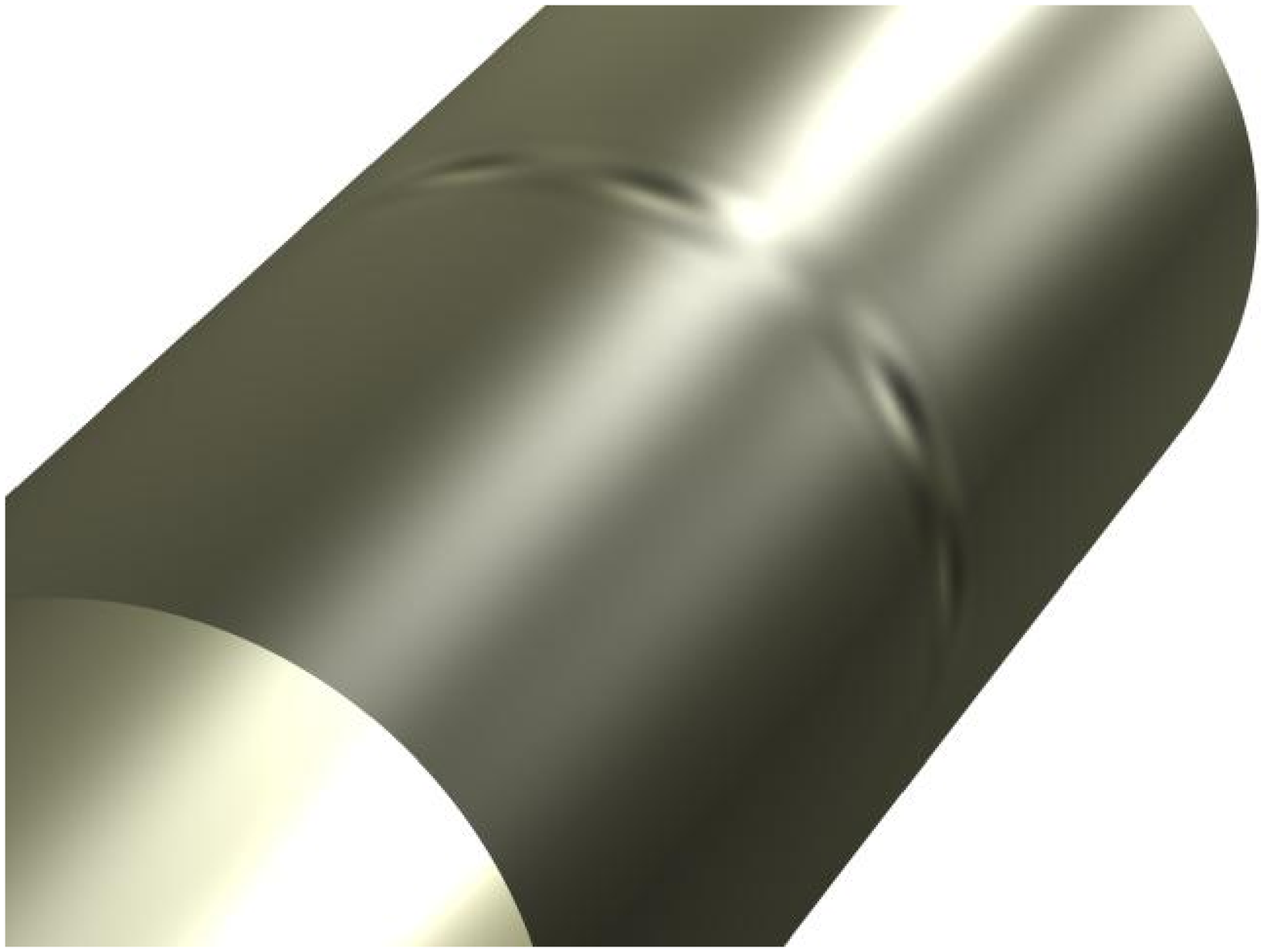}}}}
\put(94,63){\frparbcenter{\scalebox{.1}{\includegraphics[width=16in]{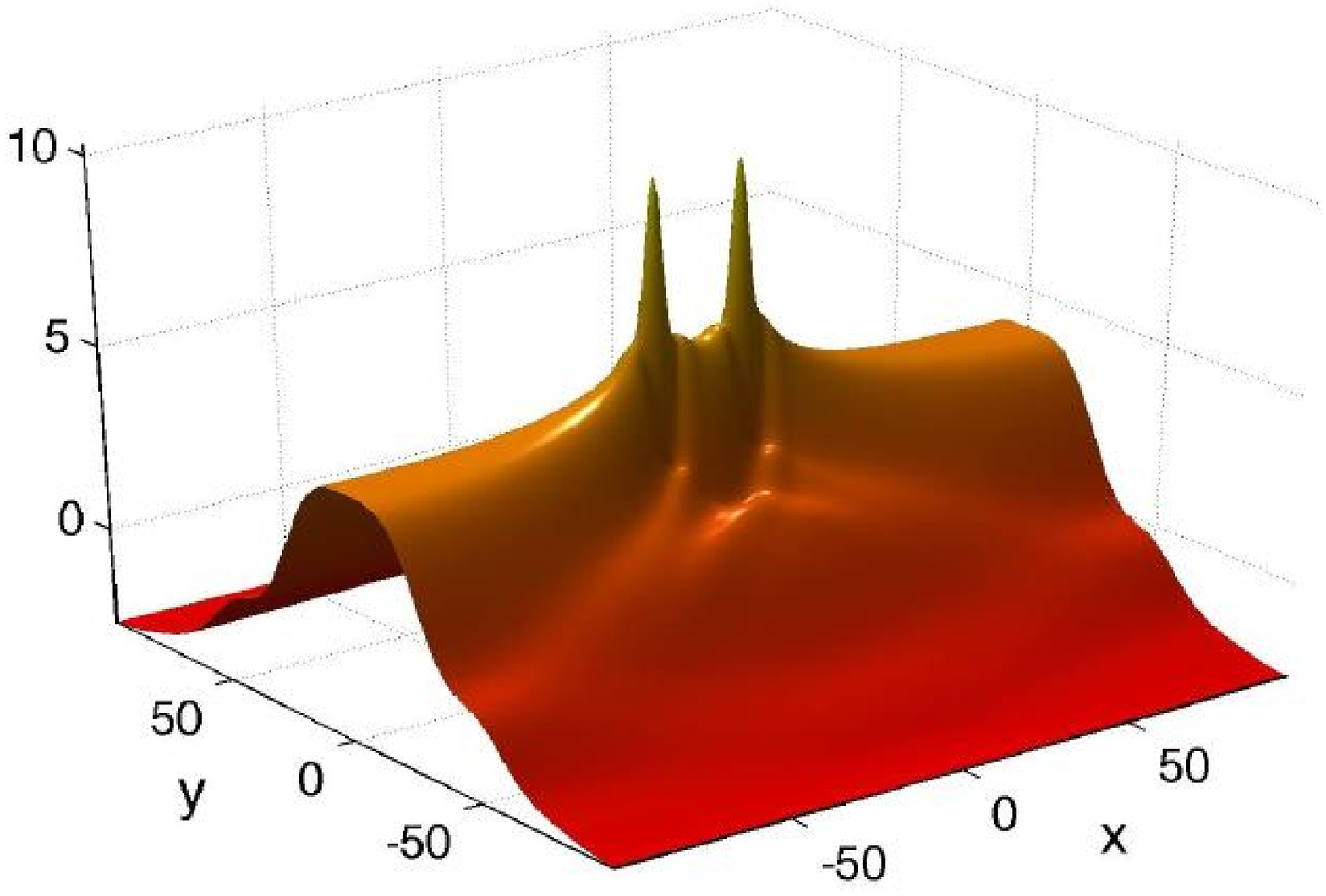}}\\\scalebox{.09}{\includegraphics[width=17in]{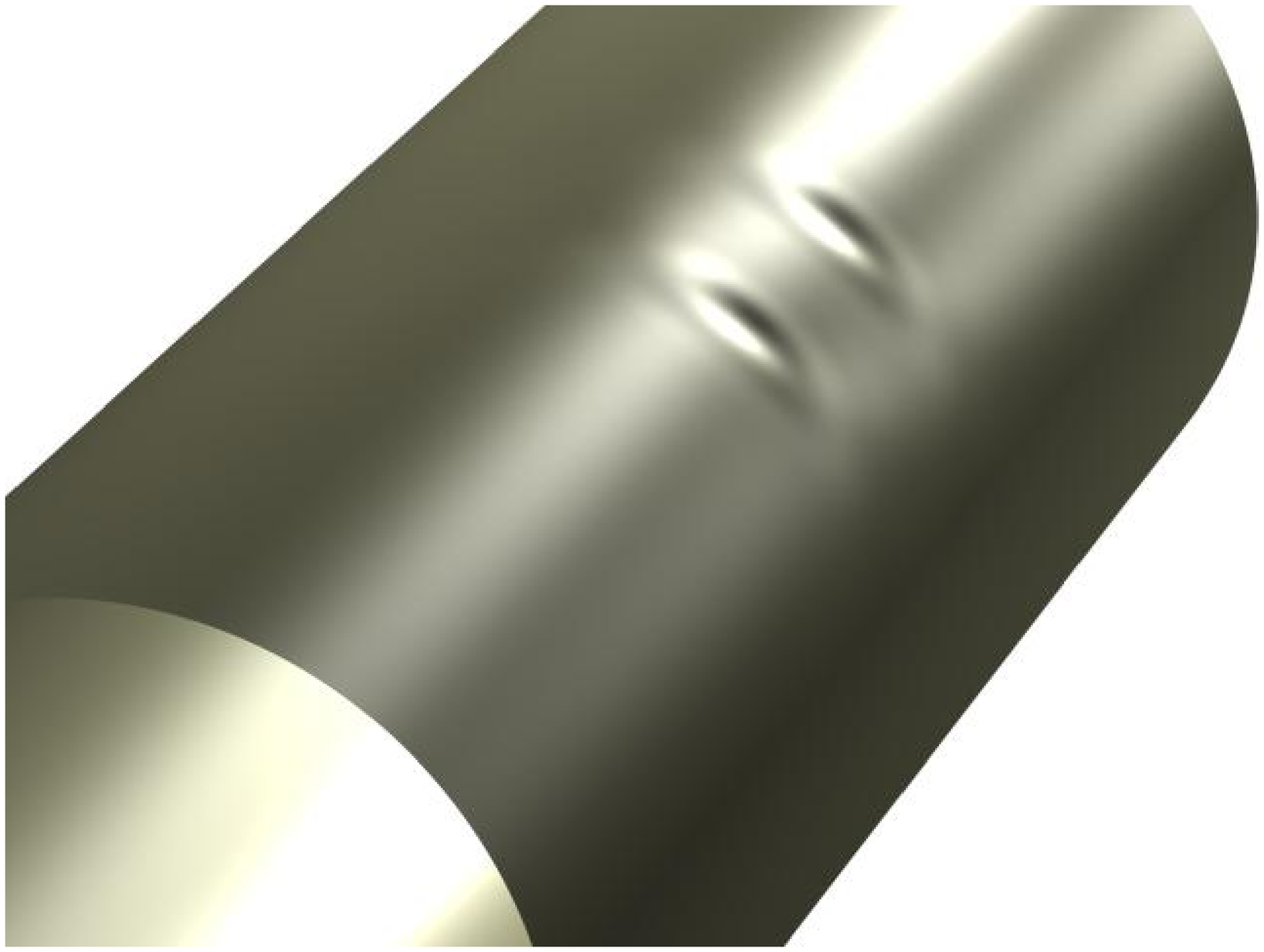}}}}
\put(0,0){\frparbcenter{\scalebox{.1}{\includegraphics[width=16in]{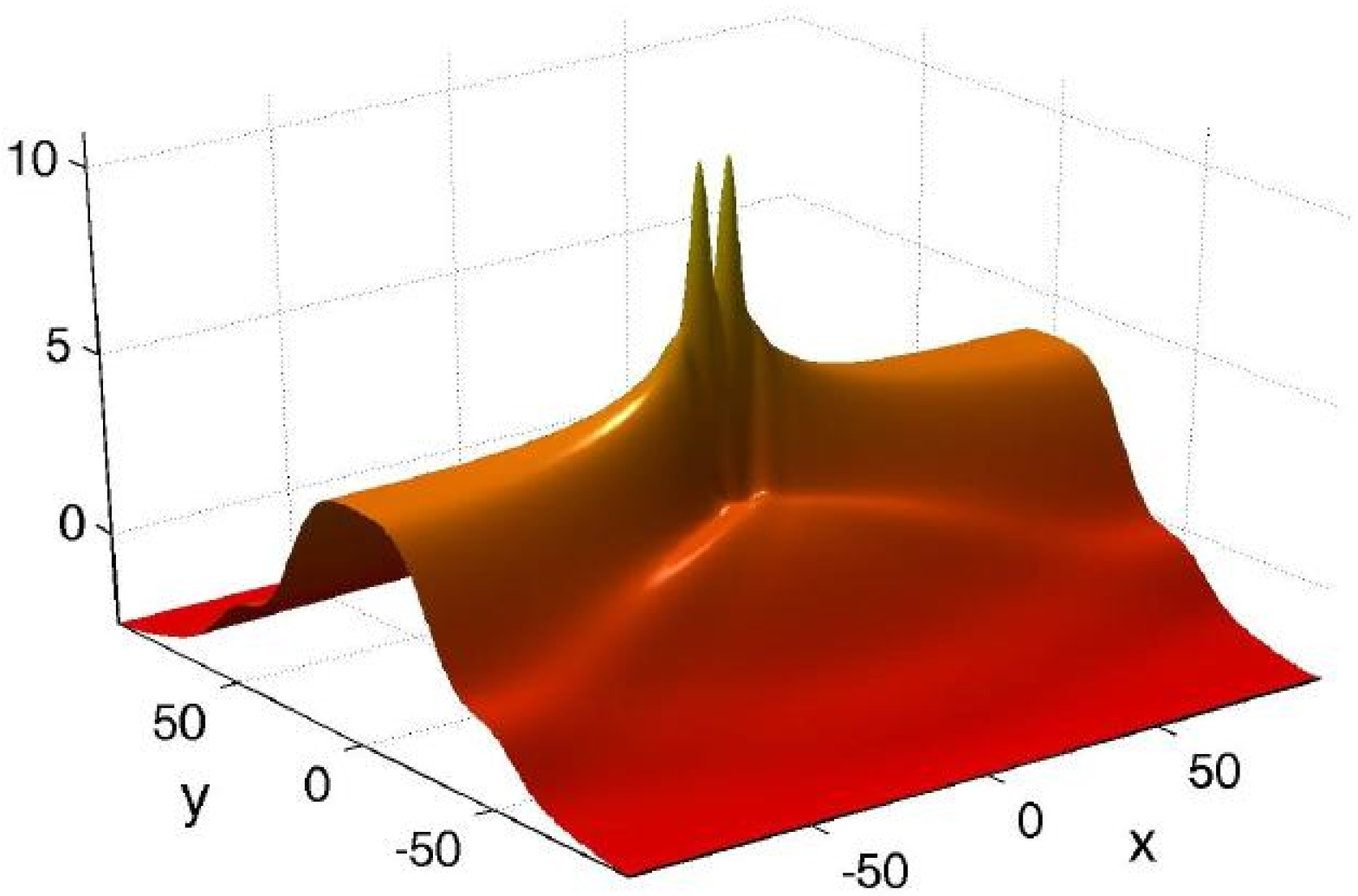}}\\\scalebox{.09}{\includegraphics[width=17in]{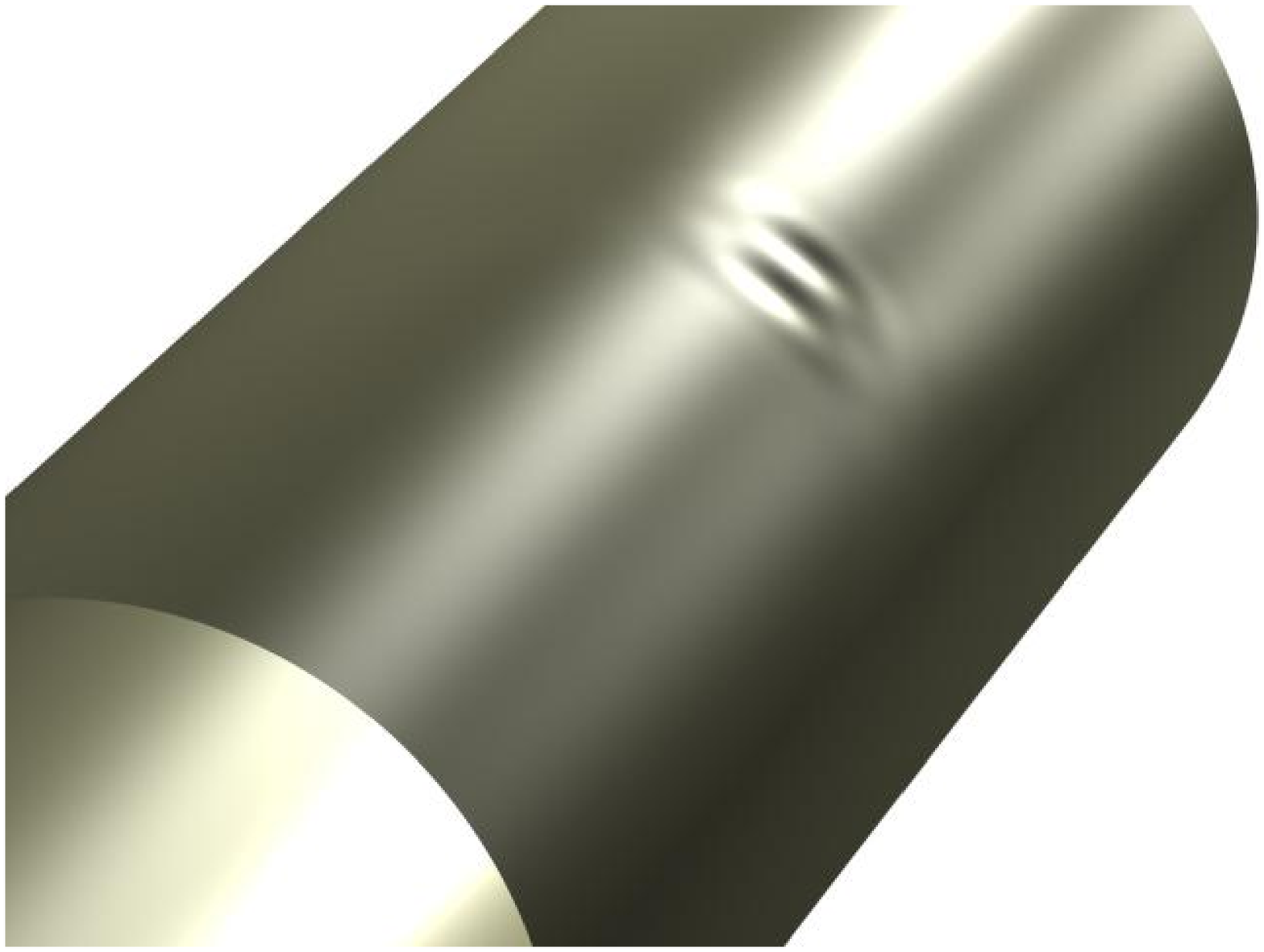}}}}
\put(47,0){\frparbcenter{\scalebox{.1}{\includegraphics[width=16in]{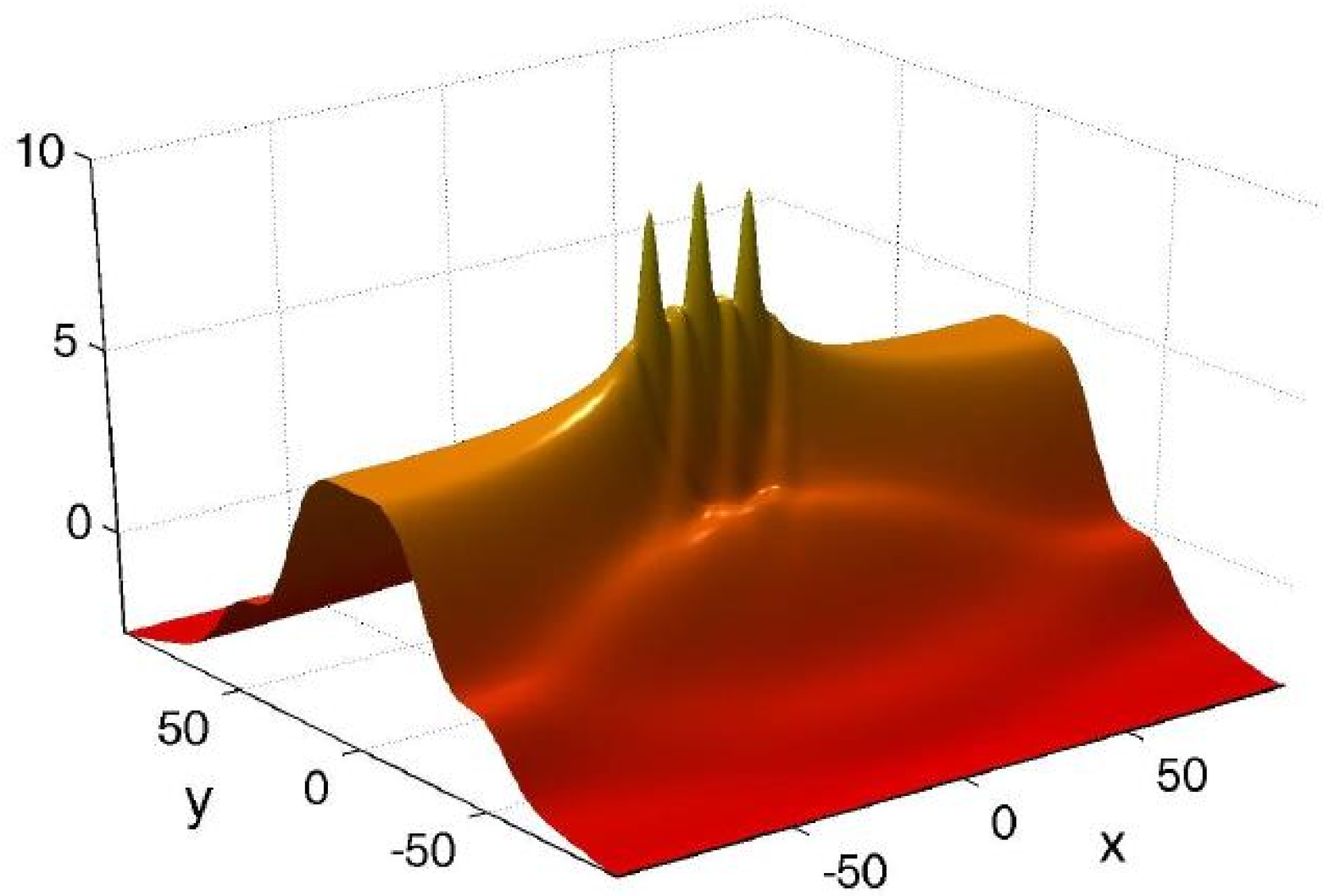}}\\\scalebox{.09}{\includegraphics[width=17in]{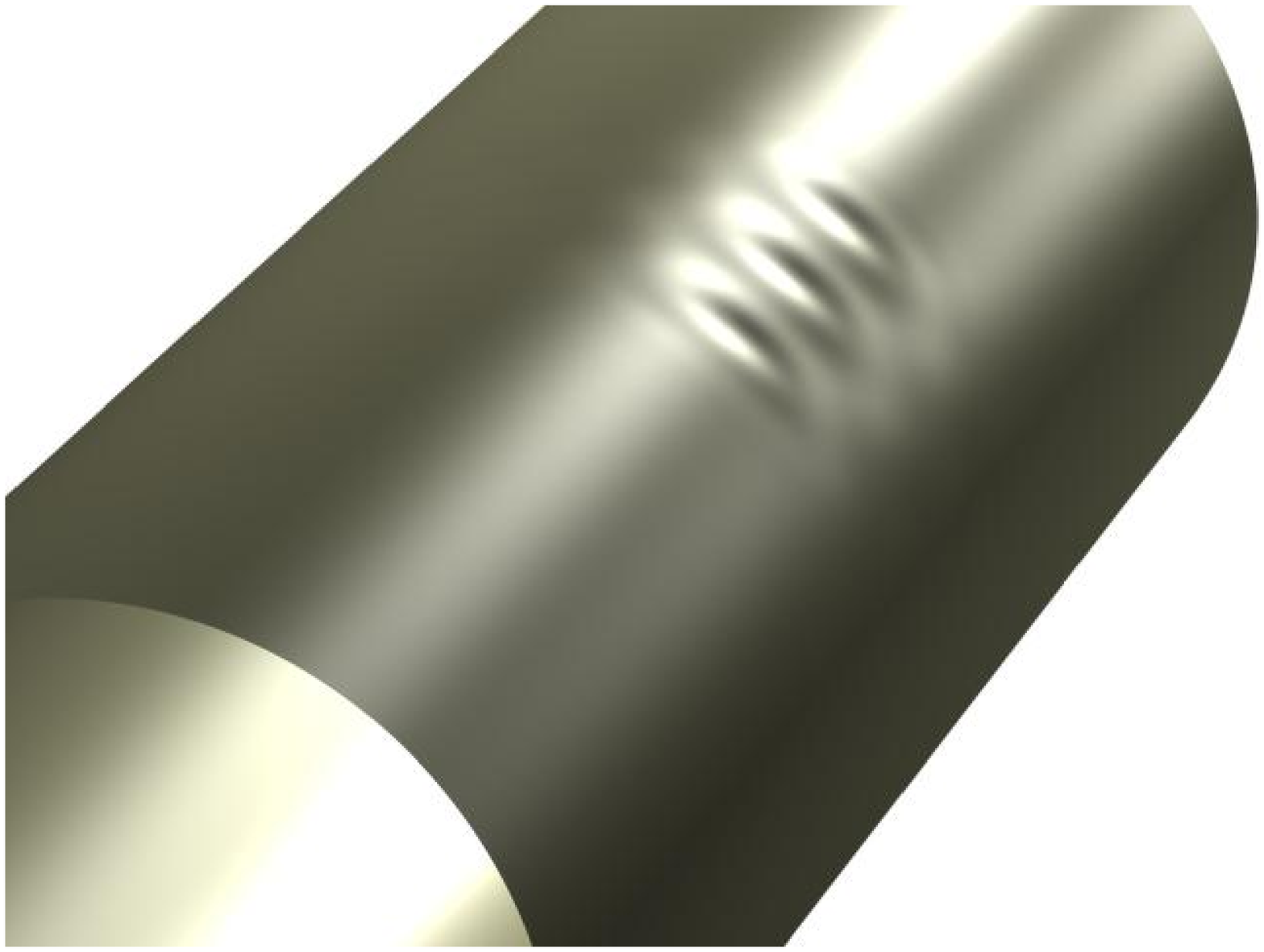}}}}
\put(94,0){\frparbcenter{\scalebox{.1}{\includegraphics[width=16in]{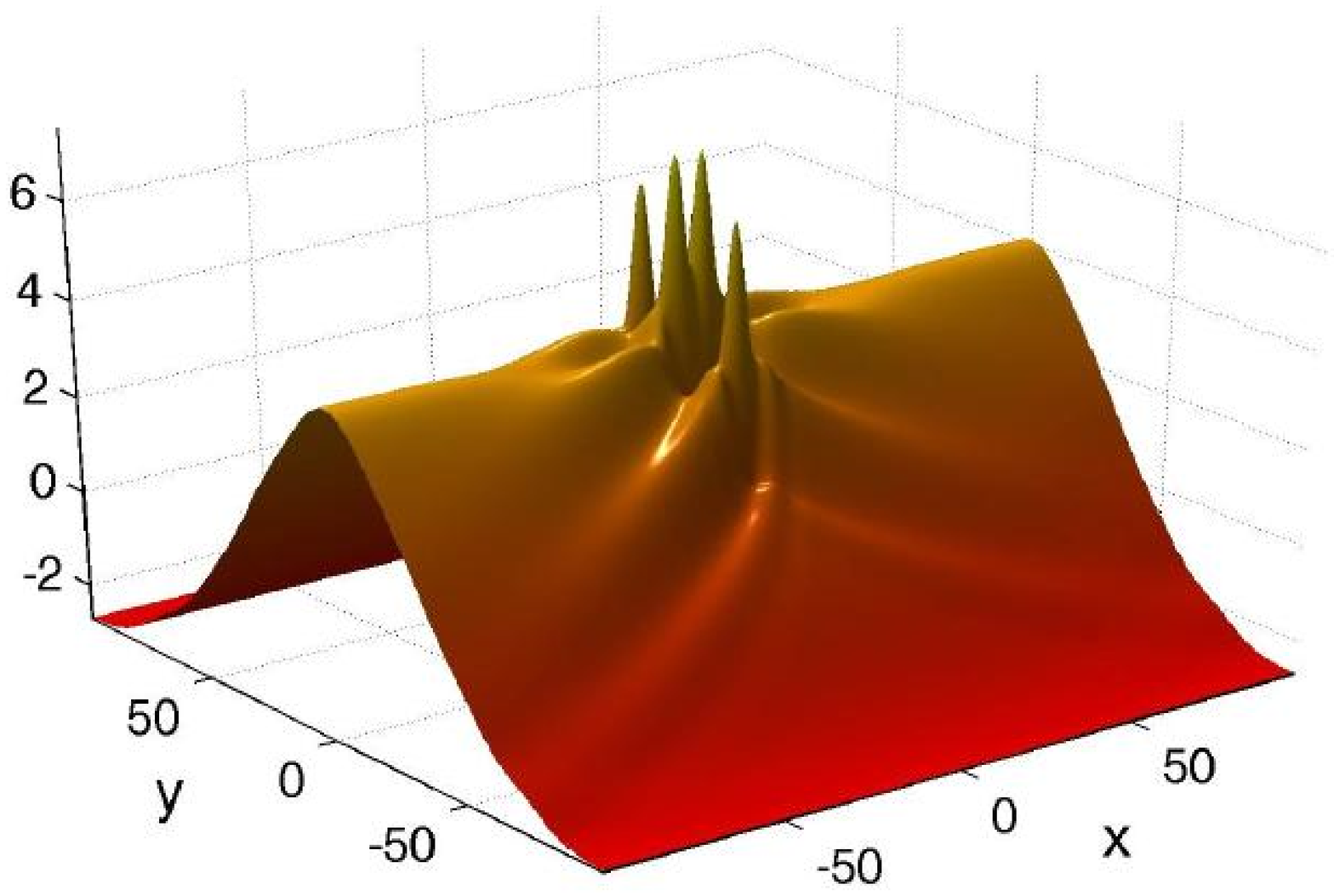}}\\\scalebox{.09}{\includegraphics[width=17in]{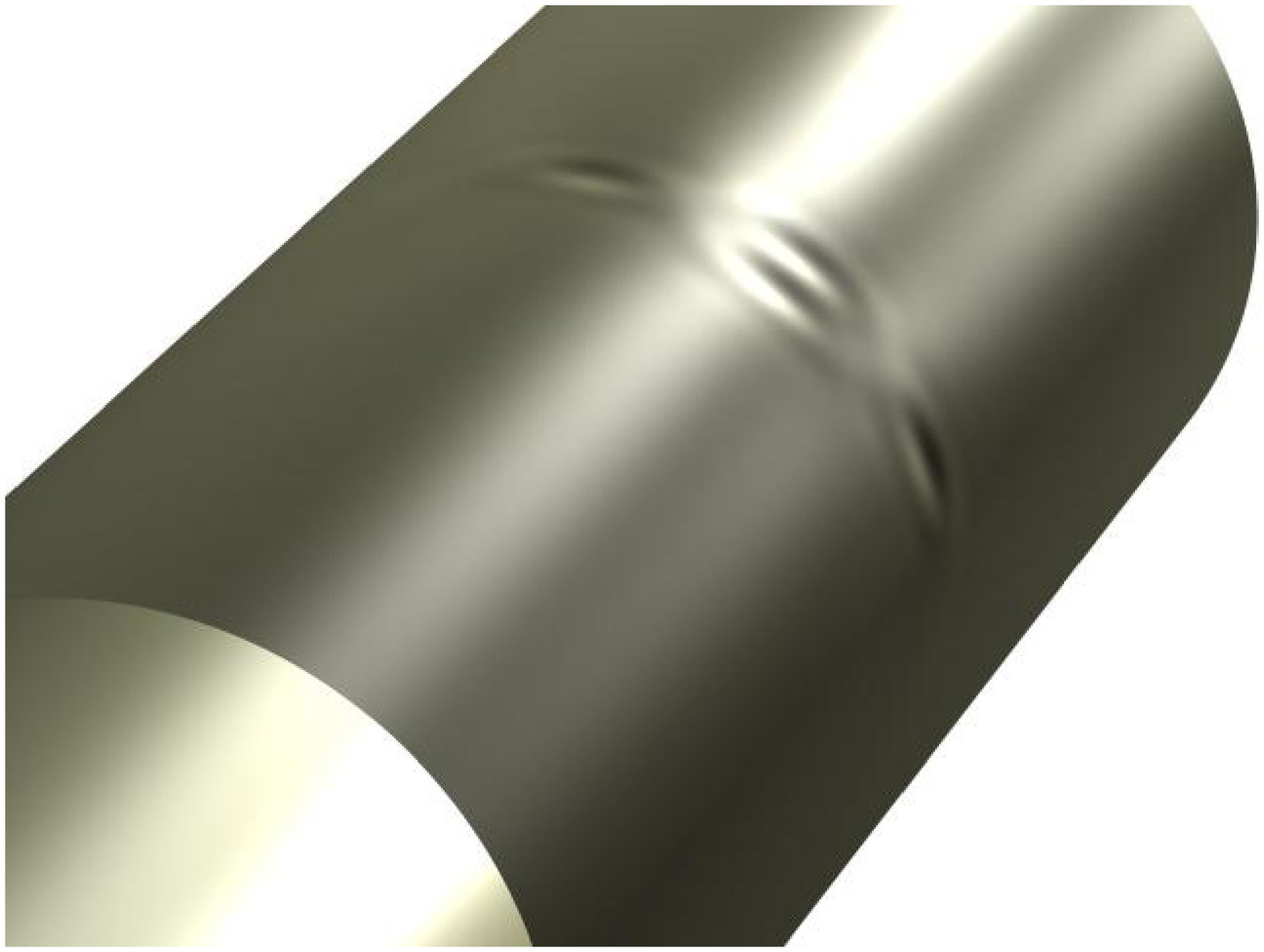}}}}
\color{black}
\put(37,128){(a)}
\put(84,128){(b)}
\put(131,128){(c)}
\put(37,65){(d)}
\put(84,65){(e)}
\put(131,65){(f)}
\put(37,2){(g)}
\put(84,2){(h)}
\put(131,2){(i)}
\end{picture}
\end{center}
\caption{Numerical solutions found using the (constrained) mountain pass
  algorithm, constrained steepest descent method, and the Newton
  algorithm. The figures show both the graph of $w(x,y)$, as well as
  its rendering on a cylinder.}
\label{fig:cylinders}
\end{figure}

\figref{fig:cylinders}(a) shows a numerical solution of
(\ref{eq:main_w2}--\ref{eq:main_phi2}) obtained by the this
algorithm with $\lambda=1.1$. The graph shows the radial displacement
$w$ as a function of $x$ and $y$. Rendered on a cylinder this solution
represents a single dimple as can be seen below the graph.

We restrict our computations to functions that are even about the $x$- and
$y$-axes, 
\emph{i.e.} to the subspace
\begin{displaymath}
  \mathscr{S}=\{\psi\in X :  \psi(x,y)=\psi(-x,y),\; \psi(x,y)=\psi(x,-y)\}\ ,
\end{displaymath}
thus reducing the computational domain $\Omega$ to one quarter, \emph{e.g.}
$(0,\tfrac12\e^{-1/2})\times(0,\tfrac{1}{2}\e^{-1/2})$.
The boundary conditions~(\ref{def:BCw}) then become
\begin{displaymath}
\textstyle
w_x = (\Delta  w)_x = 0 \text{ for }x\in\{0,\tfrac12\e^{-1/2}\}
\quad  \text{and} \quad
w_y = (\Delta  w)_y = 0 \text{ for }y\in\{0,\frac{1}{2}\e^{-1/2}\}.
\end{displaymath}
This symmetry assumption has many numerical advantages, but on the
other hand it \emph{a priori} excludes solutions that do not belong to
$\mathscr{S}$.

For the mountain-pass algorithm we always use the unbuckled state
$w_1=0$ as the first endpoint of the paths. The choice of the second
endpoint $w_2$ has a non-trivial influence on the solution to which
the mountain-pass algorithm converges. Corollary~\ref{cor:periodics}
guarantees the existence of $w_2\in \mathscr S$ with $F_\l(w_2)<0$; in
the numerical implementation, however, we found such a $w_2$ by a steepest
descent method (rather than taking the function constructed in the proof
of Lemma~\ref{th:scaling_we}): starting from a function $w_0$ that has
one peak located in the center of the domain $\Omega$, we solved the
initial value problem
\begin{equation}
\label{ode:finding_w_2}
\frac{d}{d t}w(t) = -\nabla F_\lambda(w(t)), \quad
w(0)=w_0,
\end{equation}
on an interval $(0,T)$ until $F_\lambda(w(T))<0$. We then defined
$w_2=w(T)$.

A different choice of $w_2$ (or, more precisely, of the starting
point $w_0$ of~\pref{ode:finding_w_2}) can lead to a
different solution of the problem, as \figref{fig:cylinders} (b,f)
shows. Here $w_0$ was chosen to have two peaks with centers on the axes
$x=0$ and $y=0$, respectively. The algorithm then converged to a
numerical solution with two dimples in the circumferential and
axial direction, respectively.

Note that the numerical
solution $\wMP$ selected by the mountain-pass algorithm has the mountain-pass property in a
certain neighbourhood only: there exists a ball $B_\rho(\wMP)$ and two
points $\tilde w_1, \tilde w_2 \in B_\rho(\wMP)$ such that
\begin{displaymath}
F_\l(\wMP)=\inf_{\gamma\in\tilde\Gamma} \max_{w\in\gamma}
F_\l(w) > \max\{F_\l(\tilde w_1), F_\l(\tilde w_2)\} \ ,
\end{displaymath}
where $\tilde\Gamma$ is the set of curves in $B_\rho(\wMP)$ connecting
$\tilde w_1$ and $\tilde w_2$. The reason for this is that the
algorithm deforms a certain initial path connecting $w_1$ and $w_2$
which is fixed. In order to recover the global character one would
need to run the algorithm for all possible initial paths.

The rest of the numerical solutions shown in \figref{fig:cylinders}
were obtained under a prescribed value of shortening $S$ by the
constrained steepest descent method and the constrained mountain pass
algorithm \cite{Ho1,HoLoPe2}.

\subsection{Calculation of \boldmath $F_\l(\wMP)$}
In the preceding sections we have showed that 
\begin{enumerate}
\item for a sufficiently large domain $\Omega$, a function 
$w_2$ on $\Omega$ exists with $F_\l(w_2)<0$;
\item for each such function $w_2$ and for almost all $0<\l<2$,
a mountain-pass solution $\wMP=\wMP(\l,\Omega,w_2)$ exists. 
\end{enumerate}
Different end points $w_2$ may give rise to different
mountain-pass points, as we have observed in the numerical
experiments described above. We therefore define
the mountain-pass energy function $V$ on (0,2) by
\begin{equation}
\label{def:VlO}
V(\l,\Omega) := \inf_{w_2} \bigl\{\, F_\l\bigl(\wMP(\l,\Omega,w_2)\bigr): F_\l(w_2)<0\, \bigr\}.
\end{equation}
For a given $\l$, the value of $V(\l,\Omega)$ is the lowest height (or energy level)
at which one may pass from the neighbourhood of the origin to any point with negative
total potential $F_\l$. We now 
derive some of its properties and calculate it numerically.

\begin{lemma}[Properties of $V(\l,\Omega)$]
\begin{enumerate}
\item For sufficiently large $\Omega$ 
there exists $\lambda_0(\Omega)\geq0$ such that 
$V(\lambda,\Omega)<\infty$ for almost all $\l\in(\lambda_0,2)$;
\item $V$ is a decreasing function of $\l$;
\item For sufficiently large $\Omega$, there exists $c(\Omega)>0$ such that
\[
V(\l,\Omega) \leq c(2-\l)^3
\]
for sufficiently small $2-\l>0$.
\end{enumerate}
\end{lemma}
\begin{proof}
Part 1 is a reformulation of the main result of Section~\ref{sec:MP},
making use of Corollary~\ref{cor:periodics}.
For part 2 we remark that for each fixed $w$, $F_\l(w)$ is a decreasing function
of $\l$; the infimum of a set of decreasing functions is again decreasing.

For part 3, let us set 
\[
E(w) = E_2(w) + E_3(w) + E_4(w),
\]
where
\[
E_2(w) := \frac12 \Mod w_X^2 = \frac12 \int_\Omega\bigl(\Delta w^2+\Delta \phi_1^2\bigr), 
\quad E_3(w) := \int_\Omega \Delta\phi_1 \Delta\phi_2,
\quad \text{and}\quad E_4(w) := \frac12 \int_\Omega \Delta\phi_2^2,
\]
where $\phi_1$ and $\phi_2$ are determined from $w$ by~\pref{def:phi_1} 
(see also~\pref{eq:E_develop}).
Note that $E_n$ has homogeneity $n$, \emph{i.e.} $E_n(\mu w) = \mu^n E_n(w)$.

A classical result in the engineering literature of 
cylinder buckling (see \emph{e.g.}~\cite{HLN1})
states that there exists a periodic function $w$ on $\R^2$
such that 
\[
E_2(w)  = 2S(w) \qquad\text{and}\qquad E_3(w)+E_4(w) < 0.
\]
Here and below we consider the integrals that define $E_n(w)$ and $S(w)$ 
as taken over a single period cell.
For sufficiently small $2-\l>0$ the inequality above
gives that $F_\l(w) = (2-\l)S(w)+E_3(w)+E_4(w)$ is negative,
implying that $w$ is an admissible end point $w_2$ for the definition~\pref{def:VlO}
of $V(\l,\Omega)$,
and the connecting line segment $\{\mu w: 0\leq \mu\leq 1\}$ therefore 
an admissible curve in $\Gamma$. Consequently
\[
V(\l) \leq \sup_{0\leq\mu\leq 1} F_\l(\mu w) 
= \sup_{0\leq\mu\leq 1} \mu^2 (2-\l)S(w) + \mu^3 E_3(w) + \mu^4 E_4(w)
\]
The supremum on the right-hand side is obtained at
\[
\mu = \frac {3\mod{E_3(w)}}{8E_4(w)}
  \left\{1 - \sqrt{1-\frac{32(2-\l)S(w)E_4(w)}{9E_3(w)^2}}\right\}
  = \frac{2S(w)}{3\mod{E_3(w)}}\,(2-\l) + o(1) \quad \text{as }\l\to2,
\]
implying that the claim holds for periodic functions.  The
generalization to non-periodic functions on large domains $\Omega$
(\emph{i.e.} for small $\e$) is made by filling the domain with a large
number of periodic cells of the function $w$, and connecting the
function smoothly to the boundary of $\Omega$.
\end{proof}

\begin{figure}[htbp]
  \begin{center}
    \setlength{\unitlength}{1mm}
    \begin{picture}(140,46)
      \put(4,1){\scalebox{.6}{\includegraphics{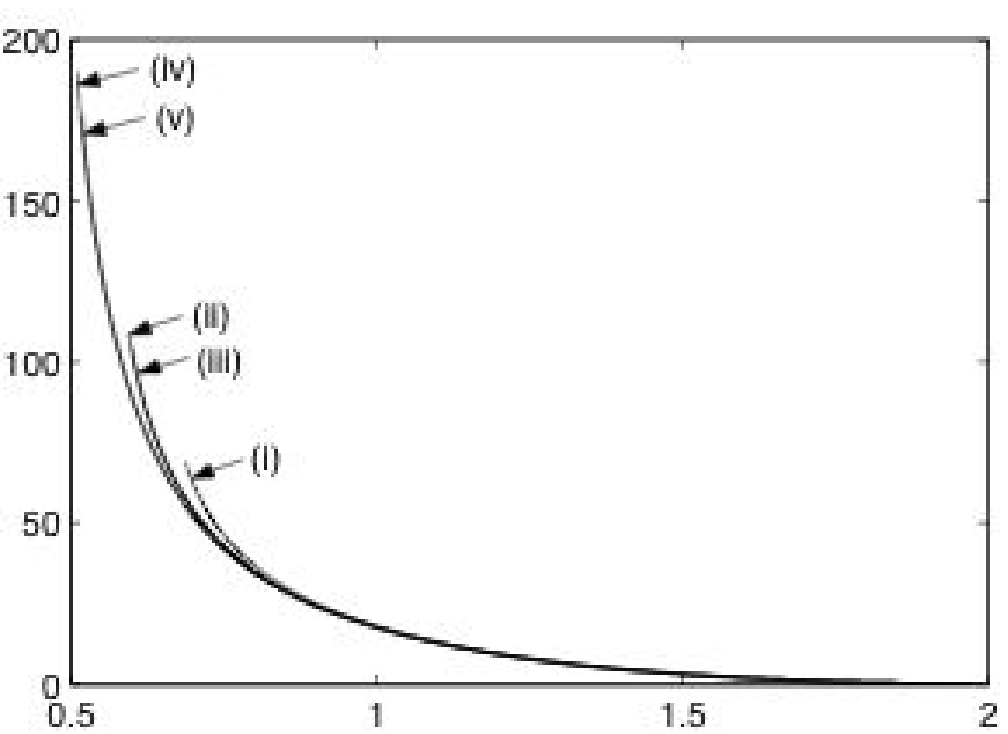}}}
      \put(75.5,1){\scalebox{.6}{\includegraphics{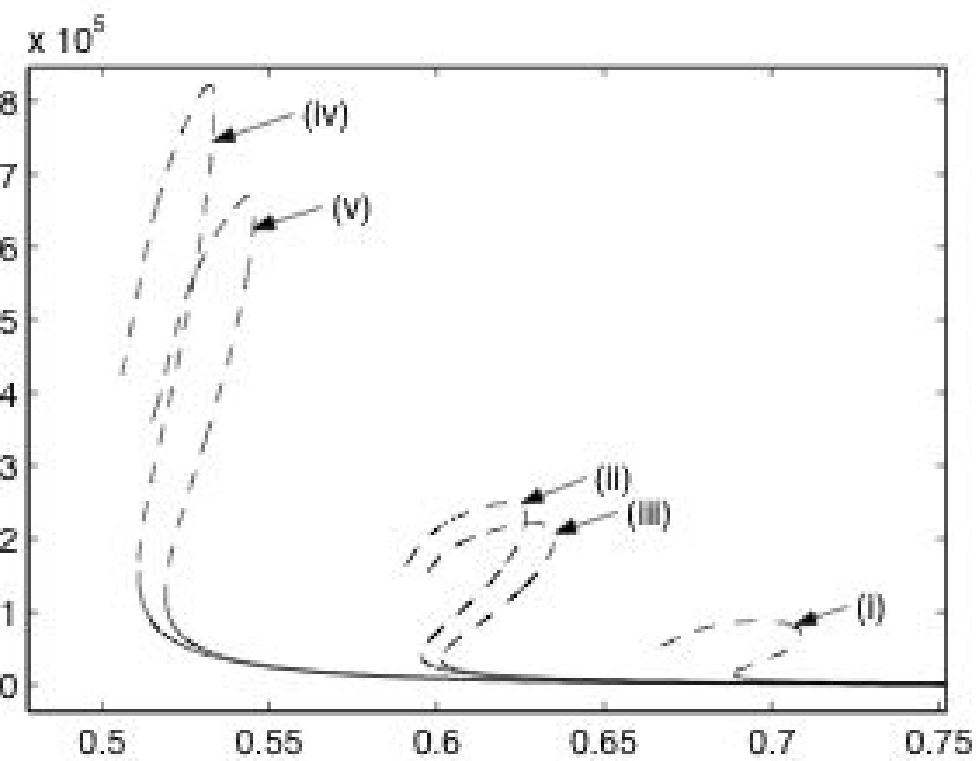}}}
      \put(29.5,23){\framebox(33,19){\parbox{35mm}{\scriptsize\centering
            $\Omega=(-a,a)\times(-b,b)$:\\
            \begin{tabular}{c@{ $\ldots$ }l}
              (i) & $a=b=50$ \\
              (ii) & $a=b=100$ \\
              (iii) & $a=100, b=200$ \\
              (iv) & $a=200, b=100$ \\
              (v) & $a=b=200$
            \end{tabular}
          }}}
      \put(1,25){\rotatebox{90}{\hbox{\small$V(\l)$}}}
      \put(57,0){\small$\l$}
      \put(70,22){\rotatebox{90}{\hbox{\small$\|\wMP(\l)\|^2$}}}
      \put(126,0){\small$\l$}
    \end{picture}
  \end{center}
  \caption{Left: the mountain-pass energy $V(\l)$ found numerically for
    various domains $\Omega$. Right: solid line shows the same
    computation as on the left, just plotted for the norm of
    $\wMP(\l)$ squared. The dashed curve was obtained by continuation
    of the solid curve; the solutions on the dashed curve do not represent
    mountain-pass points any more.}
  \label{fig:Vlambda}
\end{figure}

\figref{fig:Vlambda} shows graphs of the mountain-pass energy
$V(\l,\Omega)$ computed for various sizes of domain $\Omega$. For
each domain, the mountain-pass algorithm was employed to compute
$\wMP$ for several values of $\lambda$. These mountain-pass solutions
were then continued in $\lambda$ using numerical continuation.

\subsection{Influence of the domain}
\label{subsec:domain_influence}

The localized nature of the solutions calculated above suggests that they should
be independent of domain size, in the sense that for a sequence of domains 
of increasing size the solutions converge (for instance pointwise on
compact subsets). Such a convergence would also imply convergence of 
the associated energy levels.
Similarly, we would expect that the aspect ratio of the domain is of little importance
in the limit of large domains. 

We have tested these hypotheses by computing mountain-pass solutions 
on domains of different sizes and aspect ratios. Generally 
solutions on different domains compare well; the maximal difference in the
second derivatives of $w$ is two or three orders of magnitude smaller
than the supremum norm of the same derivative (the
details of this comparison are given in~\cite{HoLoPe2}).

Here we only include a calculation of the mountain-pass energy level
$V(\l,\Omega)$ for domains $\Omega$ of different aspect ratio and size
(\figref{fig:Vlambda}).

The comparison of solutions computed on different domains and their respective
energies suggests that for each $\l$ we are indeed dealing with a single, localized
function defined on $\R^2$, of which our computed solutions are
finite-domain adaptations. In the rest of this paper we adopt this point
of view, and consequently we will write $V(\l)$ instead of $V(\l,\Omega)$.

\begin{figure}[h]
\begin{center}
  \setlength{\unitlength}{1mm}
  \begin{picture}(130,46)
    \color[rgb]{.5,.5,.5}
    \put(0,1){\fbox{\scalebox{.13}{\includegraphics[width=17in]{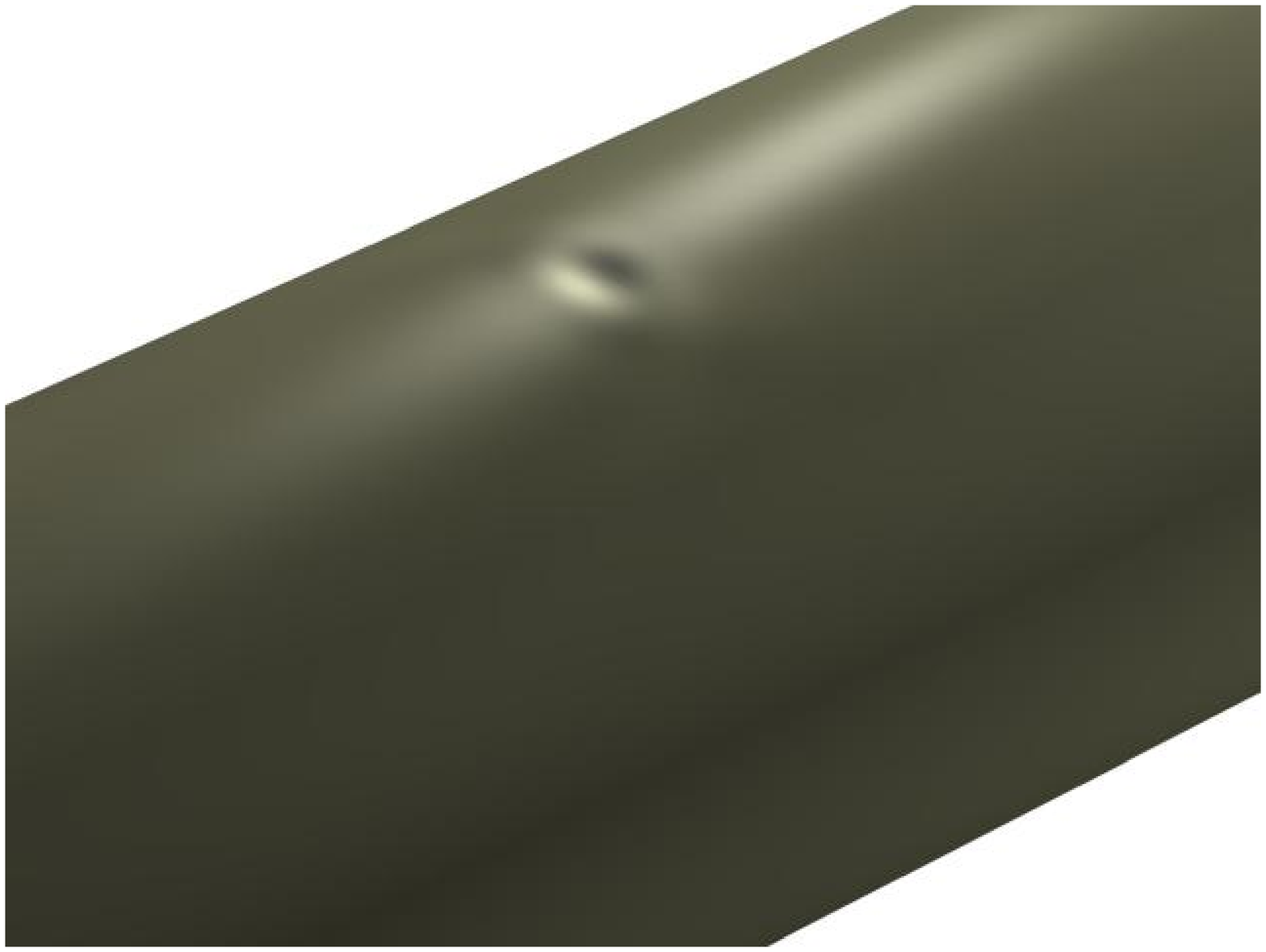}}}}
    \put(69,1){\fbox{\scalebox{.13}{\includegraphics[width=17in]{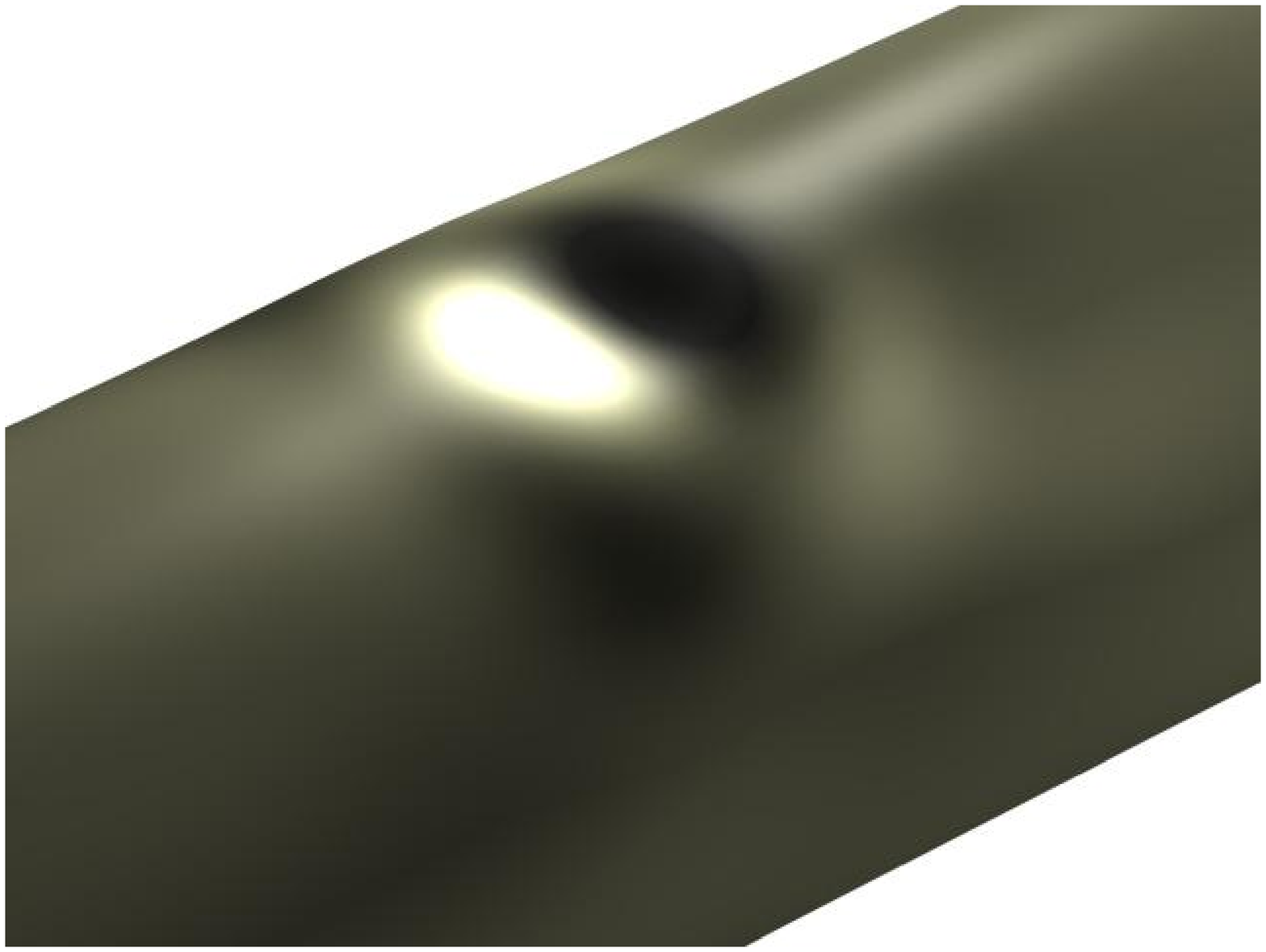}}}}
    \color{black}
    \put(55,2){(a)}
    \put(124,2){(b)}
  \end{picture}
\end{center}
\caption{The numerical mountain-pass solution $\wMP$ of the scaled equations
  (\ref{eq:main_w2}-\ref{eq:main_phi2}) for a given value of $\lambda$
  rendered on two cylinders of the same radius $R$ 
  but a different thickness $t$: (a) $t/R=0.003$, (b) $t/R = 0.04$.}
\label{fig:dimple_detail}
\end{figure}

A consequence of this point of view is that dimples
in cylinders with different geometric parameters
are mapped to the same rescaled solution, or equivalently, that the 
same single-dimple solution of (\ref{eq:main_w2}-\ref{eq:main_phi2})
corresponds to differently-sized dimples on an actual cylinder, 
as a function of the parameters (\figref{fig:dimple_detail}).

\section{Interpretation: imperfection-sensitivity}
\label{sec:discussion}

We now turn to the relevance of the mountain pass in the context of a
loading problem. This relevance can best be understood in the context
of imperfections in the loading conditions (rather than geometric
imperfections) such as in the case of a (small) lateral loading.

Under a small lateral load 
an equilibrium $w_0$, a local minimum of the functional $F_\lambda$,
may be perturbed into an equilibrium
$\tilde w_0$ of a perturbed functional $\tilde F_\lambda$. Since  $w_0$ is a
local minimum, $F_\lambda(\tilde w_0) > F_\lambda(w_0)$, \emph{i.e.} 
with respect to the unperturbed system $\tilde w_0$
has a higher total potential than $w_0$. The
level of $F_\lambda$ that is reached is a measure of the magnitude of
the imperfection---a different measure than is commonly used, but one that
has distinct advantages.

\begin{figure}[ht]
\centerline{\psfig{figure=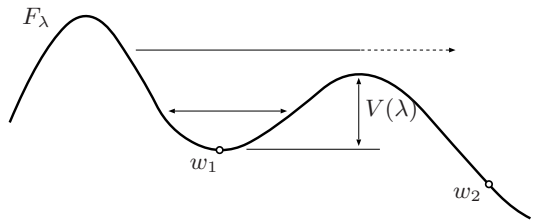,height=3cm}}
\caption{In order to leave the basin of attraction of $w_1$, the surplus
energy should exceed $V(\lambda)$}
\label{fig:trapped}
\end{figure}

By definition the number $V(\lambda)$ is 
the lowest energy level at which it is possible to move between
the basins of attraction of $w_1$ and $w_2$ (\figref{fig:trapped}).
If the loading imperfection is interpreted, as above, as a mechanism capable of
maintaining the system at a higher energy level than that of the neighbouring
fundamental minimizer, then the number $V(\lambda)$ is critical: as long as the imperfection is
so small that the energy is never raised by more than $V(\lambda)$, the new stationary point
will be part of the same basin of attraction as $w_1$. 
For larger imperfections, however, it becomes possible to leave the fundamental
basin of attraction, resulting in a large jump in state space.

\medskip

This line of reasoning provides a natural measure of the size of imperfections,
namely the maximal increase in energy (in the perfect structure) that an imperfection
can achieve. It also provides a natural measure of the stability of the unbuckled
state, since a higher mountain-pass energy level implies a larger class of
loading imperfections under which the state remains in the fundamental
basin of attraction.
This observation allows us to connect systems with different geometrical
characteristics, and compare their relative sensitivity to imperfections.

\subsection{Calibrating the mountain-pass energy}

Comparing cylinders of varying geometry requires
a common measure of imperfection sensitivity. It is not \emph{a priori}
clear which measure to take; \emph{e.g.}, one might consider either the mountain-pass energy 
itself or the average spatial density of this energy, which will 
result in different comparisons for cylinders of different wall volume.
Here we choose to rescale the mountain-pass energy level by the other 
energy level present in the loaded cylinder: 
the energy that is stored in homogeneous compression 
of the unbuckled shell.

\begin{figure}[ht]
\centering
\centerline{\psfig{height=7cm,figure=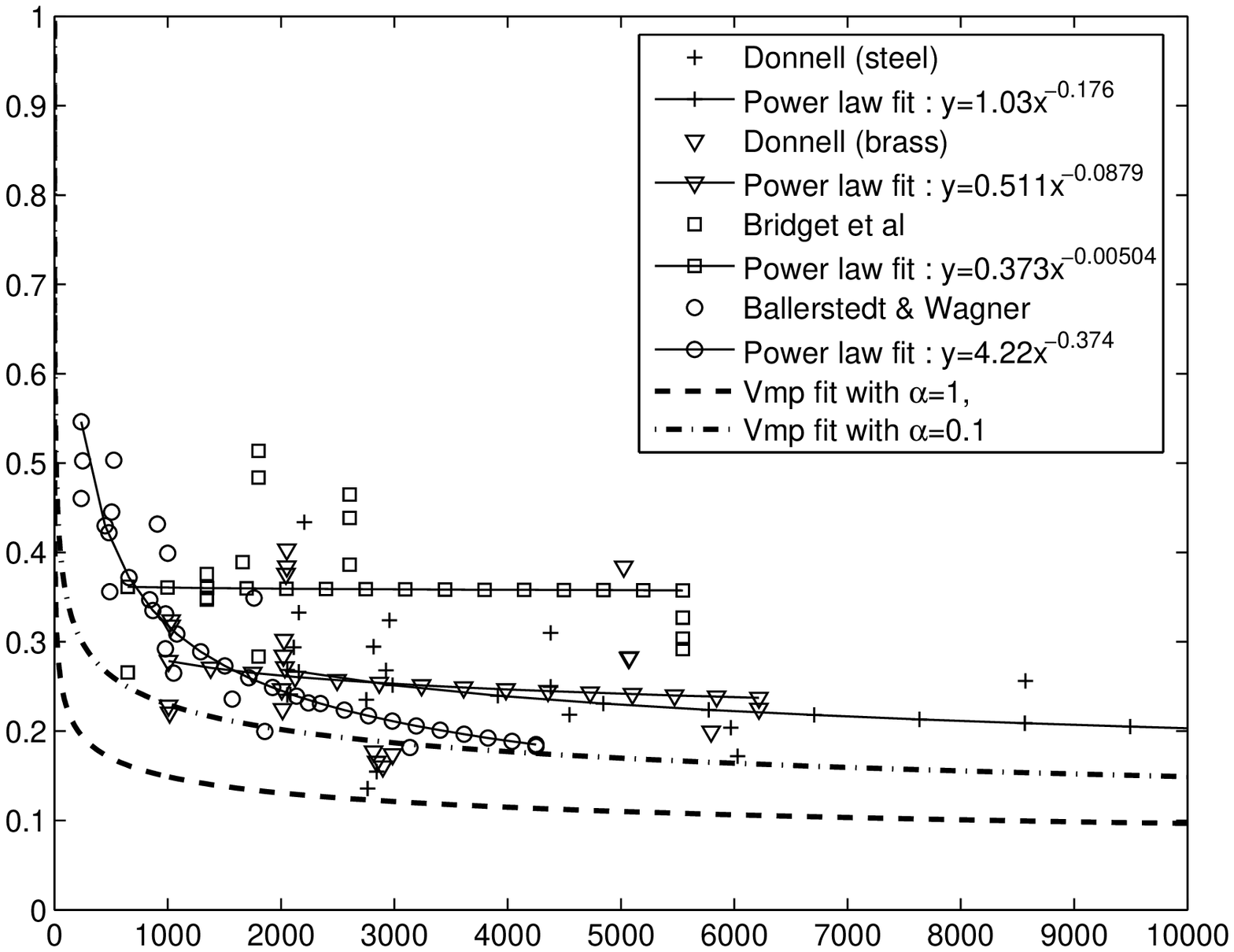}}
\caption{Same data as \figref{fig:experiments}; added are two
curves of constant $\alpha=1$, $\alpha=0.1$ where $\alpha$ is given
in~\pref{def:alpha}. Note that the load at which the mountain-pass
energy equals  the stored energy in the prebuckled cylinder
($\alpha=1$), appears to be a lower bound to the data. 
} 
\label{fig:Lt}
\end{figure}

This calculation can be done in two slightly different ways. The first and most
straightforward is to rescale the dimensional mountain-pass energy 
(see~\pref{def:rescaled_energies} and~\pref{scaling:epsilon})
\[
64\pi^6 EtR^2 \e^3 V(\lambda) = \frac {Et^4}{8(3(1-\nu^2))^{3/2}R}\; V(\lambda),
\]
by the elastic strain energy stored in the full length of the compressed cylinder
of length $L$,
\[
\frac L{4\pi ERt}\; P^2 = \frac{\pi t^3EL}{12(1-\nu^2)R}\;  \lambda ^2,
\]
to give an energy ratio, or a rescaled mountain-pass energy level,
\begin{equation}
\label{def:alpha}
\alpha = \frac{1}{2\pi\sqrt{3(1-\nu^2)}}\; \frac tL \; \frac{V(\lambda)}{\lambda^2}.
\end{equation}
From this expression and the calculation shown in 
\figref{fig:Vlambda} curves may be drawn in
a plot of load versus the ratio $L/t$ (see \figref{fig:Lt}). Note that
to obtain this figure from \figref{fig:Vlambda} the curve $V(\lambda)$
was fitted to extend the range of $\lambda$. 
This figure shows two remarkable
features:
\begin{enumerate}
\item The general trend of the constant-$\alpha$ curves is very similar
to the trend of the experimental data;
\item The $\alpha=1$ curve, which indicates the load 
at which the mountain-pass energy equals  the stored energy in the 
prebuckled cylinder, appears to be a lower bound to the data.
\end{enumerate}

\medskip

One may also consider an alternative way of rescaling energy. The cylinder 
is a long structure, and it is not clear to which extent the length of
the structure is relevant for the imperfection sensitivity. It may be
reasonable to compare the energy of the mountain pass with the stored
energy contained 
in a representative section of the cylinder; the radius $R$ provides a 
natural length scale for such a representative section. 

\begin{figure}[ht]
\centering
\centerline{\psfig{height=7cm,figure=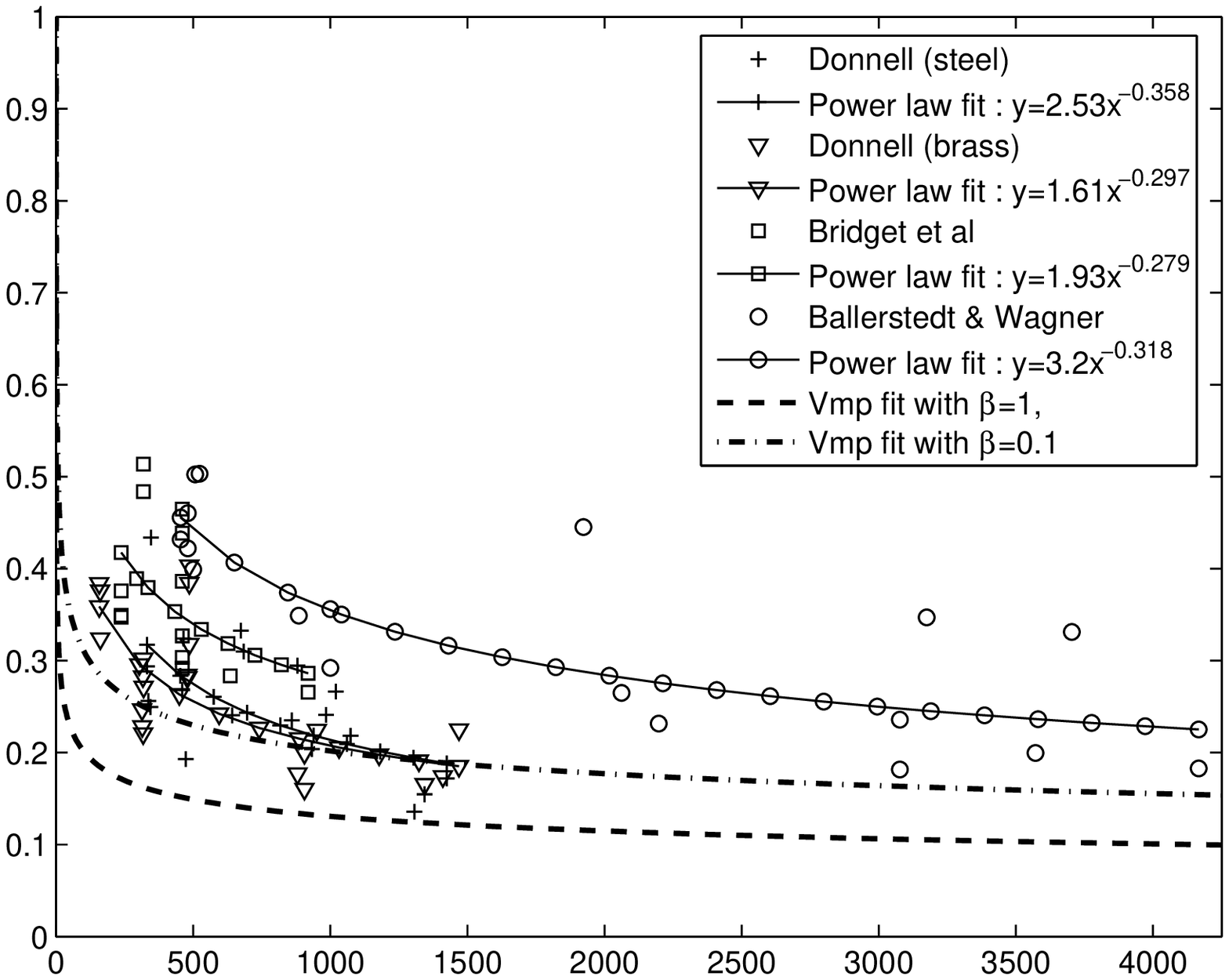}}
\caption{Experimental data and fit to $V_{mp}$ with $\beta=1$ and
  $\beta=0.1$ in ~\pref{def:beta}. Again the load at which the
  mountain-pass energy equals  the stored energy in a representative
  portion of the prebuckled cylinder
  ($\beta=1$) appears to be a lower bound to the data.
}
\label{fig:Rt}
\end{figure}

Similar to \figref{fig:Lt}, in \figref{fig:Rt} we present
curves of constant $\beta$, where $\beta$ is the ratio of 
mountain-pass energy to stored energy in
a section of length $2\pi R$:
\begin{equation}
\beta = \frac{1}{4\pi^2\sqrt{3(1-\nu^2)}}\; \frac tR \; \frac{V(\lambda)}{\lambda^2}.
\label{def:beta}
\end{equation}
Once again to obtain \figref{fig:Rt} we fitted $V(\lambda)$ from 
\figref{fig:Vlambda} to extend the range of $\lambda$.

\section{Discussion and conclusions}
\label{sec:conclusions}

The mathematical results and their interpretation in the context of
a loading problem have brought a number of new and improved insights.

\subsection{The cylinder has doubly-localized solutions}
The subcritical nature of the bifurcation in \figref{fig:bif}
strongly suggests that equilibria exist with deformation localized
to a small portion of the cylinder length. In~\cite{HLN1,LordChampneysHunt97,LCH2} 
such solutions are indeed calculated numerically and investigated analytically; 
these solutions are periodic
around the cylinder and have exponential decay in the axial direction. 

The localization in the axial direction demonstrated by these solutions
is consistent with results on simpler systems such
as the laterally supported strut~\cite{HPCWWBL,Peletier01}.
The behaviour of the cylinder in the tangential direction is not as well 
understood. The lack of localization in the simply-supported
flat plate~\cite{ArnBC2} suggests 
that the cylinder should also prefer tangentially delocalized  solutions,
as do most of the experiments.
The single- and multiple-dimple solutions of this paper, however,
clearly  demonstrate that doubly-localized solutions
do exist, and that some of these can be stable under constrained shortening.

\subsection{The mountain pass is a single-dimple solution}

The fact that the mountain-pass solution exists follows essentially
from two features, the local minimality of the unbuckled state and the existence
of a large-deflection state of lower energy.
The former is a simple consequence\footnote{On a finite domain
this consequence is indeed simple; on an infinite domain it appears
that besides the third-order term also the fourth-order term in the
energy has to be taken into account, as remarked in Section~\ref{sec:MP}}
 of the subcritical load level, but the
latter is based on an essential property of the cylinder: for a 
sequence of cylinders for which $R/t\to \infty$, the nondimensionalized
load-carrying capacity (the highest load at which the unbuckled
state is not only a local but also a global energy minimizer) decreases
to zero. This property was demonstrated implicitly by Hoff \emph{et al.}~\cite{HoffMadsenMayers66},
and Lemma~\ref{th:scaling_we} provides a simplified proof
of this result and a simple sequence of functions that illustrates
the property. 

However, the fact that the mountain-pass solution is localized, and even is the
most localized solution that is possible---a single dimple---is interesting
in its own right, and provides a complementary view of the 
discussion of localization above. A different way of formulating this
result is that ``creating the first dimple is the major obstacle''; 
afterwards one may increase the size of the dimple and add
further dimples without ever returning to the same high energy level.
In itself this interpretation points to a relationship between
single dimples and imperfection sensitivity.

\subsection{Single dimples in other contexts}
Interestingly, single dimples 
have appeared in the literature in a number of seemingly unrelated ways:
\begin{itemize}
\item In the celebrated high-speed camera images of E\ss linger~\cite{Esslngr70}
the first visible deformation is a single, 
well-developed dimple half-way between the ends of the cylinder. New dimples 
quickly appear next to this first dimple, and the deformation then spreads around
the cylinder and in axial direction. It is remarkable, though, that the
first visible deformation is a single dimple.
\item Some of the ``worst'' imperfections calculated by 
Deml and Wunderlich~\cite{DemlWunderlich97} and Wunderlich and Albertin~\cite{WunderlichAlbertin00}
are in the form of a single dimple; as the load decreases, the dimple contracts and becomes
even more concentrated.
\item H\"uhne \emph{et al.}~\cite{HuehneZimmermanRolfesGeier02} 
assert that single dimples are also
\emph{realistic} and \emph{stimulating} imperfections in the sense
of~\cite{WinterstetterSchmidt02}.
\item Zhu \emph{et al.}~\cite{ZhuMandalCalladine02} base their analysis of the scaling behaviour
of the experimental buckling load on the behaviour of a single dimple in other
structural situations (the point-loaded cylinder and the sphere under uniform
external pressure).
\end{itemize}

Note that the single-dimple appearances above are of three different types.
E\ss linger's dimple is an experimental observation; the  dimples
of Wunderlich and co-workers and of H\"uhne \emph{et al.} are geometric imperfection 
profiles; and the dimples studied by Zhu \emph{et al.} are only analogies,
since they are solutions of different problems.


\subsection{Scale-invariance of the localized solutions}

It is an interesting observation that the 
\vKD-equations can be rescaled to depend only on the (rescaled) load level.
For localized solutions, for which the boundary plays no role of
importance, this implies that the set of solutions reduces to a
one-parameter family. This allows for efficient computation of
the behaviour of such solutions but also gives
interesting insight into the relationship between dimples in cylinders of varying
geometry (see \emph{e.g.} \figref{fig:dimple_detail}).

Naturally the scale-invariance is expected to break upon replacing the
\vKD-equations by a different (probably more detailed) shell model.
Nonetheless, it may reasonably be expected that much of the understanding
of the relationship between cylinders of different geometries remains
roughly correct.

We certainly also expect that the large-scale geometry of the 
energy landscape does not depend on the specific model of the cylinder.
Using a discrete mountain-pass algorithm to find mountain-pass points therefore
does not depend on the \vKD-equations, and should give
similar results regardless which shell model is used.

\subsection{Connection with sensitivity to imperfections and ``Perturbation energy''}
Kr\"oplin and co-workers~\cite{KroeplinDinklerHillmann85,DuddeckKroeplinDinklerHillmannWagenhuber89,WagenhuberDuddeck91} were the first to suggest an estimate
of the stability of the unbuckled state in terms of the ratio of a ``perturbation
energy'' (\emph{St\"orenergie}) to the pre-buckling energy. In early 
papers~\cite{KroeplinDinklerHillmann85,DuddeckKroeplinDinklerHillmannWagenhuber89} 
the perturbations are still fixed rather than determined, but 
from both the introduction and the final results in~\cite{WagenhuberDuddeck91} it may
be deduced that
an optimization is done over all perturbations (although
this is simultaneously contradicted on page 333 of this paper). 
Unfortunately, these papers do not provide enough details to determine
exactly what the authors calculate.

\begin{figure}[ht]
\centering
\centerline{\psfig{figure=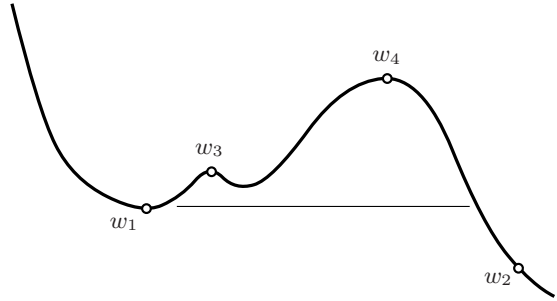,height=4cm}}
\caption{While a local algorithm to find a critical point may settle
on a minor critical point such as $w_3$, 
the mountain-pass algorithm, by its global setup, will
converge on to the essential obstacle $w_4$.
}
\label{fig:severalMP}
\end{figure}

There is one aspect in which our method can clearly be seen to differ from these earlier
approaches. The discrete mountain-pass algorithm takes into account \emph{global} features of the
energy landscape, and provides a global measure of the separation barrier
between two states that lie far apart. This is different from the
papers mentioned above, in which the method uses only local information
(reflected, for instance, in the assumption that the equilibria
in question lie on the same bifurcation branch). 
This difference is illustrated in \figref{fig:severalMP},
where a local analysis might find stationary point $w_2$, but the mountain-pass
algorithm will find the more important obstacle~$w_4$. 

\appendix

\section{Proof of Lemma~\protect\ref{th:scaling_we}}
\label{sec:Yoshimura_proof}

Lemma~\ref{th:scaling_we} states that 
{\it there exists a sequence of functions $w_\d$, $1$-periodic on $\R^2$, 
such that 
\begin{multline}
\label{app:assertion:scaling_we}
\int\limits_{[-1/2,1/2]^2} w_{\d x}^2 \sim 1,\qquad 
\int\limits_{[-1/2,1/2]^2} \Delta w_\d^2 = O(\d^{-1}), \\
\text{and}\qquad
\int\limits _{[-1/2,1/2]^2}\Delta \phi_\d^2 = O(\d^{2-\alpha})
\qquad \text{as $\delta\to0$},
\end{multline}
for any $\alpha>0$.
Here the function $\phi_\d$ solves equation~\pref{eq:main_phi2} with periodic
boundary conditions.
In addition, $w_\d$ and $\phi_\d$ satisfy~\pref{def:BC} on
the boundary of\/ $[-1/2,1/2]^2$.}

The proof consists of three parts. In the first part we construct
the functions $w_\delta$; in the second part we study the symmetry properties
and the support of the right-hand side of~\pref{eq:main_phi2}; and in
the third part we show that this sequence has the 
asserted scaling.

\subsection{Construction of \boldmath $w_\delta$}
Let $f_\d$ be given by
\[
f_\e''(s) = \begin{cases}
  \displaystyle \frac 1{4\d} \qquad &\mathrm{dist}(s,\Z)<\d \\
  0 &\text{otherwise,}
\end{cases}
\qquad\qquad
\text{with}
\qquad
f_\d(0) = f'_\d(0) = 0.
\]
Note that $f$ is even and that $f(1) = 1/4$. 
Define 
\[
w_\d(x,y) = f_\d (y+x) + f_\d(y-x) - \frac12 f_\d(2x) 
   - \frac12 y^2.
\]
We shall drop the subscript $\d$
and simply write $w$ and $f$.

\begin{figure}[ht]
\centering
\centerline{$\vcenter{\psfig{figure=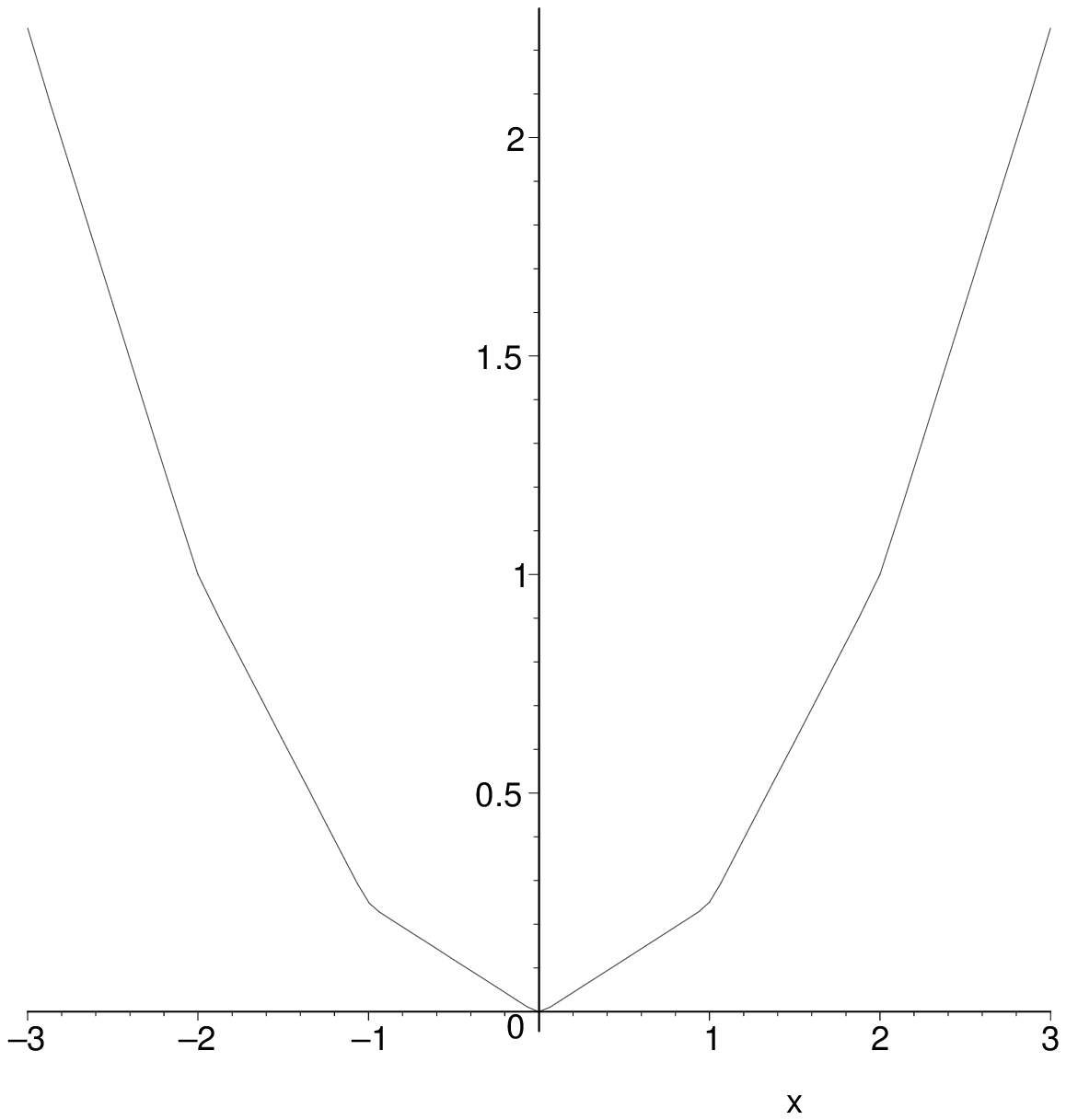,height=3.5cm}}\hskip1.5cm
\vcenter{\psfig{height=4cm,figure=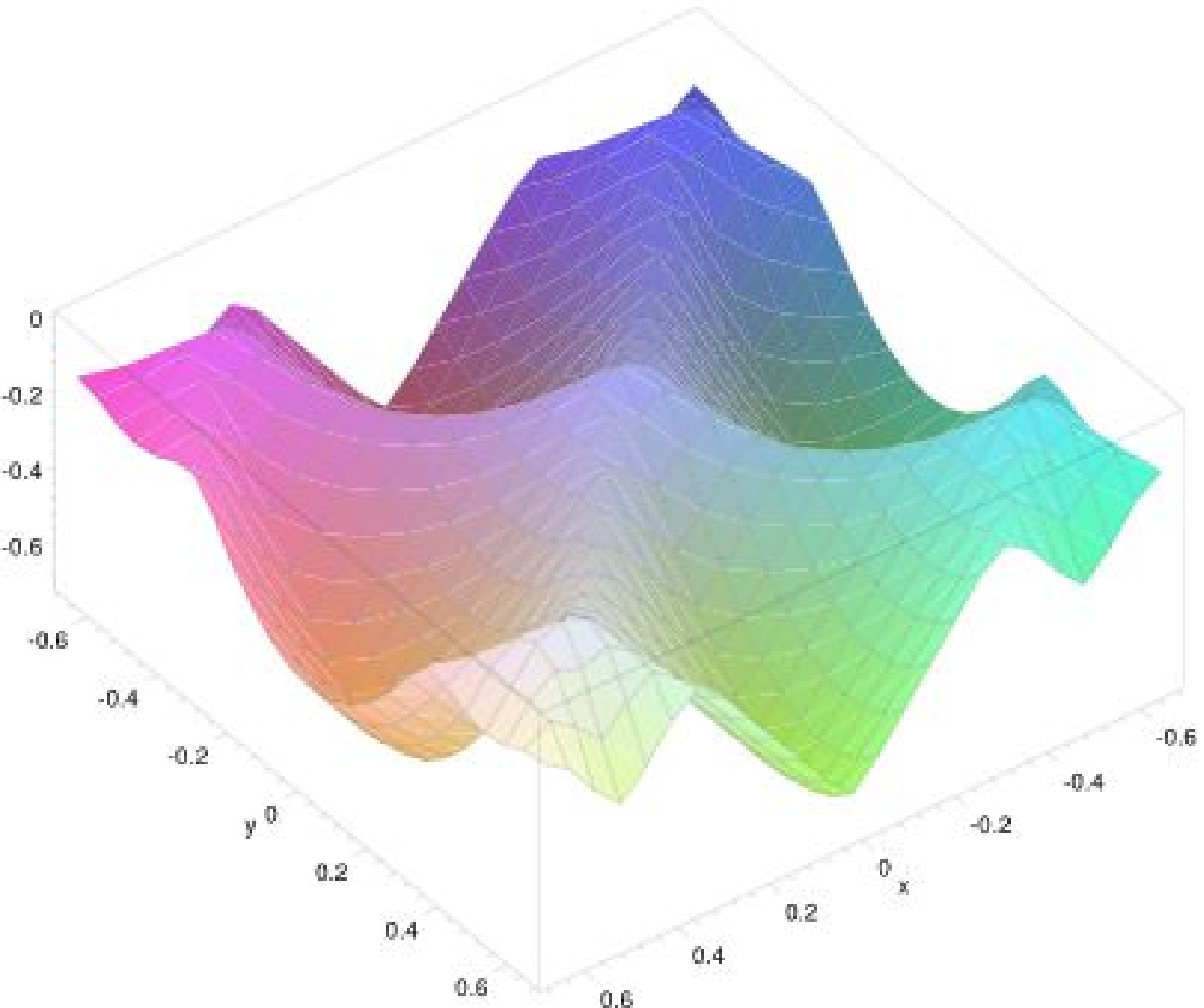}}$}
\caption{The functions $f$ and $-w$; on the right the plotting area is slightly
larger than one period.}
\end{figure}

The function $w$ is periodic on $\R^2$ with period $1$ in each 
direction. To show this, we prove that the first two derivatives
match up on opposite sides of $[-1/2,1/2]\times[-1/2,1/2]$:
\begin{itemize}
\item By the symmetry of $f$ the function $w$ is even in both $x$ and $y$.
Consequently $w$ takes the same values on $(1/2,y)$ and $(-1/2,y)$;
the same holds for $(x,\pm 1/2)$.
\item For the comparison of the first derivatives we calculate
\[
\int_{-1/2}^{1/2} w_{xx}(x,y)\, dx = \int_{-1/2}^{1/2} f''(y+x) \, dx + \int_{-1/2}^{1/2} f''(y-x) \, dx
  - 2\int_{-1/2}^{1/2} f''(2x)\,dx = 0,
\]
implying that $w_x(-1/2,y) = w_x(1/2,y)$; by the symmetry of $w$
it follows that $w_x(-1/2,y) = w_x(1/2,y)=0$. Similarly, 
using the definition of $f_\d$ 
we find that
\[
\int_{-1/2}^{1/2} w_{yy}(x,y)\,dy = \int_{-1/2}^{1/2} f''(y+x) \, dy + \int_{-1/2}^{1/2} f''(y-x) \, dy
  - 1 = 0,
\]
implying that $w_y(x,-1/2)=w_y(x,1/2)=0$.
\end{itemize}
Periodicity on $\R^2$ then follows from the remark that
all second derivatives of $w$ are periodic with period $1$ in $x$ and $y$.

\subsection{Support, symmetry, and boundary conditions}

Next we investigate the right-hand side of~\pref{eq:main_phi2}.
We find 
\begin{align*}
[w,w] + w_{xx} &= \bigl\{ (f''(y+x) + f''(y-x) - 2f''(2x))
  (f''(y+x) + f''(y-x) -1)\\
  &\qquad\qquad  - (f''(y+x)-f''(y-x))^2\bigr\} \\
  &\qquad + (f''(y+x) + f''(y-x) - 2f''(2x)) \\
  &= 4f''(y+x)f''(y-x) - 2f''(2x)(f''(y+x) + f''(y-x)).
\end{align*}
This expression has zero integral over $[-1/2,1/2]^2$. This follows
from the periodicity of $w$:
\begin{equation}
\label{eq:integral_GC_zero}
\int_{-1/2}^{1/2}\int_{-1/2}^{1/2}
      w_{xx}w_{yy}\, dxdy = 
\int_{-1/2}^{1/2}\int_{-1/2}^{1/2} w_{xy}^2\, dxdy
\end{equation}
by partial integration. More is true, however; we analyse the support of 
$[w,w]+w_{xx}$ in $[-1/2,1/2]^2$ in more detail.

The value of $f''$ is either $(4\d)^{-1}$ or zero; in order to
determine  $[w,w]+ w_{xx}$ it is
therefore sufficient to calculate the measures
of the pairwise intersections of the supports of
$f''(y+x)$, $f''(y-x)$, and $f''(2x)$:
\begin{itemize}
\item The intersection of the supports of $f''(y+x)$ and $f''(y-x)$
has total area $4\d^2$ (see \figref{fig:fpfp}).
\begin{figure}[ht]
\centering
\centerline{\psfig{figure=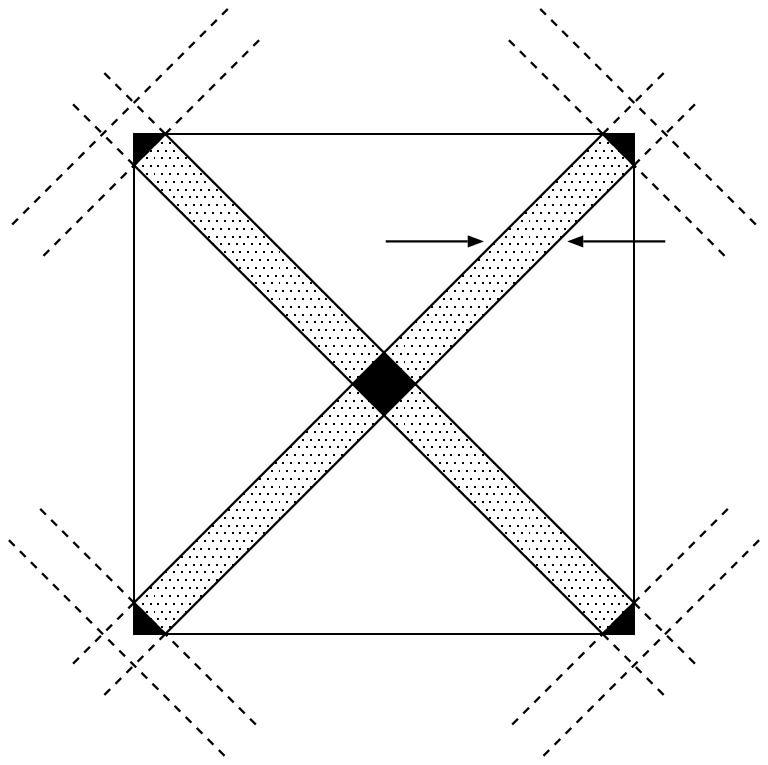,height=5cm}}
\vskip3mm
\caption{The areas of the black regions add up to $4\d^2$}
\label{fig:fpfp}
\end{figure}
\item The intersection of the supports of $f''(y+x)$ and $f''(2x)$
also has total area $4\d^2$ (see \figref{fig:fpf2}).
\begin{figure}[ht]
\centering
\centerline{\psfig{figure=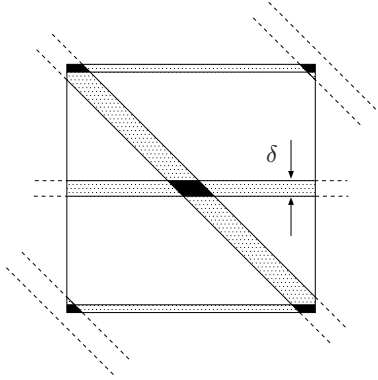,height=5cm}}
\caption{The areas of the black regions add up to $4\d^2$}
\label{fig:fpf2}
\end{figure}
\end{itemize}

Since the support of $[w,w]+ w_{xx}$ concentrates onto
a discrete set of points, let us examine the behaviour
at one of these points. For small $\d$ the support forms disjoint
sets in $[-1/2,1/2]^2$, and we can restrict our attention to
the origin alone.

If $\mod s<1/2$, then $f''_\d(s)$ can be written as
\[
f''_\d(s) = \frac 1\d \, g\left(\frac s\d\right),
\]
where 
\[
g(\sigma) = \begin{cases}
  \displaystyle \frac 14 \qquad &\mod{\sigma}<1 \\
  0 &\text{otherwise.}
\end{cases}
\]
Therefore, as long as $\mod{(x,y)}<1/4$,
\begin{align}
\lefteqn{
4f''_\d(y+x)f''_\d(y-x) - 2f''_\d(2x)(f_\d''(y+x) + f_\d''(y-x)) = }
  \hskip3cm&\notag \\
  {}&= \frac 4{\d^2}\, g\left(\frac{y+x}\d\right) g\left(\frac{y-x}\d\right)
  - \frac2{\d^2}g\left(\frac{2x}\d\right) 
      \left[ g\left(\frac{y+x}\d\right) + g\left(\frac{y-x}\d\right)
  \right] \notag\\
  {}&= \frac1{\d^2}\, F\bigl(\d^{-1}(x,y)\bigr),
  \label{def:F}
\end{align}
where we introduce a new function $F$, which does not depend on
$\d$, to summarize the line above. Note that $\supp F \subset [-2,2]^2$.
Note also that by~\pref{eq:integral_GC_zero}
the function $F$ has zero integral; in addition, since $f''_\d$ is even,
the function $4f''_\d(y+x)f''_\d(y-x) - 2f''_\d(2x)(f_\d''(y+x) + f_\d''(y-x))$
is also even in $x$ and in $y$. Therefore
\[
\int_{\R^2} x F\bigl((x,y)\bigr) \, dxdy = \int_{\R^2} y F\bigl((x,y)\bigr) \, dxdy =0.
\]
This property will be used below.

\medskip
The assertion also states that the functions $w$ and $\phi$ satisfy~\pref{def:BC}
on the boundary of $[-1/2,1/2]^2$. We first note that 
$w$ and $\phi$ are periodic in the following sense:
\begin{equation}
\label{prop:crooked_periodicity}
w(x\pm1/2,y\pm1/2) = w(x,y) \qquad\text{and}\qquad \phi(x\pm1/2,y\pm1/2) = \phi(x,y).
\end{equation}
For $w$ this is a simple consequence of the functional form of $w$; for
$\phi$ it is a consequence of the uniqueness of solutions of~\pref{eq:main_phi2}
under periodic boundary conditions. 
The periodicity of $w$ and $\phi$ in the $y$-direction
in~\pref{def:BC} then follows from a repeated application of~\pref{prop:crooked_periodicity}.
Similarly, the symmetry conditions in $x$ in~\pref{def:BC} follow from
a combination of the symmetry of $w$ and $\phi$ around $\{y=0\}$ in
combination with~\pref{prop:crooked_periodicity}.

\subsection{Scaling properties}

We now use the information gathered above to show that the sequence $w_\d$
has the scaling properties of~\pref{app:assertion:scaling_we}.
All function spaces are on $[-1/2,1/2]^2$.

First, $f'_\d$ remains  bounded on bounded sets as
$\d\to0$; therefore $\int_\Omega w^2_{\d x}$ converges to a finite, positive
value. In addition, all second derivatives of $w_\d$ remain bounded
in $L^1$, so that 
\[
\Mod{\Delta w_\d}_{L^1} \leq C.
\]
The second derivative $f''_\d$ is bounded by $1/4\d$, so that
we can estimate
\[
\Mod{\Delta w_\d}^2_{L^2} \leq 
  \Mod{\Delta w_\d}_{L^1}\Mod{\Delta w_\d}_{L^\infty}
  \leq \frac C\d.
\]

Turning to $\phi_\d$, we start by remarking that $[w_\d,w_\d]+ w_{\d xx}$
is bounded in $L^1$, since
\[
\int_{|(x,y)|<1/4} \bigl|[w_\d,w_\d]+ w_{\d xx} \bigr|
= \frac 1{\d^2}\int_{|(x,y)|<2\d} \bigl|F\bigl(\d^{-1}(x,y)\bigr)\bigr|
= O(1).
\]
Since $W^{2,p} \hookrightarrow L^\infty$ for all $p>1$,
the solution of 
\[
\Delta^2\psi = h
\]
satisfies
\[
\Mod{\Delta\psi}_{L^{p'}} 
= \sup_\zeta \frac{\int_\Omega \Delta \psi\Delta\zeta}{\Mod{\Delta\zeta}_{L^p}}
= \sup_\zeta \frac{\int_\Omega h\zeta}{\Mod{\Delta\zeta}_{L^p}}
\leq C \Mod h_{L^1}\frac{\Mod{\zeta}_{L^\infty}}{\Mod\zeta_{W^{2,p}}}
\leq C \Mod h_{L^1},
\]
so that 
\[
\Mod{\phi_\d}_{W^{2,p'}}
  \leq C \Mod{[w_\d,w_\d]+ w_{\d xx}}_{L^1} \leq C,
\]
where $1/p + 1/p' = 1$. 
Using $W^{2,p'}\hookrightarrow C^{1,1-2/p'}$
we find
\[
\Mod{\phi_\d}_{C^{1,1-2/p'}} \leq C \Mod{\phi_\d}_{W^{2,p'}} \leq C.
\]

Writing, locally at the origin, 
\[
\phi_\d(x,y) = \phi_\d(0,0) + \nabla\phi_\d(0,0) \cdot (x,y) 
  + O\bigl(\mod{(x,y)}^{2(1-1/p')}\bigr),
\]
we find by multiplying the equation~\pref{eq:main_phi2} by $\phi_\d$
and integrating,
\begin{align*}
\Mod{\Delta\phi_\d}^2_{L^2([-1/2,1/2]^2)} 
  &= 2\int_{[-1/4,1/4]^2} \phi_\d \bigl\{[w_\d,w_\d]+w_{\d xx}\bigr\}
  \\
  &= 2\frac{\phi_\d(0,0)}{\d^2} \int_{[-1/4,1/4]^2} F\bigl(\d^{-1}(x,y)\bigr) \\
  &\qquad {} +  2\frac {\nabla\phi_\d(0,0)}{\d^2}
        \cdot  \int_{[-1/4,1/4]^2} (x,y)F\bigl(\d^{-1}(x,y)\bigr) 
      + O\bigl(\d^{2(1-1/p')}\bigr) \\
  &= O\bigl(\d^{2(1-1/p')}\bigr),
\end{align*}
since the zeroth and first moments of $F$ are zero.
Since $p'$ may be chosen arbitrary large, this estimate concludes the proof.

\section{Derivation of the Von K\'arm\'an-Donnell equations}
\label{app:derivation}
\def\EFK{E_{\mathrm{FvK}}}

The common aim of the many elastic shell theories is to approximate
three-dimensional 
elasticity by a reduced description in which the unknowns are functions of
not three but two spatial variables, see for example \cite{BzntCdlni}.
\drop{

\meta{Comment on what this appendix is about: showing what the 
boundary conditions on $\phi$ should be}

The derivation of equations~(\ref{eq:main_w}--\ref{eq:main_phi}) given
by Donnell~\cite{Dnnll34} and Von K\'arm\'an and Tsien~\cite{vKT} 
follows from an 

From the three-dimensional elastic theory a two-dimensional shell theory is
derived by 
}%
For the \vKD\ cylinder the central approximation is 
the \emph{director Ansatz}, which states that 
a normal to the centre surface remains normal through deformation. 
By this Ansatz the displacement is fully characterized by
the displacement vector $(u,v,w)$, a function of the two in-plane 
spatial variables $x$ and $y$, where $u$ is the 
displacement in axial ($x$-), $v$ in tangential ($y$-),
and $w$ in the radial direction. 
Apart from some rescaling
the function $w$ is the same as the unknown $w$ in the rest of this paper.

In the formulation of Section~\ref{sec:vKD} the unknowns $u$ and $v$ are replaced
by the Airy stress function~$\phi$, 
which is derived
by minimization
with respect to the displacements $u$ and $v$ for fixed $w$.
This minimization argument is well known in the context
of the Von K\'arm\'an plate theory and can be found in many textbooks. 
Determining the boundary conditions that the function $\phi$ satisfies, however, is
not straightforward (see also the discussion in~\cite{SchaefferGolubitsky79}),
and it is for this reason that we now describe the argument in detail.
The main goal is to show that  the function $\phi$ is periodic in the tangential direction.

\medskip

\subsection{Energy and shortening}
\let\s\sigma
All quantities in this appendix are dimensional. We assume a cylinder of thickness
$t$, length $L$,  and
radius $R$, and we set $\Omega = [0,L]\times [0,2\pi R]$. The stored energy 
given by~\cite{vKT} is 
\[
E_1=
\frac{t^3E}{24(1-\nu^2)} \int_\Omega \Delta w^2 
+ \frac t{2E}\int_\Omega \bigl[(\s_{11} + \s_{22})^2 - 2(1+\nu)(\s_{11}\s_{22}-\s_{12}^2)\bigr]
\]
for a linear material of Young's modulus $E$ and Possion's ratio $\nu$. Under
the assumption of plane stress the stress and strain tensors are related by
\begin{equation}
\label{eq:Hooke}
\s = \frac E{1-\nu^2} \bigl[(1-\nu) \e + \nu \tr \e \, I\bigr] = \frac E{1-\nu^2}\begin{pmatrix}
\e_{11} + \nu \e_{22} & (1-\nu)\e_{12} \\[\jot]
(1-\nu)\e_{12} & \e_{22} + \nu \e_{11}
\end{pmatrix},
\end{equation}
and under the small-angle approximation the strain tensor can be expressed in 
the displacements as
\begin{equation}
\label{eq:straintensor}
\e = \begin{pmatrix}
u_x +\tfrac12 w_x^2 & \tfrac12 u_y + \tfrac 12 v_x + \tfrac12 w_xw_y\\[\jot]
\tfrac12 u_y + \tfrac 12 v_x + \tfrac12 w_xw_y & v_y + \tfrac12 w_y^2 -\rho w
\end{pmatrix}.
\end{equation}
These choices for the energy and for the  stress and strain tensors are very similar
to those for a  flat plate. 
The intrinsic curvature of the cylinder, of magnitude $\rho=1/R$, only
appears in the
last term of $\e_{22}$, $-\rho w$, which expresses the fact that 
radial displacement creates extensional strain in the $y$-direction. 

The average axial shortening is given by 
\[
S_1 = -\frac1{2\pi R} \int_0^{2\pi R}\bigl[ u(L,y) - u(0,y)\bigr] \, dy 
= - \frac 1{2\pi R}\int_\Omega u_x \, dxdy,
\]
and an equilibrium $(u,v,w)$ at load level $P$ is a stationary point of 
the total potential $V_1 = E_1 - P S_1$.

\subsection{Boundary conditions}

At the boundaries $y=0,2\pi R$ it is natural
to assume that $u$, $v$, and $w$ are periodic, 
but at $x=0,L$ there is a certain amount
of choice. 

The boundary conditions on $w$  (see~\pref{def:BCw}) are
\begin{equation}
\label{def:BC'}
w_x = (\Delta w)_x  = 0, \qquad \text{at }x=0,L,
\end{equation}
and these conditions signify a fixed angle ($w_x=0$) and zero radial force ($(\Delta w)_x=0$). 
They may also be understood as symmetry boundary conditions, as in the case
of a sequence of cylinders stacked on top of each other
.
For $u$ and $v$ we assume boundary conditions
\begin{equation}
\label{def:BCuv}
u_y = 0 \qquad\text{and}\qquad \s_{12} = 0,\qquad \text{at }x=0,L,
\end{equation}
which signify that the ends of the cylinder are rigid in the $x$-direction and 
that there is no friction between the cylinder and the apparatus holding it. Note
that the pair of boundary conditions $\s_{12}=0$ and $(\Delta w)_x=0$ 
together states that the loading apparatus exerts only axial forces on the cylinder.

\medskip

The boundary conditions on $w$ are invariant under addition of a constant to $w$,
\emph{i.e.}\ under replacing $w$ by $w+c$; for stationary points we may exploit this
fact.
\begin{lemma}
If $(u,v,w)$ is a stationary point of $E_1-PS_1$ under  
boundary conditions~(\ref{def:BC'}-\ref{def:BCuv}), then
\begin{equation}
\label{eq:s22}
\int_\Omega \s_{22} = 0.
\end{equation}
\end{lemma}

\begin{proof}
Under the replacement $w\mapsto w+c$, we have
\[
\frac { d\s}{dc} = -\frac {E\rho}{1-\nu^2}\begin{pmatrix}
\nu & 0 \\ 0 & 1 \end{pmatrix},
\]
and therefore
\[
0 = \frac d{dc} (E_1-PS_1) = 
-\frac {t\rho}{1-\nu^2}\int_\Omega \bigl[ (\s_{11}+\s_{22})(\nu+1) - (1+\nu)(\s_{11}+\nu\s_{22})\bigr]
= -t\rho\int_\Omega \s_{22}.
\]
\end{proof}

\subsection{Derivation of the Airy stress function \boldmath $\phi$}

The energy~\pref{def:energy} and the Airy stress function $\phi$ are derived from
the total potential $E_1-PS_1$ by minimization
with respect to the displacements $u$ and $v$ for fixed $w$.
\drop{
This minimization argument is well known in the context
of the Von K\'arm\'an plate theory and can be found in many textbooks. 
Determining the boundary conditions that $\phi$ satisfies, however, is
not straightforward (see also the discussion in~\cite{SchaefferGolubitsky79}),
and it is for this reason that we now describe the argument in detail.
The main goal is to show that  the function $\phi$ is periodic in the tangential direction.

\medskip

From minimizing the second term in $E_1$ with respect to $u$ and $v$ while
fixing $w$ 
}%
Performing this minimization on
the second term in $E_1$ yields the classical plate equilibrium equations
\[
\s_{11x} + \s_{12y} = 0 \qquad\text{and}\qquad  \s_{12x}+\s_{22y} = 0.
\]
Note that the derivative of $S_1$ with respect to $u$ only creates
boundary terms. 
By applying three times the well-known characterization of divergence-free
vector fields as rotations of scalar fields (see \emph{e.g.}~\cite[Th.~XII.3.5]{Antman}) 
we obtain the \emph{local} existence of
a function $\phi$ satisfying
\begin{equation}
\label{def:phi}
\s_{11} = E\phi_{yy}, \qquad \s_{12} = -E\phi_{xy}, \qquad \text{and} \qquad
\s_{22} = E\phi_{xx},
\end{equation}
where we use the traditional scaling of $\phi$ by Young's modulus.

\subsection{Boundary conditions on \boldmath $\phi$}

The existence of the function $\phi$ is the result of a local differential-geometric
argument, and as such gives no reason for $\phi$ to be periodic in $y$. 
The following theorem shows that after a normalization transformation the
function $\phi$ can indeed be assumed to be periodic in $y$, and may be taken to 
satisfy the same boundary conditions as the function $w$. 

\begin{theorem}
\label{th:phi}
If $u$, $v$, and $w$ are periodic in $y$ and satisfy boundary
conditions~(\ref{def:BC'}-\ref{def:BCuv}), then there exists a function $\phi$
that satisfies
\begin{equation}
\label{def:phi2}
\s_{11} - \frac 1{\mod \Omega} \int_\Omega \s_{11} = E\phi_{yy}, 
\qquad \s_{12} = -E\phi_{xy}, \qquad \text{and} \qquad
\s_{22} = E\phi_{xx},
\end{equation}
is periodic in $y$,  and satisfies boundary conditions 
\begin{equation}
\label{app:def:BCphi}
\phi_x = (\Delta \phi)_x = 0 \qquad\text{at }x=0,L.
\end{equation}
\end{theorem}

\begin{remark}
Mechanically the normalization of $\phi$ with $\int_\Omega \sigma_{11}$ 
means that $\phi$ 
represents the \emph{deviation} from the unbuckled in-plane stress state.
\end{remark}

\begin{proof}
As discussed above, there
exists a function $\phi$ satisfying~\pref{def:phi}; we will construct
in stages a new function $\hat \phi$ which satisfies~\pref{def:phi2} 
and the boundary conditions.

We first convert condition~\pref{def:phi} into~\pref{def:phi2}.
Set
\[
p(x) := \frac 1{2\pi R}\int_0^{2\pi R} \phi_{yy}(x,y)\, dy 
  = \frac1{2\pi R} \bigl[\phi_y(x,2\pi R)-\phi_y(x,0)\bigr].
\]
Since the second derivatives
of $\phi$ can be expressed in terms of derivatives of $u$, $v$, and $w$, the second 
and higher derivatives of $\phi$ are automatically periodic in $y$. Therefore
\[
\frac d{dx} p(x) =\frac1{2\pi R}\bigl[ \phi_{xy}(x,2\pi R) - \phi_{xy}(x,0)\bigr] = 0,
\]
implying that $p$ is actually independent of $x$. (A mechanical argument
provides the same result: $Etp(x) = t\dashint_0^{2\pi R}\sigma_{11}(x,y)dy$ is
the total force applied at a virtual cut at level $x$, and mechanical equilibrium 
implies that this force is independent of $x$). Therefore
\[
\mod\Omega p = 2\pi R \int_0^L p\, dx = \int_\Omega \phi_{yy} 
= \frac 1E\int_\Omega \s_{11},
\]
so that the new function 
\[
\tilde \phi(x,y) := \phi(x,y) - \frac p{2} y^2
\]
satisfies~\pref{def:phi2}. Note that this implies
\begin{equation}
\label{eq:intphiyy0}
\int_\Omega \tilde \phi_{yy} = 0.
\end{equation}

We now turn to the periodicity in the $y$-direction. It remains
to show that $\tilde\phi$, $\tilde\phi_x$, and 
$\tilde\phi_y$ are the same at $y=0$ and $y=2\pi R$. 
Again the periodicity of the second derivatives implies that 
\[
\frac{d^2}{dx^2} \bigl[ \tilde\phi(x,2\pi R) - \tilde\phi(x,0)\bigr]
= \tilde\phi_{xx}(x,2\pi R) - \tilde\phi_{xx}(x,0) = 0,
\]
so that $\tilde\phi(x,2\pi R)-\tilde\phi(x,0) = ax+b$ for some  $a,b\in\R$. Defining
\[
\hat\phi(x,y) := \tilde\phi(x,y) - \frac {by}{2\pi R} 
= \phi(x,y) - \frac p2 y^2 - \frac {by}{2\pi R}
\]
the function $\hat\phi$ still satisfies~\pref{def:phi2}, and
\[
\hat\phi(x,2\pi R) - \hat\phi(x,0) = ax. 
\]
Finally, we find that
\[
a = \hat\phi_x(0,2\pi R) - \hat\phi_x(0,0) = \int_0^{2\pi R} \hat\phi_{xy}(0,y)\, dy 
= \frac1E \int_0^{2\pi R} \sigma_{12}(0,y)\, dy 
\stackrel{\pref{def:BCuv}}= 0,
\]
and therefore that $\hat\phi(x,2\pi R) - \hat\phi(x,0) = 0$ for all $x$. The same
follows for $\hat\phi_x(x,2\pi R) - \hat\phi_x(x,0)$ by differentiation.

To show that $\phi_y$ also matches, 
\[
\frac d{dx}\bigl[\hat\phi_y(x,2\pi R) - \hat\phi_y(x,0) \bigr]
= \bigl[\phi_{xy}(x,2\pi R)-\phi_{xy}(x,0)\bigr] = 0,
\]
and therefore $\phi_y(x,2\pi R) - \phi_y(x,0)$ is constant in $x$;
by~\pref{eq:intphiyy0} this constant is zero. This proves that $\hat \phi$ 
satisfies~\pref{def:phi2} and is periodic in $y$.

We finally discuss the boundary conditions at $x=0,L$, and we follow the line
of reasoning of~\cite{SchaefferGolubitsky79}. 
By~\pref{def:BCuv} and~\pref{eq:Hooke} $\e_{12}=0$ at $x=0,L$, so that
by~\pref{def:BC'} and~\pref{def:BCuv}
\[
v_{xy} = \frac \partial {\partial y} (2\e_{12} - u_y - w_xw_y) = 0\qquad\text{at }x=0,L.
\]
Therefore
\[
\e_{22x}  = v_{xy} + w_yw_{xy} - \rho w_x \stackrel{\pref{def:BCuv}}= 0\qquad\text{at }x=0,L.
\]
Using $E\e_{22} = \s_{22} - \nu \s_{11}$, we then find
\[
\hat \phi_{xxx} - \nu \hat \phi_{xyy} 
= \phi_{xxx} - \nu \phi_{xyy} = \frac1E\, \frac d{dx} (\s_{22} - \nu \s_{11}) 
= \e_{22x}= 0,
\]
and by adding  $(1+\nu) \hat\phi_{xyy} = -(1/E)(1+\nu) \s_{12y}=0$ it follows that
\[
(\Delta\hat\phi)_x  = \hat\phi_{xxx} +\hat\phi_{xyy} = 0\qquad\text{at }x=0,L,
\]
which proves one part of~\pref{app:def:BCphi}.

From $\hat\phi_{xy} = -\s_{12}/E = 0$ we find that
\[
\hat\phi_x(0,y) = c_0 \qquad \text{and}\qquad  
\hat\phi_x(L,y) =c_L \qquad \text{for all }y\in [0,2\pi R].
\]
Writing
\[
2\pi R(c_L-c_0) = \int_\Omega \hat\phi_{xx} 
= \frac1E\int_\Omega \s_{22} \stackrel{\pref{eq:s22}}= 0
\]
we find that $c_L=c_0$. Now the function 
\[
\overline \phi(x,y) := \hat\phi(x,y) -c_0 x = \phi(x,y) - \frac p2 y^2 - \frac {by}{2\pi R}
-c_0x
\]
satisfies~\pref{def:phi2} and~\pref{app:def:BCphi} and is periodic in $y$.
This concludes the proof.
\end{proof}

\begin{remark}
It is instructive to note that the periodicity of $\phi$ is a result
of the specific choice of boundary conditions, 
and will in fact not hold if different boundary conditions
are taken. For instance, if a tangential shear stress $\tau$ is applied at the cylinder
ends (\emph{i.e.} the cylinder is loaded under torsion), then the coefficient $a$
in the derivation above will not vanish, and $\phi_x$ will not be periodic in
$y$.
\end{remark}

\subsection{Putting it all together}

By an elementary but lengthy calculation we find that  $\phi$,
as provided by Theorem~\ref{th:phi}, satisfies the 
equation
\begin{equation}
\label{appeq:phi}
\Delta^2 \phi + \rho w_{xx} + [w,w] = 0 \qquad\text{in } \Omega,
\end{equation}
and that the second term in $E_1$ can be written as 
\[
\frac {tE}2 \int_\Omega \Bigl[ \Delta\phi^2 - 2(1+\nu)[\phi,\phi]\Bigr].
\]
By the boundary conditions given by Theorem~\ref{th:phi} the second
term vanishes, and the total stored energy functional can therefore
be written as
\[
E_2(w) := \frac{t^3E}{24(1-\nu^2)} \int_\Omega \Delta w^2 
 + \frac {tE}2 \int_\Omega  \Delta\phi^2 .
\]
Note that this energy is a function of $w$ alone; the function $\phi$ in this
definition is assumed to be given by~\pref{appeq:phi}, with
the boundary conditions of Theorem~\ref{th:phi}. 
Similarly, we rewrite the average shortening as
\begin{eqnarray*}
S_2(w) \;:=\; S_1(u) &=& - \frac 1{2\pi R}\int_\Omega u_x \\
&\stackrel{\pref{eq:straintensor}}=& 
  -\frac 1{2\pi R}\int_\Omega \bigl[ \e_{11} - \tfrac12 w_x^2\bigr] \\
&\stackrel{\pref{eq:Hooke}}=& 
  -\frac1{2\pi RE} \int_\Omega \bigl[\sigma_{11}-\nu\sigma_{22}\bigr]
  + \frac1{4\pi R} \int_\Omega w^2_x\\
&\stackrel{\pref{def:phi2}}=& 
  -\frac1{2\pi R}\int_\Omega\bigl[\phi_{yy} - \nu\phi_{xx}\bigr]
  + \frac1{4\pi R} \int_\Omega w_x^2 \\
&=& \frac1{4\pi R}\int_\Omega w^2_x.
\end{eqnarray*}

A stationary point $w$ of $E_2-PS_2$ satisfies the Euler equation
\[
\frac{t^2}{12(1-\nu^2)} \Delta^2 w + \frac P{2\pi R E t} w_{xx} 
  - \rho \phi_{xx} - 2[w,\phi] = 0,
\]
where again $\phi$ is related to $w$ by~\pref{appeq:phi}.
With the nondimensionalization 
\[
w = 4\pi^2R \, \bw, \qquad
\phi = 16\pi^4 R^2\, \bphi, \qquad
x \mapsto 2\pi R\, x,\qquad
y \mapsto 2\pi R\, y,
\]
we then obtain equations~\pref{eq:main_w2} and~\pref{eq:main_phi2}, 
and the dimensional energy $E_2$ and average shortening $S_2$ above can be expressed
in these variables as
\begin{equation}
\label{def:rescaled_energies}
E_2 = \frac{\pi^2t^3E}{6(1-\nu^2)} \int \Delta \bw^2 
 + 32\pi^6tER^2\int  \Delta\bphi^2 , \qquad
S_2 = 4\pi^3 R \int \bw^2_x.
\end{equation}

\bibliography{refs}
\end{document}